\documentstyle[12pt,amssymb,amsbsy,righttag]{amsart} 
\renewcommand{\a}{\alpha}
\renewcommand{\b}{\beta}
\newcommand{\g}{\gamma}
\renewcommand{\d}{\delta}
\newcommand{\e}{\varepsilon}
\newcommand{\z}{\zeta}

\renewcommand{\l}{\lambda}
\newcommand{\r}{\rho}

\renewcommand{\t}{\tau}

\newcommand{\f}{\varphi}
\renewcommand{\o}{\omega}

\newcommand{\G}{\Gamma}
\newcommand{\D}{\Delta}

\renewcommand{\O}{\Omega}

\newcommand{\B}{{\cal B}}

\newcommand{\cd}{{\cal D}}
\newcommand{\F}{{\cal F}}
\newcommand{\h}{{\cal H}}

\newcommand{\K}{{\cal K}}
\newcommand{\cL}{{\cal L}}
\newcommand{\M}{{\cal M}}
\newcommand{\p}{{\cal P}}
\newcommand{\cp}{{\cal P}}

\newcommand{\X}{{\cal X}}

\newcommand{\C}{{\Bbb C}}
\newcommand{\T}{{\Bbb T}}

\newcommand{\dd}{{\Bbb D}}
\newcommand{\R}{{\Bbb R}}
\newcommand{\Z}{{\Bbb Z}}

\newcommand{\0}{{\Bbb O}}

\newcommand{\bs}{\boldsymbol}

\newcommand{\m}{{\boldsymbol m}}
\newcommand{\bS}{{\boldsymbol S}}

\newcommand{\rf}[1]{(\ref{#1})}

\newcommand{\df}{\stackrel{\mathrm{def}}{=}}

\newcommand{\re}{\operatorname{Re}}

\newcommand{\supp}{\operatorname{supp}}

\newcommand{\trace}{\operatorname{trace}}
\newcommand{\rank}{\operatorname{rank}}
\newcommand{\const}{\operatorname{const}}

\newcommand{\eeq}{\end{equation}}
\newcommand{\beq}{\begin{equation}}
\newcommand{\bay}{\begin{eqnarray}}
\newcommand{\ey}{\end{eqnarray}}
\newcommand{\bey}{\begin{eqnarray*}}
\newcommand{\eey}{\end{eqnarray*}}

\newcommand{\eq}{\Leftrightarrow}
\newcommand{\imp}{\Rightarrow}
\newcommand{\be}{\infty}

\newcommand{\bl}{\blacksquare}

\newcommand{\Pf}{{\bf Proof. }}
\newcommand{\im}{\operatorname{Im}}
\renewcommand{\re}{\operatorname{Re}}
\newcommand{\ov}{\overline}

\newtheorem{thm}{\hspace{\parindent}Theorem}[section]

\newtheorem{cor}[thm]{\hspace{\parindent}Corollary}
\newtheorem{lem}[thm]{\hspace{\parindent}Lemma}

\newenvironment{remark}{\medbreak{\bf Remark.}}{\medbreak}
\newenvironment{example}{\medbreak{\bf Example.}}{\medbreak}
\newenvironment{definition}{\medbreak{\bf Definition.}}{\medbreak}
\newenvironment{proof}{\Pf}{\qed}
\newenvironment{proofx}{\Pf}{}

\begin{document}
\renewcommand{\theequation}{\thesection.\arabic{equation}}

\newcommand\bsQ{{\mathbf{Q}}}
\newcommand\bsX{{\frak{X}}}
\newcommand\bsh{{\mathbf{h}}}
\newcommand\bsf{{\phi}}

\newcommand\fa{\f_{[a]}}
\newcommand\vvp{~} 

\newcommand{\Mp}{{\frak M}_p}
\newcommand\inner[1]{\langle#1\rangle}
\newcommand\fzh{f_{z,h}}
\newcommand\fzhi{f_{{z_1},h}}
\newcommand\fzhii{f_{{z_2},h}}
\newcommand\qmu{{\bsQ}^+_\mu}   
\newcommand\qq{Q_}
\newcommand\OI{[0,1]}
\newcommand\set[1]{{\{#1\}}}
\newcommand\ggl{g_\l}
\newcommand\Ro{R_\o}
\newcommand\gTo{\Theta_\o}
\newcommand\gTn{\Theta_\nu}
\newcommand\ntoo{{n\to\be}}
\newcommand\norm[1]{\|#1\|}

\newcommand\limxoo{\lim_{x\to\infty}}

\newcommand\xtoo{{x\to\infty}}

\newcommand\CSineq{Cauchy--Schwarz inequality}

\newcommand\into{\int_{0}^{\infty}}
\newcommand{\Leb}{\operatorname{Leb}}

\newcommand\iid{i.i.d.\spacefactor=1000}     
\newcommand\ie{i.e.\spacefactor=1000}
\newcommand\eg{e.g.\spacefactor=1000}
\newcommand\viz{{viz.}\spacefactor=1000}
\newcommand\cf{{cf.}\spacefactor=1000}
\newcommand\aex{a.e.\spacefactor=1000}

\newcommand\ett{\boldsymbol1}
\newcommand\wS[1]{\bS_{#1,\infty}}
\newcommand\ettx[1]{\ett[#1]}

\renewcommand{\theenumi}{\textup{(\roman{enumi})}} 
\renewcommand{\labelenumi}{\theenumi}

\newcommand\Lloc{L_{\textup{loc}}} 
\newcommand\qf{Q_\f}
\newcommand\rmi{{\mathrm{i}}}
\newcommand\xfrac[2]{#1/#2}
\newcommand\parfrac[2]{\Bigl(\frac{#1}{#2}\Bigr)}
\newcommand\bigpar[1]{\bigl(#1\bigr)}
\newcommand\Bigpar[1]{\Bigl(#1\Bigr)}
\newcommand\biggpar[1]{\biggl(#1\biggr)}
\newcommand\Biggpar[1]{\Biggl(#1\Biggr)}
\def\rompar(#1){\textup(#1\textup)}    
\newcommand\floor[1]{[#1]}

\newcommand\ga{\alpha}
\newcommand\gb{\beta}
\newcommand\gd{\delta}
\newcommand\gf{\varphi}
\newcommand\gG{\Gamma}
\newcommand\gl{\lambda}
\newcommand\go{\omega}
\newcommand\gT{\Theta}
\newcommand\eps{\varepsilon}
\newcommand\tk{\tilde k}
\newcommand\ii{{\mathrm i}}

\newcommand\upto{\uparrow}
\renewcommand\qed{$\bl$}
\newcommand\nopf{$\bl$}

\newcommand\oxx[2]{\go_{#1}^{(#2)}}  
\newcommand\oxz[1]{\oxx{#1}{\be}}
\newcommand\oxi[1]{\oxx{#1}{1}} 
\newcommand\oxii[1]{\oxx{#1}{2}} 
\newcommand\oxp[1]{\oxx{#1}{p}} 
\newcommand\oyx[2]{\widetilde\go_{#1}^{(#2)}}  
\newcommand\oyp[1]{\oyx{#1}{p}} 
\newcommand\oyii[1]{\oyx{#1}{2}} 

\newcommand\oxzf{\oxz{\f}}  
\newcommand\oxpf{\oxp{\f}}  
\newcommand\oxiif{\oxii{\f}}  
\newcommand\oypf{\oyp{\f}}  

\newcommand\Holder{H\"older}

\newcommand{\refS}[1]{\S \ref{#1}} 
\newcommand{\refT}[1]{Theorem~\ref{#1}}
\newcommand{\refC}[1]{Corollary~\ref{#1}}
\newcommand{\refL}[1]{Lemma~\ref{#1}}
\newcommand{\refP}[1]{Proposition~\ref{#1}}
\newcommand{\refand}[2]{\ref{#1} and~\ref{#2}}

\newcommand\pfitem[2]{\emph{Step #1. #2.}}

\author{A.B. Aleksandrov, S. Janson, V.V. Peller and R. Rochberg}

\address{Alexei B. Aleksandrov,
St-Petersburg Branch,  
Steklov Institute of 
Mathematics, Fontanka 27, 
191011 St-Petersburg, 
Russia} 
\email{alex@@pdmi.ras.ru}

\address{Svante Janson, 
Department of Mathematics, 
Uppsala University, 
751 06 Uppsala,
Sweden}
\email{svante@@math.uu.se}

\address{Vladimir V. Peller,
Department of Mathematics,
Kansas State University,
Manhattan, Kansas 66506,
USA}
\email{peller@@math.ksu.edu}

\address{Richard Rochberg,
Department of Mathematics,
Washington University,
St.\ Louis, MO 63130,
USA}
\email{rr@@math.wustl.edu}

\thanks{The first author is partially supported by Grant 99-01-00103
of Russian Foundation of Fundamental Studies and by Grant 326.53 of
Integration. 
The third author is partially supported by NSF grant DMS 9970561. 
The fourth author is partially supported by NSF grant DMS 9970366}

\title{An Interesting Class of Operators with unusual Schatten--von
Neumann behavior}

\dedicatory{We dedicate this paper to Jaak Peetre on occasion of his
65th birthday and to the memory of
Tom Wolff.
Both helped shape the mathematics of our time and profoundly
influenced our mathematical thoughts.
Each, through his singular humanity, helped our hearts grow.}

\begin{abstract} 
We consider the class of integral operators $Q_\f$ on $L^2(\R_+)$
of the form $(Q_\f f)(x)=\int_0^\be\f (\max\{  x,y\})f(y)dy$.
We discuss necessary
and sufficient conditions on $\varphi$ to insure that $Q_{\varphi}$ is
bounded, compact, 
or in the Schatten--von Neumann class $\bS_p$,
$1<p<\infty$. 
We
also give necessary and sufficient conditions for $Q_{\varphi}$ to be a finite
rank operator. 
However, there is a kind of cut-off at $p=1$, and 
for membership in $\bS_{p}$, $0<p\leq1$, the
situation is more complicated.
Although we give 
various necessary conditions and sufficient conditions 
relating to $Q_{\varphi}\in\bS_{p}$ in that range, we do not have
necessary and sufficient conditions. 
In the most important case $p=1$, 
we have a necessary 
condition and a sufficient condition, using $L^1$ and $L^2$ modulus of
continuity, respectively, with a rather small gap in between.
A second cut-off occurs at $p=1/2$: if $\f$ is sufficiently smooth
and decays reasonably fast, then $\qf$ belongs to the weak 
Schatten--von Neumann class $\wS{1/2}$, but never to 
$\bS_{1/2}$ unless $\f=0$.

We also obtain results
for related families of operators acting on $L^2(\R)$ and $\ell^2(\Z)$.

We further study operations acting on bounded linear operators on 
$L^{2}(\R^{+})$ 
related to the class of operators $Q_\f$.
In particular we
study Schur multipliers given by functions of the form $\varphi\left(
\max\left\{  x,y\right\}\right) $ and 
we study properties of the averaging projection (Hilbert--Schmidt projection)
onto the operators of the form $Q_\f$.
\end{abstract}

\maketitle

\

\begin{center}
{\Large Contents}
\end{center}


\begin{enumerate}
\item[1.] Introduction\dotfill \pageref{intr}
\item[2.] Preliminaries\dotfill \pageref{prel}
\item[3.] Boundedness, compactness, and $p>1$\dotfill \pageref{S:p>1}
\item[4.] Positive operators\dotfill \pageref{S:positive}
\item[5.] A sufficient condition, $1/2<p\le1$\dotfill \pageref{S:yp}
\item[6.] $p=1$, first results\dotfill \pageref{S:p1}
\item[7.] Schur multipliers of the form $\psi(\max\{x,y\})$,
$x,y\in\R_+$\dotfill \pageref{Schur} 
\item[8.] The case $p=1/2$\dotfill \pageref{S:1/2}
\item[9.] Sturm--Liouville theory and $p=1/2$\dotfill \pageref{S:Sturm}
\item[10.] More on $p=1$\dotfill \pageref{more}
\item[11.] Averaging projection\dotfill \pageref{aver}
\item[12.] Finite rank\dotfill \pageref{S:finite}
\item[13.] A class of integral operators on $L^2(\R)$\dotfill \pageref{io}
\item[14.] Matrix representation\dotfill \pageref{S:matrix}
\item[15.] Necessary conditions for $Q_\f\in\bS_1$\dotfill \pageref{neces}
\item[16.] Dilation of Symbols\dotfill \pageref{S:dilation}
\item[] Acknowledgement\dotfill \pageref{ack}
\item[] References\dotfill \pageref{refer}
\end{enumerate}

\

\setcounter{equation}{0}
\section{\bf Introduction}
\label{intr}

\

For a function $\f\in \Lloc^1(\R_+)$, 
which means that $\f$
is a locally integrable function on $\R_+=(0,\be)$,
we define the operator
$Q_\f$ on the set of bounded compactly supported functions $f$ in
$L^2(\R_+)$ by
\begin{equation}\label{q}
(Q_\f f)(x)=\int_0^\be\f\left(\max\{x,y\}\right)f(y)dy;
\end{equation}
equivalently,
\begin{equation}\label{q2}
(Q_\f f)(x)=\f(x)\int_0^xf(y)dy+\int_x^\be\f(y)f(y)dy.
\end{equation}
We are going to study when  $Q_\f$ is (\ie, extends to) a bounded
operator in $L^2(\R_+)$, and when this operator is
compact, or belongs to Schatten--von Neumann classes $\bS_p$.

We will also consider the corresponding Volterra operators
$\qf^+$ and $\qf^-$ defined by
\begin{align*}
(\qf^+ f)(x)&=\f(x)\int_0^xf(y)dy\\
(\qf^- f)(x)&=\int_x^\be\f(y)f(y)dy; 
\end{align*}
thus $\qf=\qf^++\qf^-$.

It is straightforward 
to see (and proved more generally in \refT{ker+}) 
that if any
of these three operators is bounded on $L^2(\R_+)$, then
$\int_a^\be|\f|^2<\be$ for any $a>0$ and thus the integrals above
converge, and define the operators, for any $f\in L^2(\R_+)$ and every $x>0$.
				
We find in \refS{S:p>1} simple necessary and sufficient conditions for
$\qf$ to be bounded or compact, and for $\qf\in\bS_p$, $1<p<\be$.
The conditions are $\f\in X_\be$, $\f\in X_\be^0$ and $\f\in X_p$,
respectively, where 
the spaces $X_p$, $X_\be$, and $X_\be^0$ are defined as follows.

\begin{definition}
If\/ $0<p<\infty$, 
let $X_p$ be the linear space of all measurable functions on $\R_+$
that satisfy the equivalent conditions
\begin{gather}
\sum_{n\in\Z}2^{np/2}
 \left(\int_{2^n}^{2^{n+1}}|\f(x)|^2dx\right)^{p/2}<\be;
\label{xpi}\\
\sum_{n\in\Z}2^{np/2}
 \left(\int_{2^n}^{\be}|\f(x)|^2dx\right)^{p/2}<\be;
\label{xpii}\\
 x^{1/2}\left(\int_{x}^\be|\f(y)|^2dy\right)^{1/2} \in L^p(dx/x).
\label{xpiii}
\end{gather}
Similarly, let $X_\be$ 
be the linear space of all measurable functions on $\R_+$
that satisfy the equivalent conditions
\begin{gather}
\sup_{n\in\Z}2^{n}
 \left(\int_{2^n}^{2^{n+1}}|\f(x)|^2dx\right)<\be;
\label{xinfi}\\
 \sup_{x>0} x\int_{x}^\be|\f(y)|^2dy <\be;
\label{xinfii}
\end{gather}
Let $X^0_\be$ 
be the subspace of $X_\be$ consisting of the functions
that satisfy the equivalent conditions
\begin{gather}
\lim_{n\to\pm\be}2^{n}
 \left(\int_{2^n}^{2^{n+1}}|\f(x)|^2dx\right)=0;
\label{xinf0i}\\
\lim_{x\to0} x\int_{x}^\be|\f(y)|^2dy 
= \lim_{x\to\be} x\int_{x}^\be|\f(y)|^2dy 
=0.
\label{xinf0ii}
\end{gather}
\end{definition}

The equivalence of the different conditions is an
exercise.
For $1\le p\le \be$, $X_p$ is a Banach space with the norm
\begin{equation*}
\|\f\|_{X_p}=
\biggl\| x^{1/2}
\biggpar{\int_{x}^\be|\f(y)|^2dy}^{1/2}\biggr\|_{ L^p(dx/x)};
\end{equation*}
for $0<p<1$, 
this is a quasi-norm and
$X_p$ is a quasi-Banach space.
$X_\be^0$ is a closed subspace of $X_\be$, and thus a Banach space too.
Note that $X_p\subset X_q$ if $0<p\le q\le\be$.

\begin{remark}
It is well 
known \cite{Pee} that $\f\in X_p$ if and only if the Fourier
transform $\F\f$ belongs to the Besov 
space $B^{1/2}_{2\,\,p}$ (here we identify $\f$ with the function
extended to $\R$ by zero on $\R_-$).
\end{remark}

Note that the operators $Q_\f$ appear in a natural way in \cite{MV}
when studying the boundedness problem  
for the Sturm--Liouville operator $\cL$ from 
$\overset \circ L{}^1_2(\R_+)\to L_2^{-1}(\R_+)$
defined by $\cL u=-u''+qu$. To be more precise, 
Maz'ya and  Verbitsky studied in \cite{MV} the problem
of identifying potentials $q$ for which the inequality
$$
\left|\int_{\R_+}|u(t)|^2q(t)dx\right|\le\const\int_{\R_+}|u'(t)|^2dt
$$
holds for any $C^\be$ compactly supported function $u$ on
$(0,\be)$. This inequality is in turn 
equivalent to the boundedness of the quadratic form
\bay
\label{kvf}
\left|\int_{\R_+}u(t)\ov{v(t)}q(t)dt\right|
\le\const\|u'\|_{L^2(\R_+)}\|v'\|_{L^2(\R_+)}.
\ey
In \cite{MV} under the assumption that the limit
$$
\lim_{y\to\be}\int_x^yq(t)dt=\int_x^\be q(t)dt\df\f(x)
$$
exists for any $x>0$ the problem of boundedness (compactness) of the
quadratic form \rf{kvf} was reduced 
to the problem of boundedness (compactness) of the operator $Q_\f$ on
$L^2(\R_+)$. 
Note that in \cite{MV} the authors also obtained boundedness and
compactness criteria for the operators $Q_\f$ 
in terms of conditions \rf{xinfii} and \rf{xinf0ii}.

The conditions in the above definition
are conditions on the size of $\f$ only, and define Banach lattices of
functions on $\R_+$. 
Thus, if $|\psi|\le|\f|$ and $\qf$ is bounded, compact, or belongs to
$\bS_p$, $p>1$,  then $\qq{\psi}$ has the same property and, for
example, $\|\qq\psi\|_{\bS_p}\le C_p \|\qf\|_{\bS_p}$.
Moreover, we will see that the same conditions are necessary and sufficient for
these properties for
the operators $\qf^+$ and $\qf^-$ too; thus, if one of the three
operators has one these properties, then all three have it.
 These results for $\qf^+$ are not new, see for example
\cite{ES, EEH, No, NS, St}, and the results for $\qf$ can easily be
derived. Nevertheless we give complete proofs, by another method,
as a background for the case $p\le1$.

At $p=1$, there is a kind of threshold, and  
for $\bS_p$, $p\le1$, the situation is much more complex.
First, $\qf^+$ and $\qf^-$ never belong to $\bS_1$, except in the
trivial case $\f=0$ a.e., when the operators vanish (\refT{zero}).
Secondly, although 
$\f\in X_1$ is a necessary condition for $\qf\in\bS_1$, 
it is not sufficient.
Indeed, for
$p\le1$, we have not succeeded in finding both necessary
and sufficient conditions for $\qf\in\bS_p$, 
and any such conditions would have to be
fairly complicated. For one thing, the property $\qf\in\bS_p$ does not
depend on the size of $\f$ only; although $\f_0(x)=\chi_{\OI}$, the
characteristic function of the unit interval,
yields an operator $\qq{\f_0}$ of rank 1,
we show (see the example following Theorem \ref{E:exp}) that there
exists a function $\psi$ with 
$|\psi|=|\f_0|$ 
such that $\qq{\psi}\notin\bS_1$.
In the positive direction we show (\refS{S:yp} and \refS{more}) 
that if $\f$ is
sufficiently smooth and decays sufficiently rapidly at infinity, then
$\qf\in\bS_p$, $1/2<p\le1$.
Conversely, we give in \refS{neces} (\refT{mnc0}) 
a necessary condition on the $L^1$
modulus of continuity for $\qf\in\bS_1$. 

At $p=1/2$ there is a second threshold. 
 We prove in \refS{S:1/2} and \refS{S:Sturm}, by two different methods,
that if $\qf$ is smooth (locally absolutely continuous is enough), then
$\qf$ never belongs to $\bS_{1/2}$ except when $\f=0$ a.e.
More precisely, if $\f$ is sufficiently smooth and decays sufficiently
rapidly at infinity, and does not vanish identically, then
the singular numbers $s_n(\qf)$ decay
asymptotically exactly like $n^{-2}$ (\refT{T:lim}).

On the other hand, $\qf$ may belong to $\bS_{1/2}$ for non-smooth functions:
It is easily seen that if $\f$ is a step function, then $\qf$ has finite
rank, and thus $\qf\in\bS_p$ for every $p>0$. Taking suitable infinite
sums of step functions we find also other functions
in $\bS_p$, $p\le1/2$.

The role of smoothness is thus complicated and not well 
understood. 
It seems to be a help towards 
$\qf\in\bS_p$ for $1/2<p\le1$, but it is not necessary and 
it completely prevents $\qf\in\bS_p$ for $p\le1/2$.
On the other hand, it is irrelevant for $p>1$. 

As said above, $\qf$ has finite rank when $\f$ is a step function.
We show in \refT{T:Kronecker} that this is the only case when $\qf$
has finite rank.

The kernel in Definition \eqref{q} is symmetric, 
and thus  $\qf$ is self-adjoint if and only if $\f$ is real.
In \refS{S:positive} we show that $\qf$ is a positive operator if and
only if $\f$ is a non-negative non-increasing function.
In this special case, for each $p>1/2$,  $\qf\in\bS_p$ if and only if
$\f\in X_p$.  
In this case we also give an 
even simpler necessary and sufficient conditions 
for boundedness, compactness and $\qf\in\bS_p$, $p>1/2$
(\refT{T:mon}).
In particular, for positive operators
we have a necessary and sufficient condition for 
$p=1$ too.

When $\f$ is real and thus $\qf$ is self-adjoint, the singular values
are the absolute values of the eigenvalues. In \refS{S:Sturm}, we
study the eigenvalues, which leads to a Sturm--Liouville problem
that we study.
We include one example
(see \refS{S:Sturm}), where the singular values can be calculated exactly by
this method.
We give also another example (Theorem \ref{E:exp}) where the
singular values 
are calculated within a constant factor by Fourier analysis.

In \refS{io} and \refS{S:matrix}, 
we consider related families of operators acting on 
$L^2(\R)$ and $\ell^2(\Z)$; the latter operators 
include some given by weighted Hankel matrices.

We further study operations acting on bounded linear operators on 
$L^{2}(\R^{+})$ 
related to the class of operators $Q_\f$.
We study Schur multipliers given by functions of the form 
$\varphi\left(\max\left\{  x,y\right\}\right) $ in \refS{Schur} and 
properties of the averaging projection 
onte the operators of the form $\qf$
in \refS{aver}.

We give in this paper several necessary conditions and sufficient
conditions for 
properties such as $\qf\in\bS_p$. In all cases there are corresponding
norm estimates, which follow by inspection of the proofs or by the
closed graph theorem, although we usually do not state these estimates
explicitly.

We denote by $|I|$ the length of an interval $I$. 
We also use $|S|$ for the cardinality of a finite set $S$; there is no
danger of confusion.

We use $c$ and $C$, sometimes with subscripts or superscripts, to denote
various unspecified constants, not necessarily the same on different
occurrences. These constants are universal unless we indicate otherwise
by subscripts.

\setcounter{equation}{0}
\section{\bf Preliminaries}
\label{prel}

\

Definition \eqref{q} shows that the adjoint $\qf^*=\qq{\bar \f}$;
in particular, $\qf$ is self-adjoint if and only if $\f$ is real.
Similarly, $(\qf^+)^*=\qq{\bar \f}^-$;
which has the same norm and singular numbers as $\qf^-$.
Hence, we will mainly consider $\qf^+$;
all results obtained in this paper for $\qf^+$ immediately holds for
$\qf^-$ too.

\subsection*{Schatten classes}
We denote the singular numbers of a bounded operator $T$ on a Hilbert
space 
(or from one Hilbert space into another)
by $s_n(T)$, $n=0,1,2,\dots$; thus
$s_n(T)\df\inf\set{\|T-R\|:\rank(R)\le n}$.
We will frequently use the simple facts
\begin{equation}
\label{snsum}
s_{m+n}(T+R)\le s_{m}(T)+s_{n}(R), \qquad m,n\ge0,
\end{equation}
and
\begin{equation}
\label{snpr}
s_{m+n}(TR)\le s_{m}(T)s_{n}(R), \qquad m,n\ge0,
\end{equation}

Recall that the \emph{Schatten--von Neumann classes} $\bS_p$, $0<p<\be$
are defined by 
$$
\bS_p=\left\{T:~\sum_n  s_n(T)^p\le\be\right\},
$$
and the \emph{Schatten--Lorentz classes} $\bS_{p,q}$ are defined by
\begin{alignat*}2
\bS_{p,q}&=\left\{T:~\sum_{n\ge0}(s_n(T))^q(1+n)^{q/p-1}<\be\right\},
&&\quad0<p<\be,\;0<q<\be,
\\
\wS{p}&=\left\{T:~s_n(T)\le\const(1+n)^{-1/p}\right\},&&\quad0<p<\be.
\end{alignat*} 
See for example \cite{GK1} and \cite{BS}.

\subsection*{Other intervals}
We have defined our operators for the interval $\R_+=(0,\infty)$. More
generally, for any interval $I\subseteq\R$ and a function
$\f\in\Lloc^1(I)$, we define $\qf^I$ to be the integral operator on
$L^2(I)$ with kernel $\f\bigpar{\max(x,y)}$.

It is easily seen that if $I=(-\be,a)$ with $-\be<a\le\be$, then
$\qf^I$ is bounded only for $\f=0$ a.e. By translation invariance,
it remains only to consider the cases $I=(0,\be)$, as above, and
$I=(0,a)$ for some finite $a$. The latter case will be used sometimes
below, but it can always be reduced to the case $(0,\be)$. Indeed, if
we extend $\f$ to $(0,\be)$ letting $\f=0$ on $(a,\be)$, then $\qf^I$
and $\qf$ may be identified. (Formally, they are defined on different
spaces, and $\qf^I$ is the restriction of $\qf$ to $L^2(I)$, but the
complementary restriction to $L^2(a,\infty)$ vanishes. In particular,
$\qf^I$ and $\qf$ have the same singular numbers.)

The case of a finite interval can also be reduced to \OI{} by the
following simple homogeneity result.
\begin{lem}
\label{L:homo}
If $t>0$, and $\f_t(x)=t\f(tx)$, then $\qf$ and $Q_{\f_t}$ are
unitarily equivalent.
Similarly, for a subinterval $I\subset(0,\be)$, $\qf^I$ and
$Q_{\f_t}^{t^{-1}I}$ are unitarily equivalent.
\end{lem}

\Pf
The mapping $T_t: f(x)\mapsto t^{1/2}f(tx)$ is a unitary operator in
$L^2(\R_+)$, and $Q_{\f_t}=T_t\qf T_t^{-1}$.
$\bl$

Note that the spaces $X_p$ have the homogeneity property exhibited in
this lemma:
if $\f\in X_p$, then $\f_t\in X_p$ with the same norm.
Of course, it is natural to have this property for any necessary or
sufficient condition for $\qf\in\bS_p$.

\subsection*{Distributions} 
We can also define the operators $Q_\f$, $Q_\f^+$, and $Q_\f^-$ in the
case when $\f$ is 
a distribution. 

For an open subset $G$ of $\R^n$ we denote by $\cd(G)$ the space of
compactly supported $C^\be$ 
functions  in $G$ and denote by $\cd\,^\prime(G)$ the space of
\emph{distributions} on $G$, i.e., continuous linear functionals on $\cd(G)$.
We refer the reader to \cite{Sch} for basic facts about distributions.
We use the notation $\langle\f,f\rangle$ for $\f(f)$, where
$\f\in\cd\,^\prime(G)$ and $f\in\cd(G)$.

Suppose now that $\Sigma$ and $\Omega$ are open subsets of $\R$. In
this paper we usually 
consider the case when $\Sigma=\Omega=\R_+$ or $\Sigma=\Omega=\R$. Let
$\Phi\in\cd\,(\Sigma\times\Omega)$. We say that $\Phi$ determines a
bounded linear operator from 
$L^2(\Omega)$ into $L^2(\Sigma)$ if there exists a constant $C$ such that 
$|\langle\Phi(x,y),f(y)\overline{g(x)}\rangle|\le 
C\|f\|_{L^2}\|g\|_{L^2}$ for any $f\in\cd(\Omega)$ and any
$g\in\cd(\Sigma)$;
the corresponding operator $T$ then is given by 
$\langle Tf,g\rangle=\langle\Phi(x,y),f(y)\overline{g(x)}\rangle$
and $\Phi$ is called the \emph{kernel} of $T$.
We say that $\Phi$ determines an operator in $\bS_p$ if this operator
is an operator  
from $L^2(\Omega)$ into $L^2(\Sigma)$ of class $\bS_p$.
Note that for any bounded operator $T:L^2(\Omega)\to L^2(\Sigma)$ there exists
a distribution $\Phi\in\cd\,'(\Sigma\times\Omega)$ which determines
the operator $T$
(a special case of Schwartz's kernel theorem). 
Indeed, it is easy to see that any function in
$\cd(\Omega\times\Sigma)$ defines 
an operator from $L^2(\Sigma)$ into $L^2(\Omega)$ of class $\bS_1$, and we have
a continuous imbedding
$j:\cd(\Omega\times\Sigma)\to\bS_1(L^2(\Sigma),L^2(\Omega))$. 
We may define the distribution $\Phi_T\in\cd\,(\Sigma\times\Omega)$ by
the following 
formula $\langle\Phi_T(x,y),f(x,y)\rangle\df\trace(TA)$, where
$f\in\cd(\Sigma\times\Omega)$ and 
$A$ is the integral operator with kernel function $f(y,x)$. Clearly,
$\Phi_T$ determines 
the operator $T$.

We also consider the space ${\cal S}(\R^n)$ of infinitely smooth
functions whose derivatives 
of arbitrary orders
decay at infinity faster than any power of $(1+|x|)^n$ and the dual
space ${\cal S}\,'(\R^n)$  
of \emph{tempered distributions} (see \cite{Sch} for basic facts).
Recall that the Fourier transform 
\begin{equation}
\label{fourier}
f\mapsto (\F f)(x)\df\int\limits_{\R^n}f(t)e^{-2\pi{\rm i}(t,x)}\,dt
\end{equation}
where $(t,x)$ is the scalar product of $x$ and $t$ in $\R^n$,
is an isomorphism of ${\cal S}(\R^n)$ onto itself, and that it
can be extended to ${\cal S}'(\R^n)$ by duality.

We need the following elementary facts.

\begin{lem}
\label{temp}
Suppose that a distribution $\Phi\in\cd\,^\prime(\R^2)$ determines a
bounded operator 
on $L^2(\R)$. Then $\Phi$ is a tempered distribution.
\end{lem}

\Pf It is easy to see that any function $\Phi\in{\cal S}(\R^2)$ determines
an operator of class $\bS_1$ 
and that 
the corresponding imbedding of ${\cal S}(\R^2)$ 
into $\bS_1$ is continuous. The result follows now by duality.
$\bl$

\begin{lem}
\label{fue}
Let $\Phi\in{\cal S}^\prime(\R^2)$. Consider the distribution $\Psi$
on $\R^2$ defined by 
$\Psi(x,y)\df(\F\Phi)(x,y)$.
Then $\Phi$ determines a bounded  operator on $L^2(\R)$ if and only
if $\Psi$ does. 
Moreover, these two operators are unitarily equivalent.
\end{lem}

\Pf It suffices to observe that 
\bey
\langle\Psi(x,y),f(y)\overline{g(x)}\rangle&=&
\langle\Phi(x,y),\F(f(y)\overline{g(x)})\rangle\\[.2cm]
&=&\langle\Phi(x,y),(\F f)(y)\overline{\F g(-x)})\rangle.\quad\bl
\eey

Now we are ready to define the operators $Q_\f$, $\qf^+$ and $\qf^-$ in the
case where $\f$ is distribution.

It is not hard 
to see that the operator $f\mapsto\int\limits_0^xf(x,y)dy$
is a continuous operator 
from $\cd(\R_+\times\R_+)$ into $\cd(\R_+)$. Hence, with any 
$\f\in\cd\,^\prime(\R_+)$ we can associate the distributions 
$\Lambda_\f^+$, $\Lambda_\f^-$,and $\Lambda_\f$ in
$\cd\,^\prime(\R_+\times\R_+)$ defined 
by 
\begin{align}
\langle\Lambda_\f^+,f(x,y)\rangle&\df\langle\f,\int\limits_0^xf(x,y)dy\rangle,
\label{lp+}\\ 
\langle\Lambda_\f^-,f(x,y)\rangle&\df\langle\f,\int\limits_0^xf(y,x)dy\rangle
\label{lp-}\\ 
\intertext{and}
\Lambda&\df\Lambda_\f^++\Lambda_\f^-.
\label{lp}
\end{align}

For a distribution $\f$ in $\cd\,'(\R_+)$, we can consider now the
operators $Q_\f^+$, $Q_\f^-$, and 
$Q_\f$ determined by the distributions $\Lambda_\f^+$,
$\Lambda_\f^-$,and $\Lambda_\f$ respectively. 
It is easy to see that in case $\f\in \Lloc^1(\R_+)$, the new
definition coincides with the 
old one. The following theorem shows however that if one of those
operators is bounded on 
$L^2(\R_+)$, then $\f$ must be a locally integrable function on
$\R_+$, and so we have not 
enlarged the class of bounded operators of the form $Q_\f$.

\begin{thm}
\label{ker+} 
Let $\f\in\cd\,^\prime(\R_+)$. Suppose that at least one of the
distributions $\Lambda_\f^+$, 
$\Lambda_\f^-$, or $\Lambda_\f$ determines a bounded operator on
$L^2(\R_+)$. Then $\f\in 
\Lloc^2(\R_+)$.
\end{thm}

\Pf We consider the cases of the distributions $\Lambda_\f^+$ and $\Lambda_\f$.
For $\Lambda_\f^-$ the proof is the same as for $\Lambda_\f^+$.
Let $a\in\R_+$. Fix a function $f_0\in\cd(\R_+)$ such that
$\supp f_0\subset[0,a]$ and $\int_{\R_+}f_0(x)\,dx=1$. We have
$$
\langle\Lambda_\f^+,f_0(y)\overline{g(x)}\rangle=
\langle\Lambda_\f,f_0(y)\ov{g(x)}\rangle=\langle\f,\bar g\rangle
$$ 
for any $g\in\cd(\R_+)$  with $\supp g\subset[a,+\infty)$.  Therefore 
$$
|\langle\f,\overline g\rangle|\le C
\|f_0\|_{L^2(\R_+)}\|g\|_{L^2(\R_+)}
$$ for any $g\in\cd(a,\be)$. 
Thus, $\f\big|(a,\infty)\in L^2(a,\infty)$. $\bl$

\subsection*{Triangular projection}
On the class of operators on $\bS_p(L^2(\R_+))$, $p<\be$, we define the
operator of triangular projection 
$\cp$ as follows. Consider first the case $p\le2$. Let $T$ be an
operator on $L^2(\R_+)$ of class 
$\bS_p$, $p\le2$. Then $T$ is an integral operator with kernel function $k_T$:
$$
(Tf)(x)=\int_0^\be k_T(x,y)f(y)\,dy,\quad f\in L^2(\R_+).
$$
Then by definition 
\bay
\label{pr}
(\cp Tf)(x)=\int_0^x k_T(x,y)f(y)\,dy,\quad f\in L^2(\R_+).
\ey
It is well known that 
\bay
\label{tri}
\|\cp T\|_{\bS_p}\le c_p\|T\|_{\bS_p},\quad1<p\le2,
\ey
where $c_p$ depends only on $p$. This allows one to extend
by duality the definition of 
$\cp$  and inequality \eqref{tri} to the case $2\le p<\be$.
Note also that $\cp$ has weak type $(1,1)$, i.e.,
\bay
\label{P11}
s_n(\cp T)\le\const(1+n)^{-1}\|T\|_{\bS_1},\quad T\in\bS_1.
\ey

We will need these results on the triangular projection $\cp$ in a
more general situation. 
Let $\mu$ and $\nu$ be regular Borel measures on $\R_+$ and let $K$ be
a Borel function 
on $\R_+\times\R_+$. As above we can associate with any operator $T$ from  
$L^2(\mu)$ to $L^2(\nu)$ of class $\bS_2$ the operator $\cp T$ by
multiplying the 
kernel function of $T$ by the characteristic function of the set
$\{(x,y)\in\R_+^2:~0<y<x\}$. 

\begin{thm}
\label{itr}
$\cp$ is a bounded linear projection on $\bS_p(L^2(\mu),L^2(\nu))$ for
\linebreak$1<p<\be$ 
and $\cp$ has weak type $(1,1)$, i.e., $\cp$ is a bounded linear operator from 
$\bS_1(L^2(\mu),L^2(\nu))$ to $\bS_{1,\be}(L^2(\mu),L^2(\nu))$.
\end{thm}

\refT{itr} is well known at least when $\mu=\nu$. Let us explain how
to reduce the proof of 
Theorem \ref{itr} to the case of the triangular projection onto the
upper triangular matrices. 

Let $\{\K_j\}_{j\ge0}$ and $\{\h_k\}_{k\ge0}$ be Hilbert spaces.
Put $\K\df\bigoplus\limits_{j\ge0}\K_j$ and
$\h\df\bigoplus\limits_{k\ge0}\h_k$. 
We identify operators $A\in\B(\h,\K)$ with its block matrix
representation $\{A_{jk}\}_{j,k\ge0}$, where 
$A_{jk}\in\B(\h_k,\K_j)$. We define the triangular projection $\bs{\cp}$ by
$(\bs{\cp}A)_{jk}\df A_{jk}$ for $j>k$ and $(\bs{\cp}A)_{jk}\df0$ for
$j\le k$.  

\begin{lem}
\label{mtr}
Let $1<p<\be$. Then $\bs{\cp}$ is bounded on $\bS_p(\h,\K)$ and has
weak type $(1,1)$. 
Moreover, the norms of $\bs{\cp}$ 
can be bounded independently of $\h$ and $\K$. 
\end{lem} 

In the case $\dim\K_j=\dim\h_k=1$, this is the Krein--Matsaev theorem
(it is equivalent to Theorems 
III.2.4 and III.6.2 of \cite{GK2}, see also Theorem IV.8.2 of \cite{GK1}). 
In general the result can be
reduced easily to this special  
case. Indeed, it is easy to reduce the general case to the case when
$\dim\h_j=\dim\K_j<\be$. 
Then it is easy to see that if $A\in\bS_p$, $1\le p<\be$, then the
diagonal part of 
$A$ also belongs to $\bS_p$, and so we may assume without loss of
generality that $A_{jj}=0$, 
$j\in\Z_+$. We can take an orthonormal basis in each $\h_j$ and
consider the orthonormal basis 
of $\h$ that consists of those basis vectors of $\h_j$,
$j\in\Z_+$. Then we can consider the matrix 
representation of $A$ with respect to this orthogonal basis. We have
now two triangular projections: 
with respect to the orthogonal basis and the projection $\bs{\cp}$,
the triangular with respect to the 
the decomposition $\h=\bigoplus\limits_{k\ge0}\h_k$. It is not hard to
check that  
since the diagonal part of 
$A$ is zero, both triangular projections applied to $A$ give the same
result. This reduces the general 
case to the Krein--Matsaev theorem mentioned above. $\bl$

Now it is easy to deduce Theorem \ref{itr} from Lemma \ref{mtr}.

{\bf Proof of Theorem \ref{itr}.} Let $T$ be an integral operator with
kernel function $k$. 
For $\e>0$ we put $k_\e(x,y)\df k(x,y)\chi_{\{(x,y)\in\R^2_+:[\frac
x\e]\e>y>0\}}$, where $[a]$ denote 
the largest integer that is less than or equal to $a$. Suppose that $p>1$. 
It is sufficient to consider the case $1<p\le2$ and then use
duality. Let $T_\e$ be 
the integral operator with kernel function $k_\e$.
It follows easily from Lemma \ref{mtr} that
$$
\|T_\e\|_{\bS_p(L^2(\mu),L^2(\nu))}\le c_p\|T\|_{\bS_p(L^2(\mu),L^2(\nu))}
$$
for any $\e>0$. Clearly, $T_\e\to T$ in the weak operator topology
as $\e\to0$. 
It
follows that  
$\|\cp T\|_{\bS_p(L^2(\mu),L^2(\nu))}\le
c_p\|T\|_{\bS_p(L^2(\mu),L^2(\nu))}$. The case $p=1$ may be 
considered in the same way.
$\bl$

We have $Q_\f^+=\cp Q_\f$, which  together with the equivalence
$$
\qf^+\in\bS_p \iff \qf^-\in\bS_p
$$
yields
$$
\qf^+\in\bS_p \iff \qf\in\bS_p
$$ 
for $1<p<\be$. We will give a direct
proof of this in \refT{T:p>1}.

We introduce a more general operation. Let $A$ be a measurable subset of 
\linebreak$\R_+\times\R_+$. For an operator $T$ on $L^2(\R_+)$ of
class $\bS_2$ with 
kernel function $k_T$ we consider the integral operator $\cp_AT$ whose
kernel function 
is $\chi_A k_T$, where $\chi_A$ is the characteristic function of
$A$. In other words, 
$$
(\cp_ATf)(x)=\int_0^\be\chi_A(x,y)k_T(x,y)f(y)dy.
$$
If $0<p\le2$ and $\cp_A$ maps $\bS_p$ into itself, it follows from the
closed graph theorem that 
the linear transformation $\cp_A$ is a bounded linear operator on
$\bS_p$. If $1<p\le2$ and $\cp_A$ is a 
bounded linear operator on $\bS_p$, then by duality we can define in a
natural way the bounded linear 
operator $\cp_A$ on $\bS_{p'}$. If $\cp_A$ is bounded on $\bS_1$, we
can define by duality $\cp_A$ 
on the space $\B(L^2(\R_+))$ of bounded linear operators on
$L^2(\R_+)$. Note that the projection $\cp$ 
defined by \eqref{pr} is equal to $\cp_A$ with
$A=\{(x,y):~x\in\R_+,~y\le x\}$. 

\

\setcounter{equation}{0}
\section{\bf Boundedness, compactness, and $\bs{p>1}$}
\label{S:p>1}

\

Recall the spaces $X_p$ defined in the Introduction.

\begin{thm}
\label{bou}
Let $\f\in L^2_{\textup{loc}}(\R_+)$. 
The following are equivalent:
\begin{enumerate}
\item 
$Q_\f$ is bounded on $L^2(\R_+)$;
\vspace*{.2cm}
\item 
$Q^+_\f$ is bounded on $L^2(\R_+)$;
\vspace*{.2cm}
\item 
$\f\in X_\infty$.
\end{enumerate}
\end{thm}

Recall that the equivalence of (i) and (iii) was also 
established in \cite{MV} by a different method.

\Pf 
Let us show that (ii) implies (i). If $Q^+_\f$ is bounded, then the
integral operator 
$Q_\f^-\df Q_\f-Q_\f^+$ is also bounded, since its kernel function is
the reflection of the kernel 
function of $Q_\f^+$ with respect to the line $\{x=y\}$. Hence, $Q_\f$
is bounded. 

Let us deduce now (iii) from (i). For $n\in\Z$ put
$$
A_n=[2^n,2^{n+1}]\times[2^{n-1},2^n].
$$
Certainly, 
if $Q_\f$ is bounded, then
$$
\sup_{n\in\Z}\|\cp_{A_n}Q_\f\|<\be.
$$
It is easy to see that $\cp_{A_n}Q_\f$ is a rank one operator and
$$
\|\cp_{A_n}Q_\f\|=\left(2^{n-1}\int_{2^n}^{2^{n+1}}|\f(x)|^2dx\right)^{1/2}
$$
which implies \eqref{xinfi}.

It remains to prove that \eqref{xinfi} implies 
(ii). Put
\bay
\label{B}
B=\bigcup_{n\in\Z}B_n,
\ey
where
\bay
\label{Bn}
B_n=\{(x,y):~2^n\le x\le 2^{n+1},~2^n<y<x\}.
\ey
We also define the sets
\bay
\label{Ank}
A_n^{(k)}=[2^n,2^{n+1}]\times[2^{n-k},2^{n-k+1}]
\ey
and
\bay
\label{Ak}
A^{(k)}=\bigcup_{n\in\Z}A_n^{(k)}.
\ey
Clearly,
$$
\{(x,y):~x>0,~0<y<x\}=B\cup\bigcup_{k\ge1}A^{(k)},
$$
and so
\bay
\label{Q+}
\|Q_\f^+\|\le\|\cp_BQ_\f\|+\sum_{k\ge1}\|\cp_{A^{(k)}}Q_\f\|.
\ey

Since the projections of the $B_n$ onto the coordinate axes are
pairwise disjoint, it is straightforward 
to see that
$$
\|\cp_BQ_\f^+\|=\sup_{n\in\Z}\|\cp_{B_n}Q_\f^+\|.
$$
Let $R_n$ be the integral operator with kernel function 
\bay
\label{Rn}
k_{R_n}(x,y)=\f(x)\chi_{[2^n,2^{n+1}]}(x)\chi_{[2^n,2^{n+1}]}(y).
\ey
Obviously,
$\rank R_n=1$ and 
$\|R_n\|_{\bS_2}=\|R_n\|
=\left(2^n\int_{2^n}^{2^{n+1}}|\f(x)|^2dx\right)^{1/2}$.
It is also evident that 
$\cp R_n=\cp_{B_n}Q_\f^+$, and  since $\cp$ an orthogonal projection
on $\bS_2$, we have 
$$
\|\cp_{B_n}Q_\f^+\|=\|\cp R_n\|=\|\cp R_n\|_{\bS_2}\le\|R_n\|_{\bS_2}=
\left(2^n\int_{2^n}^{2^{n+1}}|\f(x)|^2dx\right)^{1/2}\le\const_\f.
$$

Next, it is also easy to see that
$$
\|\cp_{A^{(k)}}Q_\f\|=\sup_{n\in\Z}\|\cp_{A_n^{(k)}}Q_\f\|.
$$
Also, $\cp_{A_n^{(k)}}Q_\f$ has rank one and norm
$\left(2^{n-k}\int_{2^n}^{2^{n+1}}|\f(x)|^2dx\right)^{1/2}$, and so
\bey
\sum_{k\ge1}\|\cp_{A^{(k)}}Q_\f\|&=&
\sum_{k\ge1}\sup_{n\in\Z}
 \left(2^{n-k}\int_{2^n}^{2^{n+1}}|\f(x)|^2dx\right)^{1/2}\\
&=&\sum_{k\ge1}2^{-k/2}\sup_{n\in\Z}
 \left(2^n\int_{2^n}^{2^{n+1}}|\f(x)|^2dx\right)^{1/2}\\
&=&\const\sup_{n\in\Z}
 \left(2^n\int_{2^n}^{2^{n+1}}|\f(x)|^2dx\right)^{1/2}.
\eey
The result follows now from \eqref{Q+}. $\bl$

\begin{thm}
\label{com}
Let $\f\in \Lloc^1(\R_+)$. The following are equivalent:
\begin{enumerate}
\item  
$Q_\f$ is compact on $L^2(\R_+)$;
\vspace*{.3cm}
\item  
$Q^+_\f$ is compact on $L^2(\R_+)$;
\vspace*{.3cm}
\item  
$\f\in X_\be^0$.
\end{enumerate}
\end{thm}

Recall that the equivalence of (i) and (iii) was also established in
\cite{MV} by a different method. 

\Pf
It is easy to see that the estimates given in the proof of Theorem
\ref{bou} actually 
lead to the proof of Theorem \ref{com};
for the step (i)$\implies$(iii) we observe that 
if $Q_\f$ is compact, then
$\lim_{n\to\pm\be}\|\cp_{A_n}Q_\f\|=0$.
$\bl$

\begin{thm}
\label{T:p>1}
Let $1< p<\be$ and let $\f\in \Lloc^1(\R_+)$. 
Then the following conditions are equivalent:
\begin{enumerate}
\item \label{tpi}
$Q_\f\in\bS_p$;
\vspace*{.3cm}
\item \label{tpii}
$Q^+_\f\in\bS_p$;
\vspace*{.3cm}
\item \label{tpiii}
$\f\in X_p$.
\end{enumerate}
\end{thm}

Note that the fact that (ii)$\Leftrightarrow$(iii) was proved in
\cite{No} by a different method, 
see also \cite{NS} and \cite{St} for the case of more general Volterra
operators. 

\Pf
The fact that (ii)$\imp$(i) can be proved exactly as in the
proof of Theorem \ref{bou}. 
Let 
us show that (i) implies 
(iii). Consider the sets $A_n=A_n^{(1)}$ introduced in the proof of
Theorem \ref{bou} (see \rf{Ank}).
Recall that $A^{(1)}=\bigcup_{n\in\Z}A_n^{(1)}$. It is easy to see
that
$$
\|\cp_{A^{(1)}}Q_\f\|_{\bS_p}\le\|Q_\f\|_{\bS_p}.
$$
Clearly,
$$
\|\cp_{A^{(1)}}Q_\f\|_{\bS_p}^p=\sum_{n\in\Z}\|\cp_{A_n}Q_\f\|_{\bS_p}^p,
$$
the operator $\cp_{A_n}Q_\f$ has rank one and 
$$
\|\cp_{A_n}Q_\f\|_{\bS_p}=
 \left(2^{n-1}\int_{2^n}^{2^{n+1}}|\f(x)|^2dx\right)^{1/2}.  
$$
This implies \eqref{xpi}.

Let us show that \eqref{xpi} implies
(ii). Consider the sets, $B_n$, $B$, 
$A_n^{(k)}$, $A^{(k)}$ defined in \rf{Bn}, \rf{B}, \rf{Ank}, and
\rf{Ak}. 
Clearly,
$$
\|Q^+_\f\|_{\bS_p}
\le\|\cp_B Q^+_\f\|_{\bS_p}+\sum_{k\ge1}\|\cp_{A^{(k)}}Q^+_\f\|_{\bS_p}.
$$
Let us first estimate $\|\cp_B Q^+_\f\|_{\bS_p}$. Clearly,
$$
\|\cp_B Q^+_\f\|_{\bS_p}^p=\sum_{n\in\Z}\|\cp_{B_n}Q^+_\f\|_{\bS_p}^p.
$$
Consider the rank one operators $R_n$ defined in \eqref{Rn}. As in the
proof of Theorem \ref{bou} 
we have $\cp R_n=\cp_{B_n}Q_\f^+$ and since $\cp$ is bounded on
$\bS_p$, we obtain 
$$
 \|\cp_{B_n}Q_\f^+\|_{\bS_p}=\|\cp R_n\|_{\bS_p}\le\const_p\|R_n\|_{\bS_p}=
\const_p\left(2^n\int_{2^n}^{2^{n+1}}|\f(x)|^2dx\right)^{1/2},
$$
and so 
$$
\|\cp_B Q^+_\f\|_{\bS_p}^p\le
\const_p\sum_{n\in\Z}\left(2^n\int_{2^n}^{2^{n+1}}|\f(x)|^2dx\right)^{p/2}.
$$

It is also easy to see that 
$$
\|\cp_{A^{(k)}}Q^+_\f\|_{\bS_p}^p
=\sum_{n\in\Z}\|\cp_{A_n^{(k)}}Q^+_\f\|_{\bS_p}^p
$$
and, since $\cp_{A_n^{(k)}}Q_\f$ has rank 1, 
$$
\|\cp_{A_n^{(k)}}Q^+_\f\|_{\bS_p}=\|\cp_{A_n^{(k)}}Q_\f\|_{\bS_p}=
\left(2^{n-k}\int_{2^n}^{2^{n+1}}|\f(x)|^2dx\right)^{1/2},
$$
and so
\begin{equation}\label{3.7a}
\begin{split}
\|\cp_{A^{(k)}}Q^+_\f\|_{\bS_p}
&= \left(\sum_{n\in\Z}\left(2^{n-k}
  \int_{2^n}^{2^{n+1}}|\f(x)|^2dx\right)^{p/2}\right)^{1/p}\\
&= 2^{-k/2}\Biggpar{\sum_{n\in\Z}\left(2^{n}
  \int_{2^n}^{2^{n+1}}|\f(x)|^2dx\right)^{p/2}}^{1/p}
\end{split}
\end{equation}
and 
\begin{equation}\label{upper}
\sum_{k\ge1}\|\cp_{A^{(k)}}Q^+_\f\|_{\bS_p}
\le\const
\Biggpar{\sum_{n\in\Z}\left(2^{n}
  \int_{2^n}^{2^{n+1}}|\f(x)|^2dx\right)^{p/2}}^{1/p}
\end{equation}
which completes the proof.  $\bl$

\begin{remark}
The same proof shows that for $p=1$, 
$\ref{tpii}\imp\ref{tpi}\imp\ref{tpiii}$,
but the final part of it fails because the triangular projection is not
bounded on $\bS_1$.
We will see later that, indeed, none of the implications can be
reversed 
for $p=1$.

In the Hilbert--Schmidt case $p=2$, the result simplifies further.
Indeed, we have the equalities
\begin{equation*}
\begin{split}
\|\qf\|_{\bS_2}
=\left(\into\into|\f(\max(x,y))|^{2}\,dx\,dy\right)^{1/2}
=\left(2\into x|\f(x)|^{2}\,dx\right)^{1/2}
\end{split}
\end{equation*}
and 
$$
\|\qf^+\|_{\bS_2}=2^{-1/2}\|\qf\|_{\bS_2}=\|x^{1/2}\f(x)\|_2=\|\f\|_{X_2}.
$$
\end{remark}

\

\setcounter{equation}{0}
\section{\bf Positive operators}
\label{S:positive}

\

We consider the special case when $\qf$ is a positive operator, \ie{}
$\inner{\qf f,f}\ge0$ for every $f\in L^2(\R_+)$. In this case we
obtain rather complete results. We first characterize the corresponding
symbols $\f$.

\begin{thm}
\label{T:positive}
Suppose that $\f\in \Lloc^1(\R_+)$
is such that $\qf$ is a bounded operator.
Then $\qf$ is a positive operator if and only if $\f$ is a.e.\ equal
to a non-increasing, non-negative function.
\end{thm}

{\bf Proof.}
Suppose that $\qf$ is positive.
Define, for $z,h>0$, $\fzh=h^{-1}\chi_{(z,z+h)}$
and let $\Leb(\f)$ be the set of  Lebesgue points of $\f$.
Then, if $z\in\Leb(\f)$,
\begin{equation*}
\begin{split}
|\inner{\qf\fzh,\fzh}-\f(z)|
&=
h^{-2}\left|\int_z^{z+h}\int_z^{z+h}\big(\f\{\max(x,y)\}-\f(z)\big)
  \,dx\,dy\right|\\
&\le
h^{-2}\int_z^{z+h}\int_z^{z+h}\bigpar{|\f(x)-\f(z)|+|\f(y)-\f(z)|}
  \,dx\,dy\\
&=2h^{-1}\int_z^{z+h}|\f(x)-\f(z)| \,dx
\to0
\end{split}
\end{equation*}
as $h\to0$.
Since $\inner{\qf\fzh,\fzh}\ge0$ for all $h>0$, this implies
$\f(z)\ge0$.

Moreover, if $z_1,z_2\in\Leb(\f)$ are two Lebesgue points with $z_1<z_2$
and $0<h<z_2-z_1$, then, similarly,
\begin{equation*}
\inner{\qf\fzhi,\fzhii}
=\inner{\qf\fzhii,\fzhi}
=h^{-1}\int_{z_2}^{z_2+h}\f(y) \,dy
\to\f(z_2)
\end{equation*}
as $h\to0$, and thus, with $g_h=\fzhi-\fzhii$,
\begin{equation*}
\inner{\qf g_h,g_h}\to\f(z_1)+\f(z_2)-2\f(z_2)
=\f(z_1)-\f(z_2).
\end{equation*}
Hence $\f(z_1)\ge\f(z_2)$.

It follows that the function 
$\tilde\f(x)=\sup\set{\f(z):z\ge x,\, z\in\Leb(\f)}$ is
non-negative and non-increasing, and that $\f=\tilde\f$ a.e.

Conversely, if $\f$ is non-negative and non-increasing, then
$\lim_{x\to\be}\f(x)=0$, since a positive lower bound is impossible by
\refT{bou}. Thus there exists a measure $\mu$ on $(0,\be)$ such that
$\f(x)=\mu(x,\be)$ a.e.
If, say, $f$ is bounded with compact support in $(0,\be)$, then by
Fubini's theorem
\begin{equation*}
\begin{split}
\inner{\qf f,f}
&=\iint \f\bigpar{\max\{x,y\}} f(x)\bar f(y)\,dx\,dy
=\iiint\limits_{\max\{x,y\}<z} f(x)\bar f(y)\,dx\,dy\,d\mu(z)\\
&=\into\left|\int_0^z f(x)\,dx\right|^2d\mu
\ge0.
\end{split}
\end{equation*}
Hence $\qf$ is a positive operator.
\qed

We have used the fact that $\qf$ is a sum of the Volterra operators
$\qf^+$ and $\qf^-$. Operators of the type $\qf{}$  also appear as the
composition of Volterra operators.

\begin{thm}
\label{fact}
Suppose that $\psi_1$ and $\psi_2$ are functions on $\R_+$ such that
$Q_{\psi_1}$ and $Q_{\psi_2}$ 
are bounded linear operators. Let $\f$ be the function defined by
\bay
\label{fipsi}
\f(x)=\int_x^\be\psi_1(t)\psi_2(t)dt.
\ey
Then the operator $Q_\f$ is bounded and admits a factorization 
$$
Q_\f=Q^-_{\psi_1}Q^+_{\psi_2}.
$$
\end{thm} 

{\bf Proof.} 
Let $k_1$ be the kernel function of $Q^-_{\psi_1}$ and $k_2$ the 
kernel function of 
$Q^+_{\psi_2}$. We have
\begin{equation*}
k_1(x,t)=\begin{cases}
\psi_1(t),&t\ge x\\
0,&t<x
\end{cases}
\end{equation*}
and
\begin{equation*}
k_2(t,y)=\begin{cases}
0,&t<y\\
\psi_2(t),&t\ge y.
\end{cases}
\end{equation*}
Then the kernel function $k$ of the product $Q^-_{\psi_1}Q^+_{\psi_2}$
is given by 
$$
k(x,y)=\int_{\R_+}k_1(x,t)k_2(t,y)dt
=\int^\be_{\max\{x,y\}}\psi_1(t)\psi_2(t)dt=\f(\max\{x,y\})
$$
by the hypotheses; the integrals converge by \refT{bou} and the 
\CSineq. 
$\bl$

The function $\f$ in \eqref{fipsi} is always locally absolutely
continuous. In order to treat more general 
non-increasing $\f$, we define,
for a positive measure 
$\mu$ on $(0,\infty)$, the operator
$\qmu:L^2(0,\be)\to L^2(\mu)$ by $\qmu f(x)=\int_0^x f(y)\,dy$.
(Thus, the operator itself does not depend on $\mu$; only its range
space does.)
We have the following analogues of Theorems \refand{bou}{T:p>1}.
(We leave the corresponding criterion for compactness to the reader.)

\begin{thm}
\label{T:qmubou}
Let $\mu$ be a measure on $\R_+$. The following are equivalent:
\begin{enumerate}
\item 
$\qmu$ is bounded operator from $L^2(\R_+)$ to $L^2(\mu)$;
\vspace*{.2cm}
\item 
$\displaystyle\sup_{n\in\Z}2^n\mu[{2^n},2^{n+1})<\be$;
\vspace*{.2cm}
\item 
$\displaystyle\sup_{x>0}\,x\mu[x,\be)<\be$.
\end{enumerate}
\end{thm}

\begin{thm}
\label{T:qmu}
Let $1< p<\be$ and let $\mu$ be a measure on $\R_+$. 
The following conditions are equivalent:
\begin{enumerate}
\item 
$\qmu\in\bS_p$;
\vspace*{.3cm}
\item 
$\displaystyle\sum_{n\in\Z}2^{np/2}
\bigpar{\mu[{2^n},2^{n+1})}^{p/2}<\be$;
\vspace*{.2cm}
\item 
$\displaystyle\sum_{n\in\Z}2^{np/2}
\bigpar{\mu[{2^n},\be)}^{p/2}<\be$;
\vspace*{.2cm}
\item 
$\displaystyle
 x^{1/2}\bigpar{\mu(x,\be)}^{1/2} \in L^p(dx/x)$.
\end{enumerate}
\end{thm}

The proofs of Theorems \ref{T:qmubou} and \ref{T:qmu} are almost the same 
as the proofs of Theorems \refand{bou}{T:p>1}. 
The main difference is that we 
have to apply the theorem on the boundedness of the triangular
projection on $\bS_p$, $1<p<\be$, 
in the case of weighted $L^2$ spaces (see Theorem \ref{itr}). $\bl$

Furthermore, the factorization in \refT{fact} extends.
\begin{thm}
\label{T:qmufact}
Suppose that $\mu$ is a measure on $\R_+$ such that
$\qmu$ is a
bounded linear operator. Let $\f$ be the function defined by
$\f(x)=\mu(x,\be)$.
Then the operator $Q_\f$ is bounded and 
$
Q_\f=(\qmu)^*\qmu.
$
\end{thm} 

{\bf Proof.} 
By \refT{T:qmubou}, $0\le \f(x)\le C_\f/x$, and thus $\qf{}$ is bounded by
\refT{bou}. 

If, say, $f,g\in L^2(\R_+)$ are non-negative, then by Fubini's theorem
\begin{equation*}
\begin{split}
\inner{(\qmu)^*\qmu f,g}
&=\inner{\qmu f,\qmu g}
=\int_0^\be \int_0^z f(x)\,dx \int_0^z \bar g(y)\,dy\, d\mu(z)\\
&=\iiint\limits_{\max\{x,y\}<z} f(x)\bar g(y)\,dx\,dy\,d\mu(z)
=\inner{\qf f,g}.\quad\bl
\end{split} 
\end{equation*}

For positive operators $\qf$, we have a simple result.
(For (i), 
cf.\ the discussion of the Hille condition in 
\cite{MV}.)

\begin{thm}
\label{T:mon}
Suppose that $\f$ is a non-negative, non-increasing function on $\R_+$.
\begin{enumerate}
\item
$\qf{}$ is bounded if and only if $x\f(x)$ is bounded.
\item
$\qf{}$ is compact if and only if $x\f(x)\to0$ as $x\to0$ and as
$x\to\be$.
\item
If $1/2<p<\be$, then the following are equivalent:
\begin{enumerate}
\item
$\qf\in \bS_p$;
\item
$\f\in X_p$;
\item
$x\f(x)\in L^p(dx/x)$.
\end{enumerate}
\end{enumerate}
\end{thm}

{\bf Proof.}
The equivalence of $\f\in X_p$ and $x \f(x)\in L^p(dx/x)$ for
non-increasing, non-negative
$\f$ is elementary, using $\f(2^n)^2\ge 2^{-n}\int_{2^n}^{2^{n+1}}
|\f|^2 \ge \f(2^{n+1})^2$.
Hence, (i) follows from \refT{bou}, and (iii) for $p>1$ from
\refT{T:p>1}; furthermore, (ii) follows similarly from \refT{com}.

For (iii) for a general $p>1/2$, we first note that any of the three
conditions (a), (b) and (c) implies that $x\f(x)$ is bounded.
(This follows by (i) for (a), and by elementary estimates for (b) and
(c).) 
We 
can assume without loss of generality that $\f$ is
right-continuous on $(0,\be)$. 
 If we let $\mu$ be the measure on $\R_+$ with
$\mu(x,\infty)=\f(x)$, then by Theorems \refand{T:qmubou}{T:qmufact},
$\qmu$ is bounded and $\qf=(\qmu)^*\qmu$. 
Hence, $\qf\in\bS_p \eq \qmu\in\bS_{2p}$, and the result follows by
\refT{T:qmu}. 
$\bl$

We will see in the example given at the beginning of \refS{S:Sturm}
that \refT{T:mon}\vvp(iii)   
does not extend to $p\le1/2$.

\

\setcounter{equation}{0}
\section{\bf A sufficient condition, $\bs{1/2<p\le1}$}
\label{S:yp}

\

By linearity, we 
immediately obtain from  \refT{T:mon}  
a sufficient, but not necessary, condition for 
general symbols $\f$.

\begin{definition}
$Y_p$ is the subspace of $X_p$ spanned by non-increasing
functions. I.e., $\f\in Y_p$ if and only if $\re\f$ and $\im\f$ both
are differences of non-increasing functions in $X_p$.
\end{definition}

\begin{thm}
\label{T:yp}
Let $p>1/2$. If $\f\in Y_p$, then $\qf\in \bS_p$.
\nopf
\end{thm}

The condition $\f\in Y_p$ can be made more explicit and useful as follows. 
We denote by $\norm{\f}'_{BV(I)}$ the total variation of a function 
$\f$ over an interval $I$, and let 
$\norm{\f}_{BV(I)}\df\norm{\f}_{BV(I)}+\sup_I|\f|$. 
Moreover, we let $V_\f(x)$ denote the total variation of a function
$\f$ over the interval 
$[x,\be)$. Note that if $\f$ is locally absolutely continuous, then
$V_\f(x)=\int_x^\be|\f'(y)|\,dy$.

\begin{lem}
\label{L:yp}
Let $0<p<\be$.
If $\f$ is non-increasing, then
$$
\f\in Y_p \eq \f\in X_p \eq 
\int_0^\be|x\f(x)|^p\frac{dx}x<\be.
$$
\end{lem}

\Pf
For non-increasing $\f$, the first equivalence follows from the
definition of $Y_p$, 
while the second equivalence was noted in the proof of \refT{T:mon}.
\qed

\begin{thm}
\label{T:yp1}
Let $\f$ be a function on $\R_+$ and let $0<p<\be$.
The following are equivalent:
\begin{enumerate}
\item
$\f\in Y_p$;
\item
$V_\f\in X_p$ and $\lim\limits_{x\to\be}\f(x)=0$;
\item
$xV_\f(x)\in L^p(dx/x)$ and $\lim\limits_{x\to\be}\f(x)=0$;
\item
$\f$ has locally bounded variation, $\lim\limits_{x\to\be}\f(x)=0$ and
\newline
$\sum\limits_{n\in\Z}2^{np}
 \left(\int_{2^n}^{\be}|d\f(x)|\right)^{p}<\be;$
\item
$\f$ has locally bounded variation, $\lim\limits_{x\to\be}\f(x)=0$ and
\newline
$
\sum\limits_{n\in\Z}2^{np}
 \left(\int_{2^n}^{2^{n+1}}|d\f(x)|\right)^{p}<\be$. 
\item 
$
\sum\limits_{n\in\Z}2^{np}\norm{\f}_{BV[2^n,2^{n+1}]}^{p}<\be$. 
\item 
$
\sum\limits_{n\in\Z}\norm{x\f(x)}_{BV[2^n,2^{n+1}]}^{p}<\be$. 
\end{enumerate}
\end{thm}
\Pf 
To show that (i) implies (ii), it suffices to consider a
non-increasing $\f\in X_p$; 
it is easily seen that then $\lim_{x\to\be}\f(x)=0$ and $V_\f=\f$,
whence (ii) follows. 

Conversely, suppose that (ii) holds. By considering real and imaginary
parts, we may assume that $\f$ is real. 
Then $\f=V_\f-(V_\f-\f)$, where $V_\f$ and $V_\f-\f$ are
non-increasing functions in $X_p$;  
note that $0\le V_\f-\f\le 2 V_\f$. Consequently (i) holds.

Since $V_\f$ is non-increasing, (ii)$\eq$(iii) follows by \refL{L:yp}. 

Next, (iii)$\eq$(iv) follows easily because
$V_\f(x)=\int_x^\be|d\f(y)|$, and  
(iv)$\eq$(v) is easily verified.

If (iv) holds, then $\norm{\f}_{BV[2^n,2^{n+1}]}=
\int_{2^n}^{2^{n+1}}|d\f(x)|
+\sup_{[2^n,2^{n+1}]}|\f|
\le 2\int_{2^n}^{\be}|d\f(x)|$
and (vi) follows. Conversely, (vi) immediately implies (v).

Finally, for any functions $\f$ and $\psi$ on an interval $I$,
$\norm{\psi\f}_{BV(I)}\le\norm{\psi}_{BV(I)}\norm{\f}_{BV(I)}$,
and the equivalence (vi)$\eq$(vii) follows by taking 
$\psi(x)=x$ and $\psi(x)=1/x$.
\qed

We  can define a norm in $Y_p$ (a quasi-norm for $p<1$) by
\begin{equation}\label{yp}
\|\f\|_{Y_p}\df\biggpar{\sum_{n\in\Z}2^{np}
 \biggpar{\int_{2^n}^{2^{n+1}}|d\f(x)|}^{p}}^{1/p};
\end{equation}
an alternative is 
$$
\left(\int_0^\be|xV_\f(x)|^p\frac{dx}{x}\right)^{1/p},
$$ 
which yields an
equivalent (quasi-)norm.

We obtain as corollaries to Theorems \refand{T:yp}{T:yp1} 
the following simple sufficient conditions for 
$\qf\in \bS_p$.
\begin{cor}
\label{C:yp1}
If $\f$ is absolutely continuous on $(0,\be)$, 
$\lim\limits_{x\to\be}\f(x)=0$ and $\sup_{x>0} x^\gamma|\f'(x)|<\be$
for some $\gamma>2$, then 
$\f\in Y_p$ for every $p>0$ and thus $\qf\in \bS_p$ for every $p>1/2$.
\end{cor}

\Pf
$V_\f(x)$ is bounded and $|V_\f(x)|\le\const\cdot x^{1-\gamma}$, and thus
$xV_\f(x)\in L^p(dx/x)$ for every $p>0$.
$\bl$

\begin{cor}
\label{C:yp2}
If $\f$ has bounded variation and support in a finite interval, then
$\f\in Y_p$ for every $p>0$ and thus $\qf\in \bS_p$ for every $p>1/2$.
\nopf
\end{cor}

\

\setcounter{equation}{0}
\section{\bf $\bs{p=1}$, first results}
\label{S:p1}

\

Let us now consider the case $p=1$. We know already that $\f\in X_1$
is a necessary and $\f\in Y_1$ a sufficient condition for $\qf\in\bS_1$.
We will later see that neither condition is both
necessary and sufficient (see the example following Theorem \ref{E:exp}). 
We restate these results as follows.

\begin{thm}
If $\f$ has locally bounded variation, $\into x|d\f(x)|<\be$ and
\linebreak$\limxoo\f(x)=0$, then $\qf\in\bS_1$.
\end{thm}

\Pf
It follows from Theorem \ref{T:yp1} and the calculation
\begin{equation*}
\into V_\f(x)\,dx=\into\int_x^\be|d\f(y)|\,dx =\into y|d\f(y)|
\end{equation*}
that the assumption is equivalent to $\f\in Y_1$, so the result
follows from Theorem \ref{T:yp}. $\bl$

\begin{thm}\label{T:x1}
If $\qf\in\bS_1$, then $\f\in X_1$. 
Furthermore, $\f\in L^1(0,\be)$ and
\begin{equation*}
\trace\qf=\into\f(x)\,dx.
\end{equation*}
\end{thm}

\Pf
By the Remark at the end of \refS{S:p>1}, $\qf\in\bS_1\imp\f\in X_1$. 
Next, \linebreak$X_1\subset L^1(0,\be)$, since by the \CSineq{} and \eqref{xpi}
\begin{equation*}
\into|\f|\le
\sum_{n\in\Z}2^{n/2}
 \left(\int_{2^n}^{2^{n+1}}|\f(x)|^2dx\right)^{1/2}<\be.
\end{equation*}
Finally, the trace formula follows from \refT{T:trace} below, since
with \linebreak$k(x,y)=\f(\max\{x,y\})$ for $x,y>0$,
\begin{equation*}
\int_{-\be}^\be k(x,x+a)\,dx=\int_{|a|}^\be\f(x)\,dx.
\quad\bl
\end{equation*}

\begin{remark}
This theorem gives a formula for the trace of $Q_{\varphi}$ if that
operator has a trace; i.e., if it is in the trace class, $\bS_{1}.$ The
theorem also shows that if $Q_{\varphi}$ is in the trace class then we must
have $\varphi\in X_{1}.$ We will see later in this section that $\varphi\in
X_{1}$ is not enough to insure that $Q_{\varphi}$ is in the trace class.
However, we will see later, \refC{C:Dix}, that $\varphi\in X_{1}$ is
sufficient to insure that $Q_{\varphi}$ and related operators do have a
\emph{Dixmier trace}.
\end{remark}

In the previous theorem we used following fact from \cite{A},
improving an earlier result in \cite{B}: 
\begin{thm}
\label{T:trace}
If $T$ is an integral 
operator on $L^2(\R)$ of class $\bS_1$ with kernel function $k$, then
the function
$$
x\mapsto k(x,x+a),\quad x\in\R,
$$
is in $L^1(\R)$ for almost all $a\in\R$ and the function
$$
a\mapsto\int_\R k(x,x+a)dx,\quad a\in\R,
$$
is 
almost everywhere equal to
the Fourier transform of a function in $h\in L^1(\R)$,
in particular, it coincides \aex{} with a continuous function on $\R$.
Moreover, $$\trace T=(\F h)(0).$$
\end{thm}

\Pf It is sufficient to prove the result when $k(x,y)=f(x)g(y)$ with
$f$ and $g$ in $L^2(\R)$, 
in which case it can be verified straightforwardly,
with $h(\xi)=\hat f(\xi)\hat g(-\xi)$. 
$\bl$

We can reduce the estimation of $\|Q_\f\|_{\bS_1}$
to the estimation 
of $\|\qf^{I}\|_{\bS_1}$ for dyadic intervals $I$.

\begin{thm}
\label{S1}
Let $\f\in \Lloc^1(\R_+)$ and let $I_n=[2^n,2^{n+1}]$.
Then $Q_\f\in\bS_1$ if and only if $\f\in X_1$ and
\bay
\label{diag}
\sum_{n\in\Z}\|\qf^{I_n}\|_{\bS_1}<\be.
\ey
\end{thm}

\Pf 
Consider the sets $A^{(k)}$ defined by \rf{Ak} and consider
their symmetric images $A^{(-k)}$ about the line 
$\{(x,y):~x=y\}$. As in \eqref{upper} we have
$$
\sum_{k\ge1}\|\cp_{A^{(k)}}Q_\f\|_{\bS_1}
\le\const\sum_{n\in\Z}\left(2^n\int_{2^n}^{2^{n+1}}|\f(x)|^2dx\right)^{1/2}.
$$
Similarly,
$$
\sum_{k\ge1}\|\cp_{A^{(-k)}}Q_\f\|_{\bS_1}
\le\const\sum_{n\in\Z}\left(2^n\int_{2^n}^{2^{n+1}}|\f(x)|^2dx\right)^{1/2}.
$$
It thus follows from $\f\in X_1$ that
$$
\cp_{A}Q_\f\in\bS_1\quad\mbox{and}\quad \cp_{\breve A}Q_\f\in\bS_1.
$$
where
$$
A=\bigcup_{k\ge1}A^{(k)}\quad\mbox{and}\quad \breve
A=\bigcup_{k\le-1} A^{(k)}. 
$$
Consequently, using \refT{T:x1}, $Q_\f\in\bS_1$ 
if and only if $\f\in X_1$ and $\cp_BQ_\f\in\bS_1$,
where $B$ is defined by 
\rf{B}. Since the projections of the sets $B_n$ onto the coordinate
axes are disjoint, 
and $\cp_{B_n}Q_\f=\qf^{I_n}$,
it follows 
that $\cp_BQ_\f\in\bS_1$ if and only if \rf{diag} holds. 
$\bl$

Let $n\in\Z$. It is easy to see that $Q_\f^{I_n}\in\bS_1$ if and
only if $Q_{\f_n}\in\bS_1$, where 
$$
\f_n(x)\df\begin{cases} \f(x+2^n),&x\in[0,2^n]\\
0,&\text{otherwise}.\end{cases}
$$
Hence, the question of when $Q_\f$ belongs to $\bS_1$ reduces to the
question of estimating 
$\|Q_\f\|_{\bS_1}$ for functions $\f$ supported on finite intervals.

\begin{remark}
For $0<p<1$, 
it can similarly be shown that if
$\f\in X_p$ and
$
\sum_{n\in\Z}\|\qf^{I_n}\|_{\bS_p}^p<\be
$,
then $\qf\in\bS_p$.
We do not know whether the converse holds.
\end{remark}

We next show that $\f\in X_1$ is not sufficient for $\qf\in \bS_1$.

\begin{thm}
\label{E:exp}
Let $\f_N(x)=e^{2\pi{\rm i}Nx}\chi_{\OI}(x)$ for $N=1,2,\dots$
Then
\begin{equation*}
s_n(\qq{\f_N}) \asymp \min\left\{\frac 1{n+1},\frac{N}{(n+1)^2}\right\},
\end{equation*}
and so
\begin{equation}\label{eexp1}
\|\qq{\f_N}\|_{\bS_1}\asymp \log(N+1).
\end{equation}
\end{thm}

Note that $\asymp$ means that the ratio of the two sides are bounded from
above and below by positive constants (not depending on $n$ or $N$).
Clearly, 
$$
\|\f_N\|_{X_1}=\|\chi_{\OI}\|_{X_1}=C
$$ 
is independent of $N$. This shows, by the closed graph theorem, that 
$$
\f\in X_1\not\imp\qf\in\bS_1;
$$
a concrete counterexample is given in
the example following the proof of Theorem \ref{E:exp}. Moreover,
$\f_N\in X_p$ for any $p>0$, again 
with norm independent of $N$, so for every $p\le1$
$$
\f\in X_p\not\imp\qf\in\bS_1\quad\mbox{and}\quad
\f\in X_p\not\imp\qf\in\bS_p.
$$

It also follows from Theorem \ref{E:exp} that 
\begin{equation}\label{eexpp}
\|\qq{\f_N}\|_{\bS_p}\asymp N^{(1-p)/p}, 
\qquad \tfrac12 < p <1.
\end{equation}

{\bf Proof of Theorem \ref{E:exp}.}
Since multiplication by a unimodular
 function is a unitary operator, the
singular values of $\qq{\f_N}$ are the same as the singular values
$s_n(T_N)$, where $T_N$ is the integral operator on $L^2[0,1]$ with kernel
\begin{equation*}
\begin{split}
e^{-\pi \ii Nx} \f_N\bigpar{\max(x,y)} e^{-\pi\ii Ny}
&=\exp\bigpar{2\pi\ii N\max(x,y)-\pi\ii Nx-\pi\ii Ny}\\
&=\exp\bigpar{\pi\ii N|x-y|}.
\end{split}
\end{equation*}

Let $g_N(x)=e^{\pi\ii N|x|}$ for $|x|\le1$, and extend $g_N$ to a
function on $\R$ with period 2.
Let $T_N'$ be the integral operator on $L^2[-1,1]$ with kernel
$g_N(x-y)$.
If $I_+=[0,1]$, $I_-=[-1,0]$ and $A_{\ga\gb}=I_\ga\times I_\gb$,
$\ga,\gb\in\set{+,-}$, then $\cp_{A_{++}} T_N'=T_N$ and thus
\begin{equation}
\label{sntn}
s_n(\qq{\f_N}) = s_n(T_N)\le s_n(T_N').
\end{equation}
Moreover, each $\cp_{A_{\ga\gb}} T_N'$ is by a translation unitarily
equivalent to either $T_N$ or the integral operator on $L^2[0,1]$ with
kernel
$
g_N(x-y-1)=\exp\bigpar{\pi\ii N(1-|x-y|)} = (-1)^N \overline{
g_N}(x-y)$,
which has the same singular values.
Hence, by \eqref{snsum},
\begin{equation}
\label{sntn4}
s_{4n}(T_N')\le 4 s_n(T_N) =4 s_n(\qq{\f_N}).
\end{equation}

$T_N'$ is a convolution operator on the circle $\R/2\Z$, so the 
elements of the orthonormal basis 
$\set{2^{-1/2}\exp(\pi\ii kx)}_{k=-\infty}^\be$ in $L^2[-1,1]$ 
are eigenvectors with
eigenvalues
\begin{equation*}
\begin{split}
\hat g_N(k)&=\int_{-1}^1 g_N(x)e^{-\pi\ii kx}dx
=\int_{-1}^1 e^{\pi\ii N|x|-\pi\ii kx}\,dx\\
&=\int_{-1}^0 e^{-\pi\ii(N+k)x}\,dx +\int_{0}^1 e^{\pi\ii(N-k)x}\,dx\\
&=
\begin{cases}
0, & k\equiv N \pmod 2,\, k\neq\pm N,\\
1, & k=\pm N,\\
\frac{-2}{\pi\ii(N+k)} + \frac{-2}{\pi\ii(N-k)} 
= \frac{4\ii N}{\pi(N^2-k^2)}, & k\equiv N+1\pmod2.
\end{cases}
\end{split}
\end{equation*}
Consequently, the singular values $s_n(T_N)$ are the absolute values
$|\hat g_N(k)|$, $k\in\Z$, arranged in decreasing order. This easily
yields
$$
s_n(T_N')\asymp \min\left\{\frac 1{n+1},\frac{N}{(n+1)^2}\right\},
$$
and the result follows by \eqref{sntn} and \eqref{sntn4}.
\qed

\begin{remark}
A related method is used in a more general context in \refS{S:matrix}.
Indeed, the estimates \eqref{eexp1} and \eqref{eexpp} follow easily
from \refT{wHT} (with the norm estimates implicit there).

Moreover, \eqref{eexp1} and \eqref{eexpp} also follow
from the results in \refS{S:dilation}, obtained by a different 
method.
\end{remark}

\begin{example}
For a concrete counterexample we let $N_k\ge2$ be integers and define
\begin{equation*}
\f=\sum_{k=1}^\be \exp(2\pi{\rm i}2^kN_kx)\chi_{(2^{-k},2^{1-k})}(x).
\end{equation*}
Then $|\f|=\chi_{\OI}$; in particular, $\f\in X_p$ for every $p>0$.

On the other hand, for every $k$, by \refL{L:homo} and \refT{E:exp},
\begin{equation*}
\|{\qf}\|_{\bS_1}
\ge
\left\|\qf^{[2^{-k},2^{1-k}]}\right\|_{\bS_1}
=\|{\qq{\f_{N_k}}}\|_{\bS_1}
\ge c2^{-k}\log N_k.
\end{equation*}
Choosing $N_k=2^{3^k}$, we find that $\qf\notin\bS_1$. 

Other choices yield further interesting examples.
Thus, $N_k=2^{k}$ yields a symbol $\f\in X_1$ but $\f\notin Y_1$
such that, by \refT{S1}, $\qf\in\bS_1$.
In fact, using the Remark followed by Theorem \ref{E:exp}, we can
conclude that  
$\qf\in\bS_p$ for every $p>1/2$.

The choice $N_k=2^{k^2}$ yields a symbol $\f\in X_1\setminus Y_1$
such that $\qf\in\bS_1$ but $\qf\notin\bS_p$ for $p<1$.
\end{example}

\begin{remark}
\refT{E:exp} implies also that
$\f\in X_1$ does not imply $\qf\in\bS_{1,q}$ 
for any
Schatten--Lorentz space $\bS_{1,q}$ with $q<\be$.
\end{remark}

Let us prove 
now that the condition $Q_\f\in\bS_1$ does not imply that $Q_\f^+\in\bS_1$. 
Moreover, as previously shown by Nowak \cite{No}, 
we show that there are no non-zero operators
$Q_\f^+$ of class $\bS_1$.
(A more refined result will be given in \refT{T:ws1}.)

\begin{thm}
\label{zero}
Suppose that $\f\in \Lloc^1(\R_+)$. If $Q_\f^+\in\bS_1$, then
$\f$ is the zero function. 
\end{thm}

\Pf
Suppose that $Q_\f^+\in\bS_1$. 
Let $k$ be the kernel function of $Q_\f^+$, extended by 0 to $\R^2$, \ie
$$
k(x,y)= \begin{cases}\f(x),&0<y<x\\0,&\text{otherwise}.\end{cases}
$$
Let $\D\subset\R_+$ be a compact interval. 
Consider the operator $P_\D Q_\f^+P_\D$, where $P_\D$ is
multiplication by $\D$. Clearly,  
$P_\D Q_\f^+P_\D$ is an integral operator with kernel function
$k_\D\df\chi_{\D\times\D}k$, (recall that 
$\chi_A$ is the characteristic function of a set $A$) and $P_\D
Q_\f^+P_\D\in\bS_1$. It follows from 
\refT{T:trace}
that the function $u$ 
$$
u(a)\df\int_\R k_\D(x,x+a)\,dx
$$
is a.e.\ equal to a continuous function $\R$. Clearly, $u(a)=0$ if $a>0$. And 
$$
u(a)\to\int_{\D} \f(x)\,dx
$$
when $a\upto0$. Hence, 
$$
\int_\D\f(x)\,dx=0,\quad\text{for any interval $\D\subset\R_+$}.
$$
Consequently, $\f=0$. 
$\bl$

\

\setcounter{equation}{0}
\section{\bf Schur multipliers of the form $\bs{\psi(\max\{x,y\})}$,
$\bs{x,y\in\R_+}$} 
\label{Schur}

\

Let $0<p\le2$. Recall that a function $\o$ on $\R^2$ is called a {\it
Schur multiplier of} 
$\bS_p$ if the integral operator on $L^2(\R)$ with kernel function 
$\o k$  belongs to $\bS_p$ 
whenever the integral operator with kernel function $k$ does. If
$2<p<\be$, the class of  
Schur multipliers of $\bS_p$ can be defined by duality: $\o$ is a
Schur multiplier of $\bS_p$ 
if $\o$ is a Schur multiplier of $\bS_{p'}$, where $p'\df p/(p-1)$. 
We say that $\o$ is a \emph{Schur multiplier of weak type} $(p,p)$,
$0<p\le2$, if the integral operator 
with kernel function $\o k$ belongs to $\bS_{p,\be}$ whenever the
integral operator with kernel function 
$k$ belongs to $\bS_p$. Note that in a similar way 
one can define Schur multipliers for an arbitrary 
measure space $(\X,\mu)$.

In this section for a function $\psi\in L^\infty(\R)$ we find a
sufficient condition   
for the function $(x,y)\mapsto\psi(\max\{x,y\})$, $(x,y)\in\R^2$, to
be a Schur multiplier of $\bS_p$. 
We also obtain a sufficient condition for this function to be a Schur
multiplier of weak type 
$(1/2,1/2)$.

\begin{thm}
\label{mp>1}
Let $1<p<\be$ and let $\psi\in L^\be(\R)$. Then the function 
$$
(x,y)\mapsto\psi(\max\{x,y\}),\quad (x,y)\in\R^2,
$$
is a Schur multiplier of $\bS_p$.
\end{thm}

\Pf Since the triangular projection $\cp$ is bounded on $\bS_p$, the
characteristic 
function of the set \mbox{$\{(x,y):~x>y\}$} is a Schur multiplier of
$\bS_p$. It remains 
to observe that
\bay
\label{khar}
\psi(\max\{x,y\})=\psi(x)\chi_{\{(x,y):x>y\}}+\psi(y)\chi_{\{(x,y):x<y\}}.
\quad\bl 
\ey

\begin{thm}
\label{slt1}
Let $\psi\in L^\be(\R)$. Then the function 
$$
(x,y)\mapsto\psi(\max\{x,y\}),\quad (x,y)\in\R^2,
$$
is a Schur multiplier of weak type $(1,1)$.
\end{thm}

\Pf The result follows from \eqref{khar} and the fact that
the triangular projection $\cp$ has weak type $(1,1)$ (see \rf{P11}). $\bl$

\begin{thm} 
\label{mlp}
Let $1/2<p<\be$ and let $\psi$ be a function of bounded variation.
Then the function $(x,y)\mapsto\psi(\max\{x,y\})$ on $\R^2$ is a Schur
multiplier of $\bS_p$. 
\end{thm}

\Pf By Theorem \ref{mp>1} we may assume that $p\le1$. 

We consider first the case when $\psi$ is absolutely continuous, i.e., 
\bay
\label{inth}
\psi(x)=\int\limits_x^{\infty}h(t)\,dt+C, \quad h\in L^1(\R).
\ey
We may assume that $C=0$. 

Let $\xi$ and $\eta$ be functions  in $L^2$ and let $T$ be the integral
operator defined by 
\bay
\label{xiepsi}
(Tf)(x)=\int_\R \xi(x)\eta(y)\psi(\max\{x,y\})f(y)\,dy,\quad f\in L^2(\R).
\ey
We have to prove that
$$
\|T\|_{\bS_p}\le C(p,h)\|\xi\|_{L^2(\R)}\|\eta\|_{L^2(\R)},
$$
where $C(p,h)$ may depend only on $p$ and $h$. 

We can factorize the function $h$ in the form
$h=uv$, where $u,\,v\in L^2(\R)$. Put

$$
k_1(x,y)\df\left\{\begin{array}{ll}0,&y<x,\\
\xi(x)u(y),&y>x,
\end{array}\right. 
$$
and
$$
k_2(x,y)\df\left\{\begin{array}{ll}\eta(y)v(x),&y<x,\\
0,& y>x.
\end{array}\right.
$$
Let $T_1$ and $T_2$ be the integral operators on $L^2(\R)$ with kernel
functions $k_1$ and $k_2$. It follow from the boundedness of the
triangular projection that if $1<q<\be$, 
then $\|T_1\|_{\bS_q}\le C(q)\|\xi\|_{L^2}\|u\|_{L^2}$ and
$\|T_2\|_{\bS_q}\le C(q)\|\eta\|_{L^2}\|v\|_{L^2}$, where $C(q)$ may
depend only on $q$. 
It is also easy to verify that $T=T_1T_2$. It follows that
$$
\|T\|_{\bS_p}\le\|T_1\|_{\bS_{2p}}\|T_2\|_{\bS_{2p}}
\le(C(2p))^2\|\xi\|_{L^2}\|\eta\|_{L^2}\|u\|_{L^2}\|v\|_{L^2}.
$$

To reduce the general case to the case of an absolutely continuous
function $\psi$, we can 
consider a standard regularization process. $\bl$

We complete this section with the following result.

\begin{thm}
\label{sl1/2}
Suppose that $\psi$ is a function of 
bounded variation. 
Then the function 
$(x,y)\mapsto\psi(\max\{x,y\})$ is a Schur multiplier of weak type $(1/2,1/2)$.
\end{thm}

We need two lemmata.

\begin{lem} 
\label{pc}
Let $0<p<1$ and
let $A\in\bS_{p,\infty}$. 
Set 
$$
\|A\|_{\bS_{p,\infty}}^*
\df\sup_{t>0}\left(t^{p-1}\sum\limits_{n\ge0}\min\{t,s_n(A)\}\right)^{1/p}.
$$ 
Then 
$$
\|A\|_{\bS_{p,\infty}}\le\|A\|_{\bS_{p,\infty}}^*
\le(1-p)^{-1/p}\|A\|_{\bS_{p,\infty}}.
$$
\end{lem}

\Pf Taking $t=s_n(A)$ in the definition of
$\|\cdot\|_{\bS_{p,\infty}}^*$, we obtain 
$$
(n+1)^{\frac1p}s_n(A)\le\|A\|_{\bS_{p,\infty}}^*.
$$
Consequently,
$\|A\|_{\bS_{p,\infty}}\le\|A\|_{\bS_{p,\infty}}^*$. Next, we have 
\bey
t^{p-1}\sum_{n\ge0}\min\{t,s_n(A)\}&\le&t^{p-1}\sum_{n\ge0}\min\{t,\|A\|_
{\bS_{p,\infty}}(n+1)^{-\frac1p}\}\\[.2cm]
&\le&
t^{p-1}\int\limits_0^\infty \min\{t,\|A\|_{\bS_{p,\infty}}x^{-1/p}\}\,dx=
(1-p)^{-1}\|A\|_{\bS_{p,\infty}}^p.
\eey
Hence, $\|A\|_{\bS_{p,\infty}}^*\le(1-p)^{-1/p}\|A\|_{\bS_{p,\infty}}$.
$\bl$

\begin{lem}
\label{pnorm}
If $0<p<1$, then $\|\cdot\|_{\bS_{p,\be}}^*$ is a $p$-norm, i.e.,
$$
\|A_1+A_2\|_{\bS_{p,\be}}^{*p}
\le\|A_1\|_{\bS_{p,\be}}^{*p}+\|A_2\|_{\bS_{p,\be}}^{*p},
\quad A_1,\,A_2\in\bS_{p,\be}.
$$
\end{lem}

\Pf By Rotfeld's theorem \cite{R}, if $\Phi$ is a concave
nondecreasing function on $\R_+$ such that 
$\Phi(0^+)=0$, then
\bay
\label{rot} 
\sum_{j=0}^m\Phi(s_j(A_1+A_2))
\le\sum_{j=0}^m\Phi(s_j(A_1))+\sum_{j=0}^m\Phi(s_j(A_2)),\quad m\in\Z_+.
\ey
For $t>0$ we define the function $\Phi_t$ on $\R_+$ by
$$
\Phi_t(x)=t^{p-1}\min\{t,x\}.
$$
Clearly,
$$
\|A\|_{\bS_{p,\be}}^{*p}=\sup_{t>0}\sum_{j\ge0}\Phi_t(s_j(A)).
$$
It remains to apply \ref{rot} for $\Phi_t$ and take the {\it supremum}
over $t>0$. $\bl$ 

Note that the fact that the space $L^{p,\be}$ has a $p$-norm is well
known (see \cite{K}).

{\bf Proof of Theorem \ref{sl1/2}.} The proof is similar to the proof
of Theorem \ref{mlp}. Again, it is 
sufficient to assume that $\psi$ has the form \rf{inth} with $C=0$. 

By Lemma \ref{pc} and Lemma \ref{pnorm},
to prove that our function is a Schur multiplier of weak type
$(1/2,1/2)$, it is sufficient 
to prove that if $T$ is defined by
\rf{xiepsi}, then $T\in\bS_{1/2,\be}$. Let $u$, $v$, $k_1$, $k_2$,
$T_1$, and $T_2$ be as in the proof of 
Theorem \ref{mlp}. Let $n\ge2$ and $n=m_1+m_2$, where
$|m_1-n/2|\le1/2$ and $|m_2-n/2|\le1/2$. 
Since the triangular projection $\cp$ has weak type $(1,1)$ (see
\rf{P11}), we have 
$$
\|T_1\|_{\bS_{1,\be}}\le\const\|\xi\|_{L^2}\|u\|_{L^2}\quad\mbox{and}
\quad\|T_2\|_{\bS_{1,\be}}\le\const\|\eta\|_{L^2}\|v\|_{L^2}.
$$
Hence, by \rf{snpr},
\bey
s_n(T)&\le&s_{m_1}(T_1)s_{m_2}(T_2)
\le\const\frac1{m_1m_2}\|\xi\|_{L^2}\|u\|_{L^2}\|\eta\|_{L^2}\|v\|_{L^2}
\\[.2cm]
&\le&\const\frac{1}{n^2}\|\xi\|_{L^2}\|u\|_{L^2}\|\eta\|_{L^2}\|v\|_{L^2}
\eey
which completes the proof. $\bl$

We let $\Mp$ denote the space of 
Schur multipliers of $\bS_p$, and put
$$
\|\omega\|_{\Mp}\df\sup\|\omega k\|_{\bS_p},
$$
where the {\it supremum} is taken over
all integral operators with kernel $k\in L^2$ such that
$\|k\|_{\bS_p}=1$. Here by 
$\|k\|_{\bS_p}$ we mean the $S_p$ norm (quasi-norm if $p<1$) of the integral
operator with kernel function $k$.
If $\omega$ is a Schur multiplier of weak type $(p,p)$, we put
$$
\|\omega\|_{{\frak M}_{p,{\rm w}}}\df\sup\|\omega k\|_{\bS_{p,\infty}},
$$ 
where the {\it supremum} is taken
over all integral operators with kernel $k\in L^2$ such that $\|k\|_{\bS_p}=1$.

\begin{remark}
It is clear from the proofs of Theorems \ref{mp>1},
\ref{slt1}, \ref{mlp}, and \ref{sl1/2} 
that under the hypotheses of Theorem \ref{mp>1} we have 
$$
\|\psi(\max\{x,y\})\|_{\Mp}\le C(p)\|\psi\|_{L^\infty}
$$ 
and  
$$
\|\psi(\max\{x,y\})\|_{{\frak M}_{1,{\rm w}}}\le C\|\psi\|_{L^\infty}
$$
while under the hypotheses of Theorem \ref{mlp} we have 
$$
\|\psi(\max\{x,y\})\|_{\Mp}\le C(p)\|\psi\|_{BV},
$$
and
$$
\|\psi(\max\{x,y\})\|_{{\frak M}_{1/2,{\rm w}}}\le C\|\psi\|_{BV}.
$$
Here we use the notation
$$
\|\f\|'_{BV}
=\int_\R|d\f|\quad\mbox{and}\quad \|\f\|_{BV}=\|\f\|'_{BV}+\|\f\|_{L^\be}.
$$
\end{remark}

The following result gives us a more accurate estimate for
$\|\psi(\max\{x,y\})\|_{\Mp}$ 
in the case $1/2<p\le1$.

\begin{thm}
\label{mpe}
Let $\psi$ be a function of bounded variation on $\R$. Then
$$
\|\psi(\max\{x,y\})\|_{{\frak M}_1}\le \const
\|\psi\|_{L^\infty}\log
 \left(2+\frac{\|\psi\|^\prime_{BV}}{\|\psi\|_{L^\infty}}\right) 
$$
and
$$
\|\psi(\max\{x,y\})\|_{\Mp}\le 
\const\|\psi\|_{L^\infty}^{2-1/p}\|\psi\|_{BV}^{1/p-1},\quad\frac12<p<1.
$$
\end{thm}

\Pf Let $\xi$ and $\eta$ be function in $L^2$ such that
$\|\xi\|_{L^2}=\|\eta\|_{L^2}=1$.  
We have to estimate the $\bS_p$-norm of the integral operator with kernel
$\psi(\max\{x,y\})\xi(x)\eta(y)$. Let $\{s_n\}_{n\ge0}$ be the
sequence of $s$-numbers of this integral 
operator.  Theorem \ref{slt1} implies that 
$$
s_n\le\const\frac{\|\psi\|_{L^\be}}{n+1}.
$$
Theorem \ref{sl1/2} implies that
$$
s_n\le\const\frac{\|\psi\|_{L^\infty}+\|\psi\|^\prime_{BV}}{(n+1)^2}.
$$ 
Consequently,
$$
s_n\le\const\min\left\{\frac{\|\psi\|_{L^\be}}{n+1},
 \frac{\|\psi\|_{L^\infty}+\|\psi\|'_{BV}}{(n+1)^2}\right\}.
$$
The rest of the proof is an easy exercise.
$\bl$

\

\setcounter{equation}{0}
\section{\bf The case $\bs{p=1/2}$}
\label{S:1/2}

\

\begin{thm} 
\label{wt}
Let $\f$ be a function of bounded variation on $[0,1]$.
Then \linebreak$Q_\f^{[0,1]}\in\bS_{\xfrac12,\,\infty}$.
\end{thm}

\Pf We may extend the function $\f$ by putting $\psi(t)=0$ for
$t\in\R_+\setminus[0,1]$. Clearly,  
the integral operator with kernel function $\chi_{[0,1]^2}$ has rank
one, and so it  
belongs to $\bS_{1/2}$. Consequently, by Theorem \ref{sl1/2}, the
integral operator with kernel function 
$\chi_{[0,1]^2}\f(\max\{x,y\})$ belongs to
$\bS_{1/2,\,\be}$.
$\bl$

To see that this result cannot be improved to   
$Q_\f^{[0,1]}\in\bS_{\xfrac12}$,
we begin with two extensions of \refT{zero}.

\begin{lem}
\label{sone}
Let $\f,\psi\in \Lloc^1(\R)$. Suppose that the integral operator with
kernel function $k$, 
$$
k(x,y)= \begin{cases}\f(x)\psi(y),&y\le x,\\0,&\text{otherwise},\end{cases}
$$
belongs to $\bS_1$.
Then $\f\psi=0$ almost everywhere. 
\end{lem}

\Pf 
First we assume that $\f,\psi\in L^2(\R)$.  By \refT{T:trace} we have
$$
\lim\limits_{a\to0-}\int\limits_{\R}k(x,x+a)dx
=\lim\limits_{a\to0+}\int\limits_{\R}k(x,x+a)dx=0, 
$$
whence $\int\limits_{\R}\f(x)\psi(x)\,dx=0$. Now let $\f$ and $\psi$
be arbitrary functions in $\Lloc^1(\R)$. 
Suppose that $f$ and $g$ are functions in $L^\infty(\R)$ such that
\mbox{$f\f\in L^2(\R)$} and  
\linebreak$g\psi\in L^2(\R)$. Consider the integral operator with
kernel function 
$$
(x,y)\mapsto f(x)k(x,y)g(y).
$$
Clearly, it belongs to $\bS_1$. It follows from what we have just proved that
\linebreak$\int\limits_{\R}\f(x)f(x)\psi(x)g(x)dx=0$. Since $f$ and
$g$ are arbitrary, this implies the result. 
$\bl$

\begin{lem}
\label{soone}
Let $\f\in \Lloc^1(\R_+)$ and let $A=(0,\infty)\times\Delta$, where
$\Delta$ is a measurable subset of $(0,\infty)$. Suppose that
$\cp_AQ^+_\f\in\bS_1$. 
Then $\f=0$ almost everywhere on $\Delta$. 
\end{lem}

\Pf The result follows easily from 
\refL{sone} 
 with $\psi=\chi_\Delta$.
$\bl$

\begin{lem}
\label{lmon}
Let $\f$ be  a nonincreasing locally absolutely continuous function on
$\R_+$ 
and let $\D$ be a measurable subset
of $\R_+$. Suppose that the integral 
operator with kernel function
$$
(x,y)\mapsto\f(\max\{x,y\})\chi_{\D}(x)\chi_{\D}(y)
$$ 
belongs to $\bS_{\xfrac12}$. Then $\f^\prime=0$ almost everywhere on $\D$.
\end{lem}

\Pf 
By replacing $\D$ with $\D\cap(a,b)$, we may assume that
$\D\subset[a,b]$ where $0<a<b<\be$. We may then subtract $\f(b)$ and
modify $\f$ outside $[a,b]$ so that $\f$ becomes constant on $(0,a]$
and zero on $[b,\be)$.
Let $\psi=(-\f')^{1/2}$.
Since $x\f(x)$ is bounded, we have $\psi\in X_\be$.
Thus $Q_\psi^+$ is bounded and by
Theorem \ref{fact}, 
$Q_\f$ admits a factorization $Q_\f=(Q^+_\psi)^*Q_\psi^+$.
Let $M$ be
multiplication by $\chi_\D$. It follows that  
$MQ_\f M=(Q^+_\psi M)^*(Q_\psi^+M)$, and so
$Q_\psi^+M\in\bS_1$. The result follows now from Lemma 
\ref{soone}. $\bl$

\begin{cor}
\label{locone}
Let $\f$ be as in Lemma \ref{lmon} and let $\psi$ be a function on
$\R_+$ such that 
$Q_\psi\in\bS_{\xfrac12}$. Set $\D\df\{x\in\R_+:~\psi(x)=\f(x)\}$.
Suppose that $\psi$ is differentiable almost everywhere on $\D$.
Then 
$\f'=\psi^\prime=0$ almost everywhere on $\Delta$.
\end{cor}

\Pf It suffices to apply Lemma \ref{lmon}.
$\bl$

\begin{thm}
\label{oc}
Let $\psi$ be an absolutely continuous function on $[0,1]$. Suppose
$Q_\psi^{[0,1]}\in\bS_{1/2}$. 
Then $\psi$ is constant.
\end{thm}

\Pf Clearly, we may assume that $\psi$ is a real function. Suppose
that $\psi$ is not constant. 
Then $\max \psi>\psi(1)$ or $\min \psi<\psi(1)$. To be definite, suppose that
$\max\psi>\psi(1)$. We use the ``sun rising method''.
Let 
$$
\D\df\{x\in[0,1]:~\psi(x)\ge \psi(t)\quad\text{for all}\quad t\ge x\}.
$$
Clearly,
$\Delta$ is closed and $1\in\D$. Moreover, the restriction $\psi\big|[\a,\b]$
is constant for any interval $(\a,\b)$ such that $\a,\b\in\Delta$ and
$\Delta\cap(\a,\b)=\varnothing$. 
Set 
$$
\f(x)=\max\limits_{[x,1]\cap\Delta}\psi.
$$ 
Clearly, $\f$ is decreasing, $\f|_\Delta=\psi|_\Delta$ and $\f\big|[\a,\b]$
is constant for any interval $(\a,\b)$ such that $\a,\b\in\Delta$ and
$\Delta\cap(\a,\b)=\varnothing$. 
Consequently, $\f$ is absolutely continuous. Thus, $\f^\prime=0$
almost everywhere on $\Delta$  
by Corollary \ref{locone}. Moreover, $\f^\prime=0$ outside $\Delta$
because $\f$ is locally constant 
outside $\Delta$. Consequently, $\f^\prime=0$ almost everywhere on
$[0,1]$, and so $\f(t)=\f(1)=\psi(1)$ 
for any $t\in[0,1]$ which contradicts the condition $\max \psi>\psi(1)$.
$\bl$

\begin{cor}
\label{oc1} 
Suppose that $Q_\psi\in\bS_{\xfrac12}$. Then $\psi$ is constant on any
interval $I\subset\R_+$ on which 
$\psi$ is absolutely continuous.
\nopf
\end{cor}

\begin{cor}
\label{oc2} 
Suppose that $\psi$ is locally  absolutely continuous and
$Q_\psi\in\bS_{\xfrac12}$. Then $\psi=0$ 
everywhere on $\R_+$.
\nopf
\end{cor}

\begin{lem}
\label{lmonm}
Let $\f$ be a nonincreasing function on $\R_+$ with 
$\lim\limits_{t\to+\infty}\f(t)=0$ and let $\e>0$. 
Then there exists a nonincreasing absolutely continuous function
$\psi$ on $\R_+$ 
such that $\lim\limits_{t\to+\infty}\psi(t)=0$ and $\m\{\f\not=\psi\}<\e$.
\end{lem}

\Pf We may assume that $\f$ is right-continuous on $(0,\be)$
and that $\f$ is bounded. 
Consider
the positive measure $\mu$ on $\R_+$ such that $\f(t)=\mu(t,\infty)$
for any $t>0$. 
Denote by $\mu_s$ the singular part of $\mu$. There exists a Borel set
$E\subset\R_+$ 
such that $\m(E)=0$ and $\mu_s((0,+\infty)\setminus E)=0$. We may find
an open set $U$ 
such that $E\subset U\subset\R_+$ and $\mu(U)<\e$. Let
$U=\bigcup\limits_{n\ge1}(a_n,b_n)$, 
where $(a_n,b_n)$ are mutually disjoint. Set
\begin{equation*}
f(t)=\begin{cases}
-\f^\prime(t),&t\in\R_+\setminus E,\\[.2cm]
\dfrac{\mu[a_n,b_n)}{b_n-a_n},&t\in (a_n,b_n).
\end{cases}
\end{equation*}
Set $\psi(t)\df\int\limits_t^\be f(s)ds$. Clearly, $\f=\psi$ outside $U$.
$\bl$

\begin{lem}
\label{lmon1}
Let $\f$ a nonincreasing function on $\R_+$ with
$\lim\limits_{t\to+\infty}\f(t)=0$. Suppose that the integral operator
with kernel function 
$$
(x,y)\mapsto\f(\max\{x,y\})\chi_{\Delta}(x)\chi_{\Delta}(y)
$$ 
belongs to $\bS_{\xfrac12}$ 
for a measurable subset $\D$ of $\R_+$. Then $\f^\prime=0$ almost
everywhere on $\Delta$. 
\end{lem}

\Pf The result follows from  Lemmas \ref{lmon} and \ref{lmonm}.
$\bl$

\begin{thm}
\label{moc}
Let $\psi$ be a function with bounded variation on $[0,1]$. Suppose that
$Q_\psi^{[0,1]}\in\bS_{1/2}$.
Then $\psi^\prime=0$ almost everywhere on $[0,1]$.
\end{thm}

\Pf Again, we may assume that $\psi$ is real. We may also make the
assumption that   
that $\psi$ is continuous at $0$ and at $1$, and
$\psi(t)=\max\{\psi(t^-),\psi(t^+)\}$ for any $t\in(0,1)$. With any
nondegenerate closed  
interval $I\subset[0,1]$ we associate the function $\f_I:I\to\R$ defined by
$\f_I(x)=\sup\{\psi(t):t\in I\quad\text{and}\quad t\ge x\}$. Set
$$
\D(I)\df\{x\in I:~\f_I(x)=\psi(x)\}.
$$
Clearly, $\Delta(I)$ is closed.  By Lemma \ref{lmon1}, 
$\psi^\prime=\f_I'=0$
almost everywhere 
on $\Delta(I)$. 

Set $E_-=\{x\in(0,1):~\psi^\prime(x)<0\}$. Let $a\in
E_-$. Clearly, 
$a\in\Delta(I)$ if $I$ is small enough and $a\in I$.
Consequently,
$$
E_-\subset\bigcup\limits_{n=2}^\be\left(\bigcup\limits_{k=1}^{n-1}
\D\left(\left[\frac{k-1}n,\frac{k+1}n\right]\right)\right).
$$
We have shown that \refL{lmon1} implies $\m(\Delta(I)\cap E_-)=0$ for
every $I$.
Consequently,
$\m(E_-)=0$. 
Thus, we have proved that $\psi^\prime\ge0$ almost everywhere. 
It remains to apply
this result to $-\psi$. 
$\bl$

The following fact is an immediate consequence of Theorem \ref{moc}.
\begin{cor}
\label{moc1} 
Suppose that $Q_\psi\in\bS_{\xfrac12}$. Then $\psi^\prime=0$ almost
everywhere on any interval $I\subset\R_+$ on which 
$\psi$ is of bounded variation.
$\bl$
\end{cor}

\

\setcounter{equation}{0}
\section{\bf Sturm--Liouville theory and $\bs{p=1/2}$}
\label{S:Sturm}

\

If $\qf$ is real, then 
$\qf$ is self-adjoint, so its singular values are the absolute
values of its eigenvalues.
Hence we next study the eigenvalues and eigenfunctions.
For simplicity we consider only the case of symbols $\f$ which vanish
on $(1,\be)$; thus it does not matter whether we consider $\qf$ on
$L^2(\R_+)$ or $\qf^{\OI}$ on $L^2\OI$.

Suppose that $\gf\in C^1\OI$ and that $\gf$ is real.
Let $\gl$ be a non-zero eigenvalue of $\qf^{\OI}$ and $g\in L^2\OI$ a
corresponding eigenfunction,
\begin{equation}
\label{a1}
\gl g(x)=\gf(x)\int_0^x g(y)\,dy + \int_x^1\gf(y)g(y)\,dy,
\qquad 0\le x\le 1.
\end{equation}  
The right hand side is a continuous function of $x$; hence $g\in C\OI$
and \eqref{a1} holds for every $x$ (and not just a.e.). By \eqref{a1}
again, $g\in C^1\OI$ with
\begin{equation}
\label{a2}
\gl g'(x)=\gf'(x)\int_0^x g(y)\,dy.
\end{equation}
Define $G(x)=\int_0^x g(y)\,dy$. Then \eqref{a2} can be written as the
system
\begin{equation}
\label{a3}
\begin{aligned}
G'(x)&=g(x)\\
g'(x)&=\gl^{-1}\gf'(x) G(x)
\end{aligned}
\end{equation}
and we have, using \eqref{a1} with $x=1$, the boundary conditions
\begin{equation}
\label{a4}
\begin{aligned}
G(0)&=0,\\
g(1)&=\gl^{-1}\gf(1) G(1).
\end{aligned}
\end{equation}
Conversely, any solution of \eqref{a3} with the boundary conditions
\eqref{a4} satisfies \eqref{a2} and \eqref{a1}, so the problem of
finding the singular values of $\qf{}$ reduces to finding the $\gl\neq0$
for which \eqref{a3} and \eqref{a4} have a solution.
Note that \eqref{a3} can be written as a Sturm--Liouville problem
\begin{equation}\label{sturm}
\gl G''(x)=\gf'(x) G(x).
\end{equation}

If $g(x_0)=G(x_0)=0$ for some $x_0\in\OI$, then \eqref{a3} shows, 
by the standard uniqueness theorem, that $g$ vanishes identically, a
contradiction. In particular, since $G(0)=0$, we have $g(0)\neq0$, and
we may normalize the eigenfunction $g$ by $g(0)=1$.

For every $\gl\neq0$, $\eqref{a3}$ has a unique solution
$(G_\gl,g_\gl)$ with $G_\gl(0)=0$, $g_\gl(0)=1$. It thus follows that
all non-zero eigenvalues of $\qf$ are simple, and that $\gl\neq0$ is an
eigenvalue if and only if
\begin{equation}
\label{a5}
g_\gl(1)=\gl^{-1}\gf(1) G_\gl(1).
\end{equation}

\begin{example}
Let $\gf(x)=1-x$, $x\in\OI$, and $\f(x)=0$, $x>1$.
Then \eqref{a3} gives $G''(x)=-\gl^{-1} G(x)$, and we find the
solutions
$g_\gl(x)=\cos \gl^{-1/2}x$, $\gl>0$,
and $g_\gl(x)=\cosh |\gl|^{-1/2}x$, $\gl<0$.

Since $\gf(1)=0$, condition \eqref{a5} is simply $g_\l(1)=0$,
and the non-zero eigenvalues are given by $\cos \gl^{-1/2}=0$ or
$\gl^{-1/2}=(n+\frac12)\pi$, $n=0,1,\dots$.
($\gl<0$ is impossible in this case; 
in other words, $\qf{}$ is a positive operator, as is also seen by
\refT{T:positive}.)
Hence, the non-zero eigenvalues are
$\set{(n+\frac12)^{-2}\pi^{-2}}_{n=0}^\infty$ and the singular values
are
$s_n=\pi^{-2}(n+\frac12)^{-2}$, $n\ge0$.
\end{example}

The behaviour $s_n(\qf)\asymp (n+1)^{-2}$ found in the above example
holds for all smooth $\f$ on \OI{} 
by Sturm--Liouville theory, as will be seen 
in \refT{T:lim}.
Hence,
for smooth $\f$ with compact support,
we have $\qf\in \wS{1/2}$ but nothing better.

Let $BV\OI$ denote the Banach space of functions on \OI{} with
bounded variation, 
with the seminorm 
$\|\gf\|'_{BV}=\int_0^1|d\gf|$ and the norm 
$\|\gf\|_{BV}=\|\gf\|'_{BV}+\sup|\f|$.

\begin{thm}
\label{T:BV}
If $\gf\in BV\OI$, and $\f=0$ on $(1,\infty)$, then $\qf\in \wS{1/2}$ and 
$\|\qf\|_{\wS{1/2}}\le C\|\gf\|_{BV}$.
More precisely,
\begin{equation}\label{bv1}
s_n(\qf)\le C_1\|\gf\|_{BV}(n+1)^{-2}, \qquad n\ge0,
\end{equation}
and
\begin{equation}\label{bv2}
s_n(\qf)\le C_2\|\gf\|'_{BV}n^{-2}, \qquad n\ge1.
\end{equation}
\end{thm}

\Pf 
Note that \eqref{bv2} follows from \eqref{bv1} since a 
symbol $\gf$ constant on \OI{} yields a rank one operator $\qf$.

We use methods from Sturm--Liouville theory, and
begin by making some simplifications.
\begin{enumerate}
\item
Replacing $\gf$ by a regularization $\gf_{\eps}$
such that $\gf_\eps\to\gf$ in $L^2\OI$ and thus $\qq{\gf_\eps}\to\qf$
in $\bS_2$ as $\eps\to0$, we see that we may assume
$\gf\in C^1\OI$.
\item
Considering real and imaginary parts separately, we may assume that
$\gf$ is real, and thus $\qf$ self-adjoint.
\item\label{bv0}
Subtracting a constant times $\chi_{\OI}$, which yields a rank 1 operator,
we may assume that $\gf(1)=0$.
\item\label{bvnorm}
By homogeneity, we may assume that $\|\gf\|'_{BV}=\int_0^1|\gf'|\le1$
and show that 
then $\|\qf\|_{\wS{1/2}}\le C$.
\item \label{bv-0}
Using 
$$
\gf=-\int_x^1\gf'(y)\,dy = -\int_x^1(\gf')_+ + \int_x^1(\gf')_-=\f_1-\f_2,
$$
and the corresponding decomposition $\qf=\qq{\f_1}-\qq{\f_2}$,
we may also assume that $\gf'\le0$ and thus $\gf\ge0$.
By \refT{T:positive}, $\qf$ is then a positive operator.
\item\label{bv-}
Similarly, writing $\gf=2\gf_1-\gf_2$ with $\gf_2=1-x$ and
$\gf_1=(\gf+\gf_2)/2$, we may further assume that $\gf'\le-1/2$ on \OI.
\end{enumerate}

Let $\gl>0$ and let, as above, $(G_\gl,g_\l)$ be the solution to
\eqref{a3} with $\ggl(0)=1$, $G_\gl(0)=0$. Thus $\gl$ is an eigenvalue
if and only if \eqref{a5} holds, \ie, by \ref{bv0}, 
if and only if $\ggl(1)=0$.

Write $\go=\gl^{-1/2}$ and express $(\ggl,\go G_\gl)$ in polar coordinates
\begin{equation}
\label{b3}
\begin{aligned}
\ggl(x)&=\Ro(x)\cos \gTo(x),\\
G_\gl(x)&=\Ro(x)\sin \gTo(x),
\end{aligned}
\end{equation}
where $\Ro(x)=\sqrt{g^2+\go^2G^2}>0$ and $\gTo$ is continuous with
$\gTo(0)=0$.
Note that $\gl$ is an eigenvalue if and only if $\cos\gTo(1)=0$, \ie{}
$\gTo(1)=n\pi+\pi/2$ for some integer $n$.

Since $\Ro(x)>0$, $\Ro$ and $\gTo$ belong to $C^1\OI$, and \eqref{b3} and
\eqref{a3} yield
\begin{equation*}
\gl^{-1/2}\Ro^2\gTo'=
gG'-Gg'
=g^2-\gl^{-1}\gf'G^2
=\Ro^2(\cos^2\gTo-\gf'\sin^2\gTo)
\end{equation*}
and thus
\begin{equation}
\label{b4}
\gTo'=\go(\cos^2\gTo-\gf'\sin^2\gTo).
\end{equation} 
In particular, since $\gf'<0$ by \ref{bv-}, $\gTo'>0$ and thus $\gTo(1)>0$.

Now suppose $0<\go<\nu$ and consider the corresponding functions
$\gTo$ and $\gTn$. We claim that
\begin{equation}
\label{b5}
\gTo(x)<\gTn(x), \qquad 0<x\le 1.
\end{equation}
Indeed, since $\gTo(0)=0=\gTn(0)$ 
and, by \eqref{b4}, $\gTo'(0)=\go<\nu=\gTn'(0)$, \eqref{b5} holds in
$(0,\gd)$ for some $\gd>0$. Hence, if \eqref{b5} fails, there exists
some $x_1\in(0,1]$ such that $\gTo(x)<\gTn(x)$ for $0<x<x_1$
but $\gTo(x_1)=\gTn(x_1)$. This would imply 
$\gTo'(x_1)\ge \gTn'(x_1)$; 
on the other hand, then
$$
\cos^2\gTo(x_1)-\gf'\sin^2\gTo(x_1)
=\cos^2\gTn(x_1)-\gf'\sin^2\gTn(x_1)>0,
$$
recalling \ref{bv-},
and \eqref{b4} would yield
$\gTo'(x_1)< \gTn'(x_1)$, a contradiction.

{}From \eqref{b5} follows in particular that 
the function $\go\mapsto \gTo(1)$ is strictly increasing. 
Hence there is for each integer $n\ge0$ at most one value of $\go$, 
$\go_n$ say, such that $\gT_{\go_n}(1)=n\pi+\pi/2$, and thus a
corresponding eigenvalue $\gl_n=\go_n^{-2}$.
(The solution to \eqref{a4} depends continuously on $\go$, with
$\gT_0(1)=0$ and $\gTo(1)\to\infty$ as $\go\to\infty$,
so $\go_n$ exists for every $n\ge1$, but we do not need that.)
Integrating \eqref{b4} we obtain by \ref{bvnorm}
\begin{equation*}
\gTo(1)=\int_0^1\gTo'(x)\,dx\le\go\int_0^1(1+|\gf'(x)|)\,dx
\le 2\go
\end{equation*}
and thus
$2\go_n\ge \gT_{\go_n}(1)=n\pi+\pi/2$, which yields
\begin{equation*}
\gl_n =\go_n^{-2}\le 4\pi^{-2}(n+{\tfrac12})^{-2}.
\quad\bl
\end{equation*}

Considering again functions on the whole half-line $\R_+$, we now can
prove an endpoint result corresponding to \refT{T:yp}.

\begin{thm}
\label{T:1/2}
If $\f\in Y_{1/2}$, then $\qf\in\wS{1/2}$.
\end{thm}
\Pf
Define $A_n^{(k)}$ and $A^{(k)}$ as in \eqref{Ank} and \eqref{Ak}, but
now for all integers $k$.
For $k\ge1$, \eqref{3.7a}  holds for every $p$, and
thus
$\f\in Y_{1/2}\subset X_{1/2}$ implies that
$$
\sum_{k\ge1}\|\cp_{A^{(k)}}\qf\|_{\bS_{1/2}}^{1/2} 
\le \sum_{k\ge1} \const 2^{-pk/2}<\infty.
$$
Moreover, by symmetry, 
$\|\cp_{A^{(-k)}}\qf\|_{\bS_{1/2}} =\|\cp_{A^{(k)}}\qf\|_{\bS_{1/2}}$ 
so
$$
\sum_{k\le-1}\|\cp_{A^{(k)}}\qf\|_{\bS_{1/2}}^{1/2} 
<\infty
$$ too.
It follows that
\begin{equation}
\label{sqf}
\qf-\cp_{A^{(0)}}\qf
=\sum_{k\neq0}\cp_{A^{(k)}}\qf\in \bS_{1/2}.
\end{equation}

Next, $\cp_{A^{(0)}}\qf$ is the direct sum of 
$\cp_{A^{(0)}_n}\qf$, $n\in\Z$, which act in the orthogonal spaces
$L^2[2^n,2^{n+1}]$. 
By translation invariance, \refL{L:homo} and \refT{T:BV}, 
with $\f_n(x)=2^n\f(2^nx+2^n)$,
\begin{equation*}
\begin{split}
\|\cp_{A^{(0)}_n}\qf\|_{\wS{1/2}}
&=\|\qf^{[2^n,2^{n+1}]}\|_{\wS{1/2}}
=\|\qq{\f_n}^{[0,1]}\|_{\wS{1/2}}
\le C\|\f_n\|_{BV\OI}\\
&\le C' 2^n\int_{2^n}^\be |d\f|.
\end{split}
\end{equation*}
By \refT{T:yp1}, we thus have
$$
\sum_{n\in\Z}\left(\|\cp_{A^{(0)}_n}\qf\|_{\wS{1/2}}\right)^{1/2}<\be
$$
and it follows from \refL{pnorm} (or as in the proof of \refL{L:limB}
below) that 
$\cp_{A^{(0)}}\qf\in\wS{1/2}$.
By \eqref{sqf}, $\qf\in\wS{1/2}$ too.
\qed

\refT{T:1/2} is the best possible; 
for any reasonably smooth $\gf$,
the singular numbers $s_n(\qf)$ decrease like $n^{-2}$ but not faster.
More precisely, we have the following very precise result.
Recall that a function in $Y_{1/2}$ has locally bounded variation and
thus is a.e.\ differentiable.
\begin{thm}
\label{T:lim}
Let $\f\in Y_{1/2}$. Then
\begin{equation}
\label{d1}
n^2s_n(\qf) \to \pi^{-2}\|\f'\|_{L^{1/2}}
=\pi^{-2}\left(\into|\f'(x)|^{1/2}dx\right)^2 
<\be
\quad\text{as }\ntoo.
\end{equation}
Equivalently,
\begin{equation}
\label{d2}
\eps^{1/2}|\set{n:s_n(\qf)>\eps}| \to \pi^{-1}\into|\f'(x)|^{1/2}dx<\be
\quad\text{as } \eps\to0.
\end{equation}
In particular, $n^2s_n(\qf)\to0$ as $\ntoo$ if and only if $\f'=0$ a.e.
\end{thm}

\Pf
Note first that by the \CSineq{}  and \eqref{yp},
\begin{equation}
\label{l1/2}
\hskip-1em
\into|\f'(x)|^{1/2}dx
\le \sum_{n\in\Z}2^{n/2} \Bigpar{\int_{2^n}^{2^{n+1}}|\f'(x)|\,dx}^{1/2}
\le \|\f\|_{Y_p}^{1/2}<\infty.
\end{equation}
For smooth and positive symbols on a finite interval, \eqref{d1}
follows by standard Sturm--Liouville theory, see \cite[\S 1.2
with the transformation in \S 1.1]{Levitan-S}.
Indeed, much more refined asymptotics of $s_n$ can be given
\cite[Chapter 5]{Levitan-S}.

We present here another proof that applies in the general case. We
prove a sequence of lemmas. The first implies that \eqref{d1} and
\eqref{d2} are equivalent.

\begin{lem}
\label{L:limA}
For any bounded operator $T$ on a Hilbert space,
$$
\limsup_{\eps\to0} \eps^{1/2}|\set{n:s_n>\eps}|
=\Bigpar{\limsup_{\ntoo} (n^2s_n)}^{1/2}
$$
and similarly with $\liminf$ instead of $\limsup$ on both sides.
\end{lem}

\Pf
If $\limsup \eps^{1/2}|\set{n:s_n>\eps}|<a$ for some $a>0$, then for
all small $\eps$, 
$|\set{n:s_n>\eps}| < a\eps^{-1/2}$. 
Taking $\eps=a^2(n+1)^{-2}$, we see 
that for large $n$, $s_n\le\eps$, and thus $(n+1)^2s_{n}\le a^2$, so 
$\limsup_{\ntoo} n^2s_n \le a^2$.
The converse is similar, and the second part follows similarly by
reversing the inequalities.
\qed

\begin{lem}
\label{L:limB}
If $T_1,\dots,T_N$ are bounded operators on Hilbert spaces
$H_1,\dots,H_N$, then 
$$
\Bigpar{\limsup_\ntoo n^2s_n(T_1\oplus\dots\oplus T_N)}^{1/2} 
\le \sum_{k=1}^N \Bigpar{\limsup_\ntoo n^2s_n(T_k)}^{1/2} 
$$
and
$$
\Bigpar{\liminf_\ntoo n^2s_n(T_1\oplus\dots\oplus T_N)}^{1/2} 
\ge \sum_{k=1}^N \Bigpar{\liminf_\ntoo n^2s_n(T_k)}^{1/2} .
$$
\end{lem}

\Pf
The singular numbers $s_n(T_1\oplus\dots\oplus T_N)$ consist of all 
$s_i(T_j)$, rearranged into a single nonincreasing sequence.
Hence,
$$
|\set{n:s_n(T_1\oplus\dots\oplus T_N)>\eps}|
=\sum_{j=1}^N|\set{n:s_n(T_j)>\eps}|
$$
and the result follows by \refL{L:limA}.
\qed

For arbitrary sums we have the following estimate.

\begin{lem}
\label{L:lim1}
If $T$ and $U$ are bounded operators in a Hilbert space, 
and $0<\gd<1$, then
\begin{align}
\label{c1}
\limsup_\ntoo n^2s_n(T+U) 
&\le (1-\gd)^{-2}\limsup_\ntoo n^2s_n(T) + \gd^{-2}\limsup_\ntoo n^2s_n(U), \\
\label{c2}
\liminf_\ntoo n^2s_n(T+U) 
&\ge (1-\gd)^{2}\liminf_\ntoo n^2s_n(T) - \gd^{-2}\limsup_\ntoo n^2s_n(U).
\end{align}
\end{lem}

\Pf
By \eqref{snsum}, $s_n(T+U)\le s_{\floor{(1-\gd)n}}(T)+s_{\floor{\gd n}}(U)$,
and \eqref{c1} follows, together with
\begin{equation}
\label{c3}
\liminf_\ntoo n^2s_n(T+U) 
\le (1-\gd)^{-2}\liminf_\ntoo n^2s_n(T) + \gd^{-2}\limsup_\ntoo n^2s_n(U) 
\end{equation}
Replacing here $T$ by $T+U$ and $U$ by $-U$, we obtain \eqref{c2} by
rearrangement. 
\qed

Letting $\gd\to0$ in \eqref{c2} and \eqref{c3}, we obtain the 
following result by Fan \cite{GK1}.
\begin{lem}
\label{L:lim2}
If $T$ and $U$ are bounded operators in a Hilbert space,
$\lim_\ntoo n^2s_n(T)$ exists and $n^2s_n(U)\to0$ as $\ntoo$, then
$\lim_\ntoo n^2s_n(T+U)=\lim_\ntoo n^2s_n(T)$.
\nopf
\end{lem}

\begin{lem}
\label{L:closed}
The set of $\f\in Y_{1/2}$ such that \eqref{d1} holds is a closed set.
\end{lem}

\Pf
Suppose that $\f_k\to\f$ in $Y_{1/2}$ and that \eqref{d1} holds for
each $\f_k$. 
By \refL{L:lim1} and \refT{T:1/2}, for every $k$ and $0<\gd<1$,
\begin{equation}
\label{lc1}
\begin{split}
\limsup_\ntoo n^2s_n(\qf) 
&\le (1-\gd)^{-2}\limsup_\ntoo n^2s_n(\qq{\f_k}) 
  + \gd^{-2}\limsup_\ntoo n^2s_n(\qq{\f-\f_k}) \\
&\le (1-\gd)^{-2}\pi^{-2}\|\f_k'\|_{L^{1/2}} + C\gd^{-2}\|\f-\f_k\|_{Y_{1/2}}
\end{split}
\end{equation}
and similarly
\begin{equation}
\label{lc2}
\liminf_\ntoo n^2s_n(\qf) 
\ge (1-\gd)^{2}\pi^{-2}\|\f_k'\|_{L^{1/2}} - C\gd^{-2}\|\f-\f_k\|_{Y_{1/2}}.
\end{equation}
Moreover, by \eqref{l1/2},
$\|(\f-\f_k)'\|_{L_{1/2}} \le \|\f-\f_k\|_{Y_{1/2}} \to0$ 
as $k\to\be$, and so 
$\|\f_k'\|_{L^{1/2}}\to \|\f'\|_{L^{1/2}}$.
Letting first $k\to\be$ and then $\gd\to0$ in \eqref{lc1} and
\eqref{lc2}, we obtain \eqref{d1}. 
\qed

\begin{lem}
\label{L:linear}
If $\f$ is linear on a finite interval $I$, then
\begin{equation*}
n^2s_n(\qf^I) \to \pi^{-2}\left(\int_I|\f'(x)|^{1/2}dx\right)^2 
\quad\text{as }\ntoo.
\end{equation*}
\end{lem}

\Pf
Let $\f(x)=\ga+\gb x$ with $\ga,\gb$ complex numbers.
Suppose first that $I=\OI$. By the example at the beginning of the
section and homogeneity, 
$$
s_n(Q^I_{-\gb+\gb x}) = |\gb|\pi^{-2}(n+{\tfrac12})^{-2}
$$
so $n^2 s_n(Q^I_{-\gb+\gb x}) \to \pi^{-2}|\gb|$ as $\ntoo$.
Since $Q^I_{\ga+\gb x} - Q^I_{-\gb+\gb x} = Q^I_{\ga+\gb}$
is a rank one operator, \refL{L:lim2}
(or, more simply, $s_{n+1}(Q^I_{-\gb+\gb x}) \le s_{n}(Q^I_{\ga+\gb x})
\le s_{n-1}(Q^I_{-\gb+\gb x})$)
yields
\begin{equation}
\label{x7}
n^2s_n(\qf^I) \to \pi^{-2}|\gb|=\pi^{-2}\left(\int_I|\f'(x)|^{1/2}dx\right)^2 .
\end{equation}

If $I=[0,a]$, we have by \refL{L:homo} and \eqref{x7}
\begin{equation*}
n^2s_n(\qf^I) 
=n^2s_n(\qq{\f_a}^{\OI}) 
\to \pi^{-2}a^2|\gb|=\pi^{-2}\left(\int_I|\f'(x)|^{1/2}dx\right)^2,
\end{equation*}
and the general case follows by
translation invariance.
\qed

\smallbreak
{\bf Completion of the proof of \refT{T:lim}.}

\pfitem1{$\f$ is piecewise linear on $\OI$ and $\f=0$ on $(1,\be)$}
Let $0=t_0<t_1<\dots<t_N=1$ be such that $\f$ is linear on every 
$I_i=[t_{i-1},t_i]$, $i=1,\dots,N$.
Let $H_i=L^2(I_i)$, so $L^2\OI=H_1\oplus\dots\oplus H_N$, and let 
$P_i:L^2\OI\to H_i$ denote the orthogonal projection.

Since each $P_i\qf P_j$, $i\neq j$, has rank 1,
$\qf-\sum_{i=1}^N P_i\qf P_i$ has finite rank 
and by  \refL{L:lim2} (or directly), it suffices to consider 
$\sum_{i=1}^N P_i\qf P_i=\qf^{I_1}\oplus\dots\oplus\qf^{I_N}$.
By \refL{L:linear}, 
\begin{equation*}
n^2s_n(\qf^{I_i}) 
\to \pi^{-2}\left(\int_{t_{i-1}}^{t_i}|\f'(x)|^{1/2}dx\right)^2 ,
\end{equation*}
and thus \refL{L:limB} yields
$$
\Bigpar{\lim_\ntoo n^2s_n\Bigpar{\sum_{i=1}^N P_i\qf P_i}}^{1/2} 
= \sum_{i=1}^N \Bigpar{\lim_\ntoo n^2s_n(\qf^{I_i})}^{1/2} 
= \pi^{-1}\int_{0}^{1}|\f'(x)|^{1/2}dx,
$$
which proves \eqref{d1}.

\pfitem2{$\f$ is absolutely continuous on $\OI$ and $\f=0$ on $(1,\be)$}
Approximate $\f'$ by step functions $h_n$ such that
$\|\f'-h_n\|_{L^1\OI} <1/n$, 
and let
\begin{equation*}
\psi_n(x)=
\begin{cases}
\f(0)+\int_0^x h_n(y)\,dy, & x\le1,\\
0, & x>1.
\end{cases}
\end{equation*}
Then \eqref{d1} holds for each $\psi_n$ by Step 1, and $\psi_n\to\f$ in
$BV\OI$ and thus in $Y_{1/2}$, 
see \refC{C:yp2}, so \eqref{d1} holds by \refL{L:closed}.

\pfitem3{$\f$ has bounded variation on $\OI$, $\f=0$ on $(1,\be)$ and
$\f$ is singular, \ie, $\f'=0$ almost everywhere} 
We may assume that $\f$ is right-continuous. Then
$\f(x)=\f(0)+\int_0^x d\mu$ for some singular complex measure
$\mu$ supported on $\OI$.
Given any $\eps>0$, there thus exists a sequence of intervals
$(I_i)_1^\be$ in $\OI$ such that $\sum_1^\be |I_i|<\eps$ and
$|\mu|\bigpar{\OI\setminus\bigcup_1^\be I_i}=0$.
Let $N$ be a positive integer such that
$|\mu|\bigpar{\OI\setminus\bigcup_1^N I_i}<\eps$.

We may assume that each $I_i$ is closed, and by combining any two of
$I_1,\dots,I_N$ that overlap, we may assume that $I_1,\dots,I_N$ are
disjoint.
The complement 
\linebreak$\OI\setminus\bigcup_1^N I_i$ 
is also a finite disjoint
union of intervals, say $\bigcup_1^M J_j$.

For each interval $I$, \refT{T:BV} and \refL{L:homo} yield
\begin{equation}
\label{axel}
\sup_{n\ge1} n^2 s_n(\qf^I)
\le C |I| \|\f\|'_{BV(I)}
\le C |I|\, |\mu|(I).
\end{equation}
Moreover, as in Step 1 of the proof, 
$$
\qf = \qf^{I_1} \oplus\dots \oplus \qf^{I_N}
\oplus  \qf^{J_1} \oplus\dots \oplus \qf^{J_M}
+R,
$$
where $R$ has finite rank. Hence, by \refL{L:limB}, \eqref{axel} and
the \CSineq,
\begin{equation*}
\begin{split}
\Bigpar{\limsup_\ntoo n^2s_n(\qf)}^{1/2} 
&\le \sum_{i=1}^N \Bigpar{\limsup_\ntoo n^2s_n(\qf^{I_i})}^{1/2} 
+\sum_{j=1}^M \Bigpar{\limsup_\ntoo n^2s_n(\qf^{J_j})}^{1/2} \\
&\le C\sum_{i=1}^N \Bigpar{|I_i|\,|\mu|(I_i)}^{1/2} 
+ C\sum_{j=1}^M \Bigpar{|J_j|\,|\mu|(J_j)}^{1/2} \\
&\le C\Bigpar{\sum_{i=1}^N|I_i|}^{1/2} \Bigpar{\sum_{i=1}^N|\mu|(I_i)}^{1/2}\\ 
&\qquad+ 
 C\Bigpar{\sum_{j=1}^M|J_j|}^{1/2}\Bigpar{\sum_{j=1}^M|\mu|(J_j)}^{1/2} \\
&\le C \eps^{1/2}\bigpar{|\mu|\OI}^{1/2} + C\cdot1\cdot\eps^{1/2}. 
\end{split}
\end{equation*}
The result $n^2 s_n(\qf)\to0$
follows by letting $\eps\to0$.

\pfitem4{$\f$ has bounded variation on $(0,a)$ and $\f=0$ on $(a,\be)$
for some $a>0$} 
By \refL{L:homo}, it suffices to consider the case $a=1$.
We can decompose $\f=\f_{\rm a}+\f_{\rm s}$ on $\OI$, with $\f_{\rm a}$ absolutely
continuous and $\f_{\rm s}$ singular; let $\f_{\rm a}=\f_{\rm s}=0$ on $(1,\be)$.
By Steps 2 and 3, 
$$
n^2 s_n(\qq{\f_{\rm a}})\to
\pi^{-2} \|\f_{\rm a}'\|_{1/2}
=\pi^{-2} \|\f'\|_{1/2}
$$
and $n^2 s_n(\qq{\f_{\rm s}})\to0$,
and the result follows by \refL{L:lim2}.

\pfitem5{$\f\in Y_{1/2}$ is arbitrary}
Define, for $N\ge1$,
\begin{equation*}
\f_N(x)=
\begin{cases}
\f(1/N)-\f(N), & 0<x\le 2^{-N},\\
\f(x)-\f(N), & 2^{-N}<x \le 2^N, \\
0, & 2^{N}<x.
\end{cases}
\end{equation*}
It is easily seen that each $\f_N$ has bounded variation and that
$\|\f-\f_N\|_{Y_{1/2}}\to0$ as $N\to\be$, \cf{} \eqref{yp}.
Thus the result follows by Step 4 and \refL{L:closed}.
\qed

As corollaries, we obtain new proofs of some results from
\refS{S:1/2}.

\begin{cor}
\label{C:bvlim}
If $I$ is a finite interval and $\f$ has bounded variation on $I$, then
\begin{equation*}
n^2s_n(\qf^I) \to \pi^{-2}\Bigpar{\int_I|\f'(x)|^{1/2}dx}^2
\quad\text{ as }\ntoo.
\end{equation*}
\end{cor}

\Pf
By translation invariance, we may assume $I=[0,a]$. Then, defining
$\f=0$ outside $I$, we have $\f\in Y_{1/2}$ by \refC{C:yp2}, and the
result follows by \refT{T:lim}.
$\bl$

\begin{cor}\label{C:moc}
If $\f$ has locally bounded variation and $\qf\in \bS_{1/2}$, then
$\f'=0$ a.e.
\end{cor}

\Pf
If $0<a<b<\be$, then $\f$ has bounded variation on $[a,b]$, and since
$$
n^2s_n(\qf^{[a,b]}) \le n^2s_n(\qf)\to0,
$$
\refC{C:bvlim} yields $\int_a^b|\f'|^{1/2}=0$.
Hence, $\f'=0$ a.e.
$\bl$

\begin{cor}
\label{C:not1/2}
If $\f$ is locally absolutely continuous and $\qf\in \bS_{1/2}$, then
$\f=0$ a.e.
\nopf
\end{cor}

\begin{remark} 
More generally, in the last two corollaries, 
$\bS_{1/2}$ can be replaced by any
Schatten--Lorentz space $\bS_{1/2,q}$ with $q<\be$.
\end{remark}

\

\setcounter{equation}{0}
\section{\bf More on $\bs{p=1}$}
\label{more}

\

Although $\f\in X_1$ does not imply $\qf,\qf^+\in \bS_1$,
the corresponding weak results holds. 
 There is, however, a striking difference between $\qf$ and $\qf^+$;
as is shown in (ii) and (iii) below, for every $\f\in X_1$ not \aex{}
equal to 0, $ns_n\to0$ for $\qf$ but not for $\qf^+$.
(Note that \refT{E:exp} implies that nothing can be said about the
rate of convergence of $ns_n(\qf)$ to 0.
In particular, if $q<\be$, then 
$\f\in X_1$ does not imply that $Q_\f\in\bS_{1,q}$.) 

\begin{thm}\label{T:ws1}
If $\f\in X_1$ then the following hold:
\begin{enumerate}
\item
$\qf,\,\qf^+,\,\qf^-\in\wS1$.
\item
$ns_n(\qf)\to0$ as $\ntoo$.
\item
$ns_n(\qf^+)=ns_n(\qf^-)\to\pi^{-1}\into|\f(x)|\,dx$ as $\ntoo$.
\end{enumerate}
\end{thm}

\Pf
Since $X_1\subset X_\be$, $\qf^+$ is bounded by \refT{bou}.
By \refT{fact},
\begin{equation*}
\bigpar{\qf^+}^*\qf^+ = \qq\Phi
\end{equation*}
where $\Phi(x)=\int_x^\be |\f(y)|^2\,dy$.
By \eqref{xpiii}, $x^{1/2}\Phi(x)^{1/2}\in L^1(dx/x)$,
so $\Phi\in Y_{1/2}$ by \refT{T:yp1} and hence
$\qq\Phi\in\wS{1/2}$ by \refT{T:1/2}.
Consequently, $\qf^+\in\wS1$. The same holds for
$\qf^-=\bigpar{\qf^+}^*$ and $\qf=\qf^++\qf^-$.

Moreover, $s_n(\qq\Phi)=s_n(\qf^+)^2$, 
and thus (iii) follows from \refT{T:lim}
applied to $\Phi$.

For (ii), we observe that $\f\mapsto\qf$ thus is a bounded linear map
$X_1\to\wS1$, and that the set of $C^1$ functions with compact support
is dense in $X_1$ and mapped (by \refT{T:1/2}) into the closed
subspace
$\wS1^0=\set{T\in\wS1: ns_n(T)\to0\text{ as }\ntoo}$ 
of $\wS1$.
Hence $\qf\in\wS1^0$ for every $\f\in X_1$.
\qed

\begin{remark} 
Under considerably more restrictive conditions on $\f$,
\refT{T:ws1}\vvp(iii) 
was obtained in \cite{EEH}. See also the related result
in \cite[Remark IV.8.3]{GK1} and 
\cite[Theorem III.2.4]{GK2}.
\end{remark}

\begin{remark}
For the related operators $\qmu$, we similarly obtain that if
$x\mu(x,\be)\in L^{1/2}(dx/x)$, then $\qmu\in \wS1$, and
$$
ns_n(\qmu)\to\pi^{-1}\into\left(\frac{d\mu}{dx}\right)^{1/2}\,dx
\quad\mbox{as }\ntoo,
$$
where $\frac{d\mu}{dx}$ is the Radon--Nikodym derivative of $\mu$. 
In particular, for such $\mu$, \linebreak$ns_n(\qmu)\to0$ if and only if $\mu$
is singular.
\end{remark}

We saw earlier that $\varphi\in X_{1}$ is not enough to insure that
$Q_{\varphi}$ is in the trace class. Furthermore the previous theorem shows
that if $\varphi\in X_{1}$ then neither $Q_{\varphi}^{+}$ nor 
$Q_{\varphi}^{-}$ will be in the trace class.  However, the
combination of size and 
regularity results for singular numbers given in the previous theorem does
insure that these operators have a well defined Dixmier trace.  Because of
the recent interest in the Dixmier trace we digress briefly to record this
observation.  For more about the Dixmier trace and its uses we refer to
IV.2.$\beta$ of \cite{C}.

Let $\ell^{\infty}$ be the space of bounded sequences indexed by non-negative
integers and let $c_1$ be the closed subspace consisting
of sequences $\left\{  a_{n}\right\}  $ for which $\lim_{n}a_{n}$ exists.  It
follows from the Hahn--Banach theorem that the functional 
$\lim(\cdot)$ which is defined on $c_1$ has a positive
continuous extension, $\lim_{\omega}(\cdot)$, to all of
$\ell^{\infty}$. By saying $\lim_{\omega}( \cdot)  $ is positive we
mean that if $a_{n}\geq0$ for $n=0,1,2,\dots$ then $\lim_{\omega}\left(
\left\{  a_{n}\right\}  \right)  \geq0$. This extension is not unique and we
are using the subscript $\omega$ to denote the particular choice.  It was
noted by Dixmier in \cite{D} that $\lim_{\omega}(  \cdot)  $ can
also be selected to have the following scaling property:
\[
\lim\nolimits_{\omega}\left(  a_{0},a_{0},a_{1},a_{1},a_{2},a_{2},\dots\right)
=\lim\nolimits_{\omega}\left(  a_{0},a_{1},a_{2},\dots\right)  .
\]
A simple proof is in \cite{C}. (Although the scaling is important for the
general theory it has no role in our discussion.)

Consider now the operator ideal 
$\bS_\O\supset\bS_{1,\be}$  
that consists of the
operators $T$ on Hilbert space 
such that 
\begin{equation}
\label{SO}
\|T\|_{\bS_\O}
\df\sup_{n\ge0}\frac{\sum\limits_{k=0}^ns_k(T)}{\sum\limits_{k=0}^n\frac1{k+1}}
<\be. 
\end{equation}
Suppose that $T\in\bS_{\O}$. 
For a fixed choice of
$\lim_{\omega}(  \cdot)  $ we define a Dixmier trace,
${\trace}_{\omega}(  \cdot)$, as follows.  For
positive $T\in\bS_{\O}$ set
\[
{\trace}_{\omega}(  T)  =\lim
\nolimits_{\omega}\left(  \left\{ \frac{1}{\log(n+2)}
\sum_{k=0}^{n}s_{k}(T)  \right\}  \right)  .
\]
Although perhaps not obvious at first glance,  
it is not difficult to see
that, in fact, if $T_{1}$ and $T_{2}$ are two \textit{positive} operators in
$\bS_{\O}$ then ${\trace}_{\omega}(
T_{1}+T_{2})  ={\trace}_{\omega}(
T_{1})  +{\trace}_{\omega}(  T_{2})  $.
A proof of this is also in \cite{C}.  Using this fact, the functional
${\trace}_{\omega}(  \cdot)  $ can be extended
uniquely by linearity to all of $T\in\bS_{\O}$. For
$T\in\bS_{\O}$ the value of 
${\trace}_{\omega}(  T)  $ need not be independent
of $\omega$. 
However, there are certain operators for which 
${\trace}_{\omega}(T)$ is 
independent of $\omega$.  Such operators are defined to be
\textit{measurable}.
In this case we will write
${\trace}_{\rm D}(  T)  $ for this common value
and refer to it as \textit{the} Dixmier trace of $T$.

\newlength\sjllgd\newlength\sjrlgd
\newlength\sjllgdb\newlength\sjrlgdb
\newcommand\SJaligna[2]{\settowidth{\sjllgd}{$\displaystyle #1$}
 \settowidth\sjrlgd{$\displaystyle #2$}
 \[#1#2\]}
\newcommand\SJalignb[2]{\settowidth{\sjllgdb}{$\displaystyle #1$}
 \settowidth\sjrlgdb{$\displaystyle #2$}
 \[\hskip\sjllgd\hskip-\sjllgdb
   #1#2
   \hskip\sjrlgd\hskip-\sjrlgdb
 \]}
\begin{cor}\label{C:Dix}\quad
\begin{enumerate}
\item  
If $\varphi\in X_1$, then the operators $|Q_\f^+|$ and
$|Q_\f^-|$ are measurable and
\SJaligna
{\trace_{\rm D}(|Q_\f^{+}|)} 
{=\trace_{\rm D}(|Q_\f^{-}|)=\frac1{\pi}\int_0^\be|\f( x)|dx.}

\item  
If $\varphi\in Y_{1/2}$, then $|Q_\f|^{1/2}$ is measurable and
\SJalignb
{\trace_{\rm D}(|Q_\f|^{1/2})}
{=\frac1{\pi}\int_0^\be|\f^{\prime}(x)|^{1/2}dx.}

\item  
If $\f\in X_1$, then $Q_\f$ is measurable and
\SJalignb
{\trace_{\rm D}(Q_\f)}
{=0.}

\item  
If $\varphi\in X_{1}$, then $Q_\f^+$ and $Q_\f^-$ are measurable and
\SJalignb
{\trace_{\rm D}(Q_\f^+)}
{=\trace_{\rm D}(Q_\f^-)=0.}
\end{enumerate}
\end{cor}

\Pf 
We start with (i). From the very definitions
$s_{n}(|Q_\f^{+}|)=s_{n}(Q_\f^+)$ 
and hence, the
previous theorem gives the asymptotic behavior of
$\{s_{n}(|Q_\f^+|)\}_{n\ge0}$. 
Those asymptotics, together
with the fact that $|Q_\f^+|$ is a positive operator, insure that
$|Q_\f^+|$ is 
measurable and has the indicated Dixmier trace. 
A similar argument applies to
$|Q_\f^{-}|$ 
and, after
noting that $s_{n}\left(|Q_\f|^{1/2}\right)=s_{n}(Q_\f)^{1/2}$ and taking note
of Theorem \ref{T:lim}, to part (ii).

We now consider (iii). By Theorem \ref{T:ws1}, we have
$\lim\limits_{n\to\be}ns_n(Q_\f)=0$.
Also $s_{n}(Q_\f)=s_{n}(Q_\f^{\ast})$.
Thus by \eqref{snsum}, 
$\lim\limits_{n\to\be}ns_{n}((Q_\f+Q_\f^\ast)/2)=0$. We now
use
the spectral projection to write $\frac12(Q_\f+Q_\f^\ast)$ as a
difference of two positive operators $\frac12(Q_\f+Q_\f^\ast)_{\pm}$ and note
that we will have 
$$
\lim_{n\to\be}ns_{n}\left(\tfrac12(Q_\f+Q_\f^\ast)_{\pm}\right)=0.
$$ 
Arguing similarly with the skew-adjoint part of
$Q_\f$, we realize $Q_\f$ as a linear combination of four
positive operators each of which have singular numbers which tend to zero
more rapidly than $n^{-1}$.
Those positive operators are certainly measurable and have Dixmier trace zero.
The result we want now follows by the linearity of $\trace_{\rm D}(\cdot)$.

For (iv) we first pick and fix a choice $\trace_\o(\cdot)$. Assume for the
moment that $\f$ is
real, supported in $[0,1]$ and in $L^{2}$. 
$R_\f\df Q_\f^+-Q_\f^{+\ast}=Q_\f^+-Q_\f^-$ has
real anti-symmetric kernel.  Thus it is normal and its eigenvalues are
imaginary and symmetric. Hence,
${\rm i}R_\f$ is symmetric and its positive and negative parts, 
$({\rm i}R_\f)_\pm$ are unitarily
equivalent. Thus
$$
\trace_\o(R_\f)
=-{\rm i}\,\trace_\o((\rmi R_\f)_+)+{\rmi}\,\trace_\o((\rmi R_\f)_-)=0.
$$ 
Taking note of the fact
that $\lim_\o(\cdot)$ is continuous on $\ell^\be$ and of the norm estimates
implicit in the previous theorem we see that we can extend this result by
linearity and continuity and conclude that $\trace_\o (R_\f)=0$ for all $\f\in
X_{1}$.
Now we use the fact that $\o$ was arbitrary to conclude $\trace_{\rm
D}(R_\f)=0$. By
linearity this result together with the result in (iii) yields (iv).
$\bl$

For a function $\f$ defined on a finite or infinite interval $I$, 
we define the standard and $L^p$ moduli of 
continuity by
\begin{equation}
\label{contmod}
\begin{aligned}
\hskip-4em
\oxzf(h;I)&\df\sup\set{|\f(x)-\f(y)|: x,y\in I,\,|x-y|\le h},\\
\hskip-4em
\oxpf(h;I)&\df\sup_{0\le s\le h}
\biggpar{\int_{I\cap(I-s)} |\f(x+s)-\f(x)|^p\,dx}^{1/p},	
\quad 1\le p<\be, \hskip-2em
\end{aligned}
\end{equation}
where 
$0<h\le|I|$ and
$I-s=\set{x-s:x\in I} =\set{x:x+s\in I}$.
It follows easily from Minkowski's inequality that
\begin{equation}\label{omega2}
\oxpf(h;I)\le 2\oxpf(h/2;I), \qquad 1\le p\le\be.
\end{equation}
Note further that for a finite interval $I$,
\begin{equation}\label{omegap}
\oxpf(h;I)\le |I|^{1/p-1/q}\oxx {\f}q (h;I), \qquad p<q\le\be.
\end{equation}

We often omit $I$ from the notation.

An alternative $L^p$ modulus of continuity is defined by
\begin{equation*}
\oypf(h;I)
\df
\left((2h)^{-1}\iint
\limits_{\textstyle{{x,y\in I }\atop{|x-y|<h}}}
|\f(x)-\f(y)|^p\,dx\,dy\right)^{1/p}.
\end{equation*}
This is equivalent to $\oxpf(h;I)$ defined above by the following lemma,
which probably is well-known to some experts. 

\begin{lem}\label{L:moduli}
Let $1\le p<\be$. Then, with $C_p$ depending on $p$ only,
$$
\oypf(h;I)
\le \oxpf(h;I)
\le C_p\oypf(h;I).
$$
\end{lem}

\Pf
The left hand inequality follows by
\begin{equation*}
\begin{split}
\bigpar{\oypf(h;I)}^p
&=\frac1h\iint
\limits_{\textstyle{{x,y\in I }\atop{0<y-x<h}}}
|\f(y)-\f(x)|^p\,dx\,dy\\
&=\frac1h\int_0^h \int_{x\in I\cap(I-s)}
|\f(x+s)-\f(x)|^p\,dx\,dy
\le \bigpar{\oxpf(h;I)}^p.
\end{split}
\end{equation*}

For the converse, we assume for convenience that $I=\OI$.
The result then follows for every finite $I$ by a linear change of
variables, and for infinite $I$ by considering $I\cap[-n,n]$ and
letting $\ntoo$. Thus $I=\OI$ and $I\cap(I-s)=[0,1-s]$.

Let $\f_s(x)=\f(x+s)$.
Assume first that $h\le1/2$. Then, for $0\le s,t\le h$, by
Minkowski's inequality,
$$
\|\f-\f_s\|_{L^p[0,1/2]}^p
\le C_p\|\f-\f_t\|_{L^p[0,1/2]}^p +C_p\|\f_s-\f_t\|_{L^p[0,1/2]}^p.
$$
Averaging over $t\in[0,h]$ we find
\begin{equation*}
\begin{split}
\|\f-\f_s\|_{L^p[0,1/2]}^p
&\le\frac{C_p}h \int_0^h\int_{0}^{1/2}
\Bigpar{|\f(x)-\f(x+t)|^p+|\f(x+s)-\f(x+t)|^p}dx\,dt\\[2mm]
&\le\frac{C_p}h 
\iint\limits_{\textstyle{{x,y\in \OI }\atop{0<y-x<h}}}
|\f(y)-\f(x)|^p\,dx\,dy
=C_p\bigpar{\oypf (h;\OI)}^p.
\end{split}
\end{equation*}
A similar argument, now taking $s-h\le t\le s$, yields the same
estimate for $\|\f-\f_s\|_{L^p[1/2-s,s]}^p$, and summing we find
\begin{equation}\label{sofie}
\|\f-\f_s\|_{L^p[0,1-s]} \le C_p \oypf(h;\OI)
\end{equation}
for every $0\le s\le h$, which proves the result for $h\le1/2$.

If $1/2<h\le1$,
the result follows from the case $h\le1/2$ and \eqref{omega2}.
\qed

For simplicity, we state the following lemma for $I=\OI$ only.

\begin{lem}
\label{L:K}
Let $1\le p<\be$. If $\f\in L^p\OI$ and $0<t\le1$, 
there exists a decomposition $\f=\f_0+\f_1$ with
$$
\|\f_0\|_{L^p\OI}\le C_p \oxpf(t)
\qquad\text{and}\qquad
\|\f_1\|_{BV\OI}'\le C_pt^{-1} \oxpf(t).
$$
\end{lem}
In other words, the Peetre
$K$-functional (see \cite{BL}), can be estimated by
$$
K(t,\f; L^p\OI, BV'\OI) \le C_p \oxpf(t),
\qquad 0<t\le1.
$$

\begin{proofx}
Take 
$
\f_1(x)=\frac1t\int_{(1-t)x}^{(1-t)x+t} \f(y)\,dy
$
and $\f_0=\f-\f_1$.
Then $\f_1$ is absolutely continuous, and thus
\begin{equation*}
\begin{split}
\|\f_1\|'_{BV}
&=\int_0^1|\f_1'(x)|\,dx
=\frac{1-t}t\int_0^1\bigl|\f\bigpar{(1-t)x+t}-\f\bigpar{(1-t)x}\bigr|\,dx\\
&=\frac{1}t\int_0^{1-t}|\f(y+t)-\f(y)|\,dy\\
&\le\frac{1}t\biggpar{\int_0^{1-t}|\f(y+t)-\f(y)|^p\,dy}^{1/p}
\le \frac1t \oxpf(t).
\end{split}
\end{equation*}
Moreover, using \Holder's inequality again,
\begin{equation*}
\begin{split}
|\f_0(x)|^p 
=\biggl| \frac1t\int_{(1-t)x}^{(1-t)x+t}\bigpar{\f(x)- \f(y)}\,dy
\biggr|^p
\le \frac1t\int_{(1-t)x}^{(1-t)x+t} |\f(x)- \f(y)|^p\,dy
\end{split}
\end{equation*}
and thus
\begin{equation*}
\int_0^1|\f_0(x)|^p 
\le \frac1t
\iint\limits_{\textstyle{{x,y\in \OI}\atop{|y-x|<t}}}
|\f(y)-\f(x)|^p\,dx\,dy
=2\bigpar{\oxpf(t)}^p.
\quad\bl
\end{equation*}
\end{proofx}

\begin{thm}
\label{T:omegafi}
If $I$ is a finite interval and $\f\in L^2(I)$, then
\begin{equation*}
s_n(\qf^I)\le C \frac{|I|^{1/2}} n \oxiif\!\left(\frac{|I|}n\right) 
\le C \frac{|I|}n\oxzf\!\left(\frac{|I|}n\right), 
\qquad n\ge1.
\end{equation*}
\end{thm}

\Pf 
By a linear change of variables, we may assume that $I=\OI$, \cf{}
\refL{L:homo}. 
Then, using the decomposition given by \refL{L:K} with $t=1/n$,
\eqref{snsum}, \refT{T:ws1} and \refT{T:BV}, we find,
for $n\ge1$,
$$
s_{2n-1}(\qf)
\le s_{n-1}(\qq{\f_0}) + s_n(\qq{\f_1})
\le C n^{-1}\norm{\f_0}_{L^2} + C n^{-2} \norm{\f_1}_{BV}'
\le C n^{-1}\oxiif(1/n),
$$
and the result follows, using \eqref{omega2} and \eqref{omegap}.
$\bl$

In particular, we see that a Dini condition implies
$\qf\in\bS_1$.

\begin{cor}
\label{C:Dini}
If $\f\in L^2\OI$ is such that
$\int_0^1\oxiif(t)\frac{dt}t<\infty$, 
in particular if $\int_0^1\oxzf(t)\frac{dt}t<\infty$, 
then $\qf^{\OI}\in \bS_1$.
\end{cor}

\Pf 
\refT{bou} shows that $\qf$ is bounded, and
Theorem \ref{T:omegafi} yields
$$
\sum_{n=2}^\be s_n(\qf) 
\le C \sum_{n=2}^\be\frac1n\oxiif\parfrac 1n
\le C \int_0^1 \oxzf(t)\frac{dt}t.
\quad\bl
$$

By a simple change of variables, 
\refC{C:Dini} applies to other finite intervals too.
Moreover, for functions $\f$ on $\R_+$,
we have the following corresponding
sufficient conditions for $\qf\in \bS_1$.

\begin{thm}\label{T:Dini2}
If $\f\in X_1$ and
$$
\sum_{n=-\be}^\be 2^{n/2}\int_0^{2^n} \oxiif\bigpar{t;[2^n,2^{n+1}]} 
\frac{dt}t <\be
$$
then
$\qf\in\bS_1$.
\end{thm}

\begin{proof}
Let $I_n=[2^n,2^{n+1}]$.
Then \refT{T:omegafi} yields
\begin{equation*}
\begin{split}
\norm{\qf^{I_n}}_{\bS_1}
&= \sum_{k=0}^\be s_k(\qf^{I_n})
\le 2 \norm{\qf^{I_n}}_{\bS_2}
+ C|I_n|^{1/2}\sum_{k=2}^\be \frac1k \oxiif\!\left(\frac{2^n}k;I_n\right)\\
&\le C 2^{n/2}\norm{\f}_{L^2(I_n)} 
+ C 2^{n/2}\int_0^{2^n} \oxiif(t;I_n) \frac{dt}t
\end{split}
\end{equation*}
and the result follows by \refT{S1} and \eqref{xpi}.
\end{proof}

\begin{cor}\label{C:Dini2}
If $\f\in X_1$ and
$$
\sum_{n=-\be}^\be 2^{n}\int_0^{2^n} \oxzf \bigpar{t;[2^n,2^{n+1}]} \frac{dt}t
<\be
$$
then
$\qf\in\bS_1$.
\nopf
\end{cor}

Note that for the functions $\f_N$ considered in \refT{E:exp},
the estimate of the singular numbers in \refT{T:omegafi}
is sharp (within a constant factor)
and the estimates of the $\bS_1$ norm implicit in \refC{C:Dini},
\refT{T:Dini2} and \refC{C:Dini2}
are of the right order.

We do not know whether the condition in \refT{T:Dini2} is necessary,
but we will give a related 
necessary condition using the $L^1$ modulus of continuity in
\refS{neces}.

We have in these applications of \refT{T:omegafi} considered $\bS_1$ only,
but the same arguments apply to $\bS_p$ for other $p$ too. In
particular, \refT{T:Dini2} extends as follows (see the remark after
\refT{S1}).

\begin{thm}\label{T:Dinip}
Let $1/2<p\le1$.
If $\f\in X_p$ and
$$
\sum_{n=-\be}^\be 2^{n(1-p/2)}\int_0^{2^n} 
 \bigpar{\oxiif\bigpar{t;[2^n,2^{n+1}]}}^p t^{p-2}\,dt <\be,
$$
then
$\qf\in\bS_p$.
\nopf
\end{thm}

Note also the following immediate consequence of \refT{T:omegafi}.
\begin{cor}
If $I$ is a finite interval and $\f$ satisfies a
\Holder\
(Lipschitz) 
condition 
$|f(x)-f(y)|\le C |x-y|^\ga$ 
for $x,y\in I$, where $0<\ga\le1$, then
\linebreak$\qf\in \bS_{1/(1+\ga),\be}$ and thus $\qf\in\bS_p$ for every
$p>1/(1+\ga)$.
\nopf
\end{cor}

\

\setcounter{equation}{0}
\section{\bf Averaging projection}
\label{aver}

\

In this section we study properties of the averaging projection onto
the set of operators of 
the form $Q_\psi$. Let us first define the averaging projection on
$\bS_2$. Let $T$ be  
an operator on $L^2(\R_+)$ of class $\bS_2$ with kernel function
$k=k_T\in L^2((\R_+)^2)$. 
We define the function $\f$ on $\R_+$ by
\bay
\label{fik}
\f(x)=\frac{1}{2x}\left(\int_0^xk(x,t)dt+\int_0^xk(s,x)ds\right),\quad x>0.
\ey
We define the averaging projection $\cal Q$ on $\bS_2$ by
$$
{\cal Q}T\df Q_\f.
$$
It is not hard 
to see that if $Q_\psi\in\bS_2$, then  ${\cal
Q}Q_\psi=Q_\psi$. It is also easy to 
see that $\|{\cal Q}T\|_{\bS_2}\le\|T\|_{\bS_2}$ for any $T\in\bS_2$,
and so $\cal Q$ is the 
orthogonal projection of $\bS_2$ onto the set of operators of the form
$Q_\psi$. 

We show in this section that $\cal Q$ is a bounded linear operator on
$\bS_p$ for \linebreak$1<p\le2$. This 
allows us to define by duality the projection $\cal Q$ on the classes
$\bS_p$ for $2\le p<\be$. 
We also show that $\cal Q$ is unbounded on $\bS_1$ but it has weak
type (1,1), i.e., 
$s_n({\cal Q}T)(1+n)\le\const\|T\|_{\bS_1}$. Finally, we use this
result to show that $\cal Q$ maps 
the Matsaev ideal into the set of compact operators.

\begin{thm}
\label{1p2}
Let $1<p\le2$. Then $\cal Q$ is a bounded projection on $\bS_p$.
\end{thm}

\Pf Let $T$ be an integral operator in $\bS_p$ with kernel function
$k$ and let $\f$ be defined by 
\rf{fik}. We have to show that $\f\in X_p$ (see the definition in the
Introduction). We can identify in 
a natural way the dual space $X_p^*$ with the space $Z_{p'}$ of
functions $f$ on $\R_+$ such that 
$$
\sum_{k\in\Z}2^{-np'/2}\left(\int_{2^n}^{2^{n+1}}|f(x)|^2dx\right)^{p'/2}<\be
$$
with respect to the pairing
\bay
\label{pair}
(\f,f)=\int_0^\be\f(x)f(x)dx.
\ey
Here $p'=p/(p-1)$. 
Suppose that $f$ is a  function on $(0,\be)$. 
Define the function $\psi$ on $\R_+$ by $\psi(x)=\frac{f(x)}{2x}$,
$x>0$.  
It is straightforward
to see from the definition of the $X_p$ spaces 
that  
$f\in X_p^*$ if and only if $\psi\in X_{p'}$ and the norm of $f$ in
$X_p^*$ and the norm of $\psi$ in $X_{p'}$ 
are equivalent. It is also easy to see that for $1<p<\be$ the space
$X_p$ is reflexive. 

Let us show that 
if $f$ is a bounded function in $X_p^*$
with compact support in $(0,\be)$, then
\bay
\label{trTQ}
(\f,f)=\trace TQ_\psi.
\ey
We have
\bey
\trace TQ_\psi&=&\iint\limits_{\R^2_+}k(x,y)\psi(\max\{x,y\})\,dxdy\\[.2cm]
&=&\int\limits_0^\be\psi(x)
 \left(\int\limits_0^xk(x,t)dt+\int\limits_0^xk(s,x)ds\right)\\[.2cm]
&=&\int\limits_0^\be2x\psi(x)\f(x)\,dx=\int\limits_0^\be\f(x)f(x)dx.
\eey

It follows that
\bey
\sup\{|(\f,f)|:~f\in X^*_p,~\|f\|_{X_p^*}\le1\}
&\le&\const\|T\|_{\bS_p}
 \sup\{\|Q_\psi\|_{\bS_{p'}}:~\|\psi\|_{X_{p'}}\le1\}\\[.2cm]
&\le&\const\|T\|_{\bS_p}
\eey
by Theorem \ref{T:p>1}. It follows that $\f\in X_p$, and again by
Theorem \ref{T:p>1}, ${\cal Q}T\in\bS_p$. 
$\bl$

Theorem \ref{1p2} allows us to define for $1<p<2$ the adjoint operator
${\cal Q}^*$ on $\bS_{p'}$. 
Since ${\cal Q}$ is an orthogonal projection on $\bS_2$, $\cal Q$ is a
self-adjoint operator on $\bS_2$. 
We denote the adjoint operator ${\cal Q}^*$ on $\bS_{p'}$ by the same
symbol $\cal Q$. 

Thus we can consider the projection $\cal Q$ on any class $\bS_p$ with
$1<p<\be$. 
It is easy to show that if $T$ is an integral
operator with kernel function $k$ of class $\bS_p$, $2<p<\be$, then
${\cal Q}T=Q_\f$, 
where $\f$ is defined by \rf{fik} and $Q_\f\in\bS_p$. We are going to
prove that for any $T\in\bS_p$, 
$2<p<\be$, the operator ${\cal Q}T$ has the form $Q_\f$ for a function
$\f\in X_p$. 

\begin{thm}
\label{p>2}
Let $T$ be an operator of class $\bS_p$, $2<p<\be$. Then there exists
a function $\f\in X_p$ 
such that ${\cal Q}T=Q_\f$.
\end{thm}

\Pf Let $\bsX_p$ be the space of operators of the form $Q_\f$ with
$\f\in X_p$. Clearly, $\bsX_p$ is 
a Banach space with norm
$$
\|Q_\f\|_{\bsX_p}=\|\f\|_{X_p}.
$$
It follows from Theorems \ref{1p2} and \ref{T:p>1}
that for $T\in\bS_2$ 
$$
\|{\cal Q}T\|_{\bsX_p}\le\const\|{\cal Q}T\|_{\bS_p}\le\const\|T\|_{\bS_p}.
$$
Since $\bS_2$ dense in $\bS_p$, it follows that ${\cal Q}T\in \bsX_p$ for
any $T\in\bS_p$. 
\qed

We consider now the behavior of $\cal Q$ on $\bS_1$. It follows from
Theorem \ref{1p2} that  
if $T\in\bS_1$, then ${\cal Q}T\in\bS_p$ for any $p>1$. The next
result shows that ${\cal Q}T$ does not have 
to be in $\bS_1$ but it has to be in $\bS_{1,\be}$.

\begin{thm}
\label{sltip}
\begin{enumerate}
\item
There exists an operator $T$ in $\bS_1$ such that ${\cal Q}T\notin\bS_1$.
\item
$\cal Q$ has weak type $(1,1)$, i.e.,  $\cal Q$ maps $\bS_1$ into
$\bS_{1,\be}$, i.e., 
$$
s_n({\cal Q}T)\le\const(1+n)^{-1}\|T\|_{\bS_1},\quad T\in\bS_1.
$$
\end{enumerate}
\end{thm}

\begin{lem}
\label{X1}
$$
{\cal Q}\bS_1=\bsX_1\df\{Q_\f:~\f\in X_1\}.
$$
\end{lem}

Let us first deduce Theorem \ref{sltip} from Lemma \ref{X1}.

{\bf Proof of Theorem \ref{sltip}.} (i) is an immediate consequence of
Lemma \ref{X1} and the Example following 
Theorem \ref{E:exp}. (ii) also follows immediately from Lemma \ref{X1}
and Theorem \ref{T:ws1}. $\bl$ 

{\bf Proof of Lemma \ref{X1}.} Let us first show that ${\cal
Q}\bS_1\subset\bsX_1$. Let $T\in\bS_1$ and 
${\cal Q}T=Q_\f$. We have to prove that $\f\in X_1$. Consider the
space $Z^0_\be$ that consists of 
functions $f$ on $\R_+$ such that 
$$
\lim_{n\to\pm\be}\left(2^{-n}\int_{2^n}^{2^{n+1}}|f(x)|^2dx\right)=0.
$$
It is not difficult 
to see that $(Z^0_\be)^*=X_1$ with respect to the pairing
\eqref{pair}. As in the proof of Theorem 
\ref{1p2} we define the function $\psi$ by $\psi(x)=\frac{f(x)}{2x}$,
$x>0$. It follows from \rf{trTQ} that 
$$
|(\f,f)|\le\const\|T\|_{\bS_1}\|Q_\psi\|
\le\const\|T\|_{\bS_1}\|\psi\|_{X_\be^0}
\le\const\|T\|_{\bS_1}\|f\|_{Z^0_\be},
$$
and so $\f$ determines a continuous linear functional on
$Z^0_\be$. Hence, $\f\in X_1$. 

To prove that ${\cal Q}\bS_1=\bsX_1$, we consider the operator
$A:\bS_1\to X_1$ defined by 
$AT=\f$, where $\f$ is the function on $\R_+$ such that ${\cal
Q}T=Q_\f$. We have to show that  
$A$ maps $\bS_1$ onto $X_1$. Consider the conjugate operator
$A^*:X_1^*\to\B(L^2(\R_+))$. 

It is easy to see that with respect to the pairing \rf{pair} the space
$X_1^*$ can be identified with the space 
$Z_\be$ that consists of functions $f$ on $\R_+$ such that 
$$
\sup_{n\in\Z}2^{-n/2}\left(\int_{2^n}^{2^{n+1}}|f(x)|^2\right)^{1/2}<\be.
$$
Consider the operator $J:X_1^*\to X_\be$ defined by
$(Jf)(x)=\frac{f(x)}{2x}$, $x>0$. It is easy to 
see that $J$ maps isomorphically $X_1^*$ onto $X_\be$. 

It can easily be verified that $A^* f=Q_{Jf}$. It follows from Theorem
\ref{bou} that  
$\|A^*f\|\ge\const\|f\|_{X_1^*}$. It follows that $A$ maps $\bS_1$
onto $X_1$. $\bl$ 

\begin{remark}
In \cite{Pel3} metric properties of the averaging projection $\p$ onto the
space of Hankel matrices were studied. In particular, it was shown in
\cite{Pel3} that $\p\bS_1\subset\bS_{1,2}$. 
However, it turns out that the averaging projection $\cal Q$ onto the
operators $Q_\f$ has different 
properties. Theorem \ref{sltip} shows that 
${\cal Q}\bS_1\subset\bS_{1,\be}$ but it follows from Lemma  
\ref{X1} and the remark preceding Theorem \ref{zero}
that ${\cal Q}\bS_1\not\subset\bS_{1,q}$ for any $q<\be$.
\end{remark}

Recall that the Matsaev ideal $\bS_\o$ consists of the operators $T$ on
Hilbert space such that 
$$
\|T\|_{\bS_\o}\df\sum_{n\ge0}\frac{s_n(T)}{1+n}<\be.
$$
It is easy to see that $\bS_p\subset\bS_\o$ for any $p<\be$.

Consider now the operator ideal $\bS_\O$ 
defined by \eqref{SO}. 
It is easy to see that $\bS_{1,\be}\subset\bS_\O$. 
It is well known (see \cite{GK1}) that  
$\bS_\o^*=\bS_\O$ with respect to the pairing
\bay
\label{trpai}
\{T,R\}=\trace TR,\quad T\in\bS_\o,\quad R\in\bS_\O.
\ey

\begin{thm}
\label{mats}
The averaging projection $\cal Q$ defined on $\bS_2$ extends to a
bounded linear operator from $\bS_\o$ 
to the space of compact operators. If \mbox{$T\in\bS_\o$}, then
\linebreak${\cal Q}T=Q_\f$ for a  
function $\f$ in $X_\be^0$.
\end{thm}

\Pf Let us prove that $\cal Q$ extends to a bounded operator from
$\bS_\o$ to the space of compact 
operators. The proof of the second part of the theorem is the same as
the proof of Theorem \ref{p>2}. 
Since the finite rank operators are dense in $\bS_\o$, it is
sufficient to show that $\cal Q$ extends 
to a bounded operator from $\bS_\o$ to $\B(L^2(\R_+))$.

Let $T\in\bS_2$ and $R\in\bS_1$. By Theorem \ref{sltip}, 
${\cal Q}R\in\bS_{1,\be}\subset\bS_\O$. We have 
$$
\{{\cal Q}T,R\}=\{T,{\cal Q}R\},
$$
and so
\bey
|\{{\cal Q}T,R\}|&\le&\const\|T\|_{\bS_\o}\|{\cal Q}R\|_{\bS_\O}\\[.2cm]
&\le&\const\|T\|_{\bS_\o}\|{\cal Q}R\|_{\bS_{1,\be}}
\le\const\|T\|_{\bS_\o}\|R\|_{\bS_1}
\eey
by Theorem \ref{sltip}. Since $\bS_1^*=\B(L^2(\R_+))$ with respect to
the pairing \rf{trpai}, it follows that 
$\|{\cal Q}T\|\le\const\|T\|_{\bS_\o}$, and so $\cal Q$ extends to a
bounded linear operator  
from $\bS_\o$ to $\B(L^2(\R_+))$. $\bl$

\

\setcounter{equation}{0}
\section{\bf Finite rank}
\label{S:finite}

\

We say that $\f$ is a \emph{step function} if there exist finitely many numbers
\linebreak$0=x_0<x_1<\dots<x_N<\be$ such that $f$ is \aex{} constant
on each interval 
$(x_{i-1},x_i)$, and zero on $(x_N,\be)$. 
The number of steps of $\f$ then is the smallest possible $N$ in this
definition. 

There is a natural correspondence between operators $\qf$ where the
symbol $\f$ is a step function, with given $x_1<\dots<x_N$,
and matrices of the form $\{a_{\max\{i,j\}}\}$. 
We need a simple result for
such matrices, but will not pursue their study further.

\begin{lem}
\label{L:kronecker}
If $a_1,\dots,a_n$ are complex numbers, then the matrix \linebreak
$\{a_{\max\{i,j\}}\}_{1\le i,j\le n}$ has determinant 
$a_n\prod_{i=1}^{n-1} (a_{i}-a_{i+1})$.
\end{lem}

\Pf
Denote this determinant by $\D(a_1,\dots,a_n)$. Subtracting the last
row from all others, we see that 
$D(a_1,\dots,a_n)=a_nD(a_1-a_n,\dots,a_{n-1}-a_n)$, and the result
follows by induction.
\qed

\begin{thm}
\label{T:Kronecker}
$\qf$ has finite rank if and only if $\f$ is a step function.
In this case, the rank of $\qf$ equals the number of steps of $\f$.
\end{thm}

\Pf
If $\f$ is a step function with $N$ steps, we have 
$\f=\sum_1^N a_i\chi_{(0,x_i)}$ for some $a_i$ and $x_i>0$, and thus
$\qf$ is a linear combination of $N$ rank one operators.

Conversely, suppose that $\rank(\qf)=M<\infty$.
Suppose that $n>M$ and that $z_1<\dots<z_n$ are Lebesgue points of
$\f$.
If $h>0$ and $\fzh=h^{-1}\chi_{(z,z+h)}$, then the
matrix $(\inner{\qf f_{z_i,h},f_{z_j,h}})_{ij}$ has rank at most $M<n$
and thus its determinant vanishes. 
As $h\to0$, as shown in the proof of \refT{T:positive},
$\inner{\qf f_{z_i,h},f_{z_j,h}}\to \f(\max\{z_i,z_j\})=\f(z_{\max\{i,j\}})$,
and thus the determinant of 
$\bigpar{\f(z_{\max\{i,j\}})}_{ij}$ vanishes too.
By \refL{L:kronecker}, this implies that either $\f(z_i)=\f(z_{i+1})$
for some $i<n$ or $\f(z_n)=0$.

Consequently, if 
$z_1<\dots<z_n$ are Lebesgue points of
$\f$ such that  $\f(z_i)\neq\f(z_{i+1})$
for $i<n$ and $\f(z_n)\neq0$, then $n\le M$. Choose such a sequence 
$z_1<\dots<z_n$ with $n$ maximal. If $z\in(z_i,z_{i+1})\cap\Leb(\f)$
for some $i<n$, then either $\f(z)=\f(z_i)$ or $\f(z)=\f(z_{i+1})$,
since $n$ is maximal. Moreover, for the same reason,
if $\f(z)=\f(z_i)$, then
$\f(z')=\f(z_i)$ for every $z'\in(z_i,z)\cap\Leb(\f)$, and 
if $\f(z)=\f(z_{i+1})$, then
$\f(z')=\f(z_{i+1})$ for every $z'\in(z,z_{i+1})\cap\Leb(\f)$.
Together with similar arguments for the intervals $(0,z_1)$ and
$(z_n,\be)$, which we leave to the reader, this easily shows that $\f$
is a step function with at most $n$ steps.
\qed

\

\setcounter{equation}{0}
\section{\bf A class of integral operators on $\bs{L^2(\R)}$}
\label{io}

\

\newcommand{\hs}{\heartsuit}

In this section we associate with the operator $Q_\f$ on $L^2(\R_+)$
an integral 
operator on $L^2(\R)$ and we study these operators.

For a function $\f\in \Lloc^2(\R_+)$ we define 
the function $\f^\hs$ on $\R$ by 
\bay
\label{ser}
\f^\hs(t)\df2\f(e^{2t})e^{2t},\quad t\in\R.
\ey
With a function $\psi$ on $L^2(\R)$ we associate the function
$\breve\psi$ on $\R\times\R$ defined 
by
\bay
\label{brpsi}
\breve\psi(s,t)=\psi(\max\{s,t\})e^{-|s-t|},\quad s,t\in\R,
\ey
and denote by $K_\psi$ the integral operator on $L^2(\R)$ with kernel
function $\breve\psi$ 
(if it makes sense):
$$
(K_\psi f)(s)=\int_\R\psi(\max\{s,t\})e^{-|s-t|}f(t)dt.
$$

\begin{thm}
\label{pr1}
Let $\f\in \Lloc^2(\R_+)$. Then
the operators $Q_\f$ and $K_{\f^\hs}$ are unitarily equivalent.
\end{thm}

Theorem \ref{pr1} certainly means that the boundedness of one of the
operators implies 
the boundedness of the other one.

\Pf Consider the unitary operator $U:L^2(\R_+)\to L^2(\R)$ defined as
follows $(Uf)(t)=\sqrt2f(e^{2t})e^t$.  It remains to observe that
$K_{\f^\hs} U=UQ_\f$.  $\bl$ 

We can identify $L^2(\R_+)$ with the subspace of $L^2(\R)$ which
consists of the functions 
vanishing on $(-\be,0)$. We can
now extend in a natural way the
operator of triangular 
projection $\cp$ to act on the space of operators on $L^2(\R)$ by
defining it in the same 
way as it has been done in \refS{prel}. 
We keep the same notation, $\cp$, for this extension.
We put $K_\psi^+\df\cp K_\psi$ and  
$K_\psi^-\df K_\psi-\cp K_\psi^+$.

It is easily seen from the proof of Theorem \ref{pr1} that the
operators $Q_\f^+$ and 
$Q_\f^-$ are unitarily equivalent to the operators $K_{\f^\hs}^+$ and
$K_{\f^\hs}^-$ respectively.

It is easy to see that, for any $p>0$, 
\begin{equation}\label{xpr}
\f\in X_p\eq\sum_{n=-\be}^\be \norm{\f^\hs}_{L^2[n,n+1]}^p<\be
\end{equation}
(and correspondingly for $X_\be$ and $X_\be^0$) and, 
using \refT{T:yp}\vvp(vii),
\begin{equation}\label{ypr}
\f\in Y_p\eq\sum_{n=-\be}^\be \norm{\f^\hs}_{BV[n,n+1]}^p<\be.
\end{equation}
(In \eqref{xpr} and \eqref{ypr},
the intervals $[n,n+1]$ can be replaced by any partition of $\R$ into
intervals of the same length.) 
We can thus translate results from the preceding sections to $K_\psi$,
for example as follows.
\begin{thm}
\label{bou1}
Let $\psi\in L^2_{\textup{loc}}(\R)$. 
The following are equivalent:
\begin{enumerate}
\item 
$K_\psi$ is bounded on $L^2(\R)$;
\vspace*{.3cm}
\item 
$K^+_\psi$ is bounded on $L^2(\R)$;
\vspace*{.2cm}
\item
$\displaystyle\sup_{n\in\Z}\int_n^{n+1}|\psi(x)|^2dx<\be$.
\end{enumerate}
\end{thm}

\Pf The theorem is a direct 
consequence of Theorems \ref{bou} and \ref{pr1}.
$\bl$

Similarly we find from Theorems \ref{T:p>1},
\ref{T:x1},
and \ref{T:yp},
respectively, the following three theorems.

\begin{thm}
Let $\psi\in L^2_{\textup{loc}}(\R)$. 
If $1<p<\be$,  the following are equivalent:
\begin{enumerate}
\item 
$K_\psi\in\bS_p$.
\item 
$K^+_\psi\in\bS_p$.
\vspace*{.2cm}
\item
$\sum_{n\in\Z}\norm{\psi}_{L^2[n,n+1]}^p<\be$. \qed
\end{enumerate}
\end{thm}

\begin{thm}
\label{T:S1R}
If $K_\psi\in\bS_1$,
then $\sum_{n\in\Z} \norm{\psi}_{L^2[n,n+1]}<\be$
and thus $\psi\in L^1(\R)$.
Moreover, then
$
\trace K_\psi=\int_{-\be}^\be\psi(x)\,dx
$. 
\qed
\end{thm}

\begin{thm}
If $1/2<p\le1$ and
$\sum_{n\in\Z} \norm{\psi}_{BV[n,n+1]}^p<\be$,
then
$K_\psi\in\bS_p$.
\qed
\end{thm}

The following two results involving the modulus of continuity
also can be obtained by changes of variables in the corresponding
Theorems \refand{T:omegafi}{T:Dinip}, 
using \eqref{omega2} and \refL{L:moduli}, but the details are 
involved and we prefer to imitate the proofs.

\begin{thm}
\label{T:omegafiR}
If $\psi\in L^2(\R)$ has support on $\OI$, then
\begin{equation*}
s_n(K_\psi)\le C \frac{1} n \oxii{\psi}\!\left(\frac{1}n\right) 
+ C \frac{1}{n^2}\norm{\psi}_{L^2}, 
\qquad n\ge1.
\end{equation*}
\end{thm}

\Pf 
We interpolate using \refL{L:K} as in the proof of
\refT{T:omegafi}.
$\bl$

\begin{thm}
\label{T:DinipR}
Let $1/2<p\le1$. If\/
$\sum_{n\in\Z} \norm{\psi}_{L^2[n,n+1]}^p<\be$
and
$$
\sum_{n=-\be}^\be \int_0^{1} 
 \bigpar{\oxii{\psi}\bigpar{t;[n,n+1]}}^p t^{p-2}\,dt <\be,
$$
then
$K_\psi\in\bS_p$.
\end{thm}

\begin{proof}
We argue as in the proofs of Theorems
\refand{T:Dini2}{T:Dinip}, using \refT{T:omegafiR}.
\end{proof}

We denote by $\F$ the Fourier transformation on $L^2(\R^n)$,
which is a unitary operator defined by \eqref{fourier} 
for $f\in L^2(\R^n)\bigcap L^1(\R^n)$.  Let $T$ be the integral
operator with kernel function $k\in L^2(\R^2)$. 
Denote by $R$ the integral operator with kernel function $\F k$. 
The following lemma has a straightforward verification.
\begin{lem}
\label{fr}
$\F T\F=R$. $\bl$ 
\end{lem}

\begin{cor} 
\label{eqp}
Let $p>0$. Then $\|T\|_{\bS_p}=\|R\|_{\bS_p}$. $\bl$
\end{cor}

Denote by $Z$ be the integral operator with kernel function
$(x,y)\mapsto(\F k)(x,-y)$. It is easy 
to see that $T$ is unitarily equivalent to $Z$. Indeed, 
the equality $\F T\F=R$ implies $\F T\F^{-1}=Z$.

\begin{lem} 
\label{tf}
If $\psi\in L^1(\R)$, then $\breve\psi\in L^1(\R^2)$, and
\bay
\label{Fbrpsi}
(\F\breve\psi)(x,y)=
(\F\psi)(x+y)\left(\dfrac1{1-2\pi{\rm i}x}+\dfrac1{1-2\pi{\rm i}y}\right).
\ey
\end{lem}

\Pf The inclusion $\breve\psi\in L^1(\R^2)$ is obvious. We have 
\begin{align*}
(\F\breve\psi)(x,y)&=
\int\limits_{\R^2}\psi(\max\{s,t\})e^{-|s-t|}
  e^{-2\pi{\rm i}sx-2\pi{\rm i}ty}dsdt\\[.2cm] &=
\int\limits_{\R}\left(\psi(s)e^{-s}e^{-2\pi{\rm i}sx}
\int\limits_{-\infty}^se^{t-2\pi{\rm i}ty}dt\right)ds\\[.2cm] 
&\qquad+
\int\limits_{\R}\left(\psi(t)e^{-t}e^{-2\pi{\rm i}ty}
\int\limits_{-\infty}^te^{s-2\pi{\rm i}sx}ds\right)dt\\[.2cm] &=
\int\limits_{\R}\frac{\psi(s)e^{-2\pi{\rm i}s(x+y)}\,ds}{1-2\pi{\rm i}y}+
\int\limits_{\R}\frac{\psi(t)e^{-2\pi{\rm i}t(x+y)}\,dt}{1-2\pi{\rm i}x}
\\[.2cm]
&= \left(\frac1{1-2\pi{\rm i}x}+\frac1{1-2\pi{\rm i}y}\right)(\F\psi)(x+y).
\quad \bl
\end{align*}

Consider the functions
\bay
\label{+-}
\breve\psi_+(s,t)\df\chi_{\{(s,t):s>t\}}\breve\psi(s,t)\quad\mbox{and}\quad
\breve\psi_-(s,t)\df\chi_{\{(s,t):s<t\}}\breve\psi(s,t).
\ey
It can easily be seen from the proof of Lemma \ref{tf} that 
\bay
\label{F+-}
(\F\breve\psi_+)(x,y)=\frac{(\F\psi)(x+y)}{1-2\pi{\rm i}y}\quad\mbox{and}\quad
(\F\breve\psi_-)(x,y)=\frac{(\F\psi)(x+y)}{1-2\pi{\rm i}x}.
\ey

It is easy to verify 
that if $\psi$ is a tempered distribution on $\R$
(see \refS{prel}), 
we can
define tempered distributions $\breve\psi$,
$\breve\psi_+$, and $\breve\psi_-$ by \rf{brpsi} and \rf{+-};
the formal definitions are by duality and analogous to
\eqref{lp+}--\eqref{lp}. 
It is
easy to check 
that 
formulas \rf{Fbrpsi} and \rf{F+-} also hold for tempered distributions $\psi$.

As a corollary to \refT{bou1}, we have the following lemma. 
\begin{lem}
\label{L:psitd} 
Let {\em$\psi\in L^2_{\text{loc}}(\R)$}. Suppose that the operator
$K_\psi$ is bounded on $L^2(\R)$. Then $\psi$ determines a tempered
distribution. $\bl$
\end{lem}

\begin{thm}
\label{boup}
Let $0<p<\be$. Suppose that
$\f\in L^2_{\textup{loc}}(\R_+)$ and $\f^\hs$ is defined by \eqref{ser}. 
The following are equivalent:
\begin{enumerate}
\item
$Q_\f\in\bS_p$;
\vspace*{.3cm}
\item
$K_{\f^\hs}\in\bS_p$;
\vspace*{.2cm}
\item
the integral operator on $L^2(\R)$ with kernel 
$$
(x,y)\mapsto(\F\f^\hs)(x+y)
 \left(\dfrac1{1-2\pi{\rm i}x}+\dfrac1{1-2\pi{\rm i}y}\right)
$$
belongs to $\bS_p$.
\end{enumerate}
\end{thm}

\Pf The theorem is a consequence of Theorem \ref{pr1}, Lemma \ref{fr}
and Lemma \ref{tf}. $\bl$

Note that if in (iii) we have a tempered distribution rather than a
function, by the integral operator 
we mean the operator determined by this tempered distribution (see 
\refS{prel}).  

In the same way one can prove the following result.
\begin{thm}
\label{boup1}
Let $0<p<\be$ and let $\f\in L^2_{\textup{loc}}(\R_+)$. 
The following are equivalent:
\begin{enumerate}
\item
$Q^+_\f\in\bS_p$;\vspace*{.3cm}
\item
$Q^-_\f\in\bS_p$;\vspace*{.3cm}
\item
$K^+_{\f^\hs}\in\bS_p$;\vspace*{.3cm}
\item
$K^-_{\f^\hs}\in\bS_p$;\vspace*{.2cm}
\item
the integral operator on $L^2(\R)$ with the kernel 
$$
(x,y)\mapsto\frac{(\F\f^\hs)(x+y)}{1-2\pi{\rm i}x}
$$ 
belongs to $\bS_p$. $\bl$
\end{enumerate}
\end{thm}

It is straightforward to show that if $p>0$ and the integral operator on $L^2(\R)$
with kernel function 
$\frac{h(x+y)}{x+\a}$ belongs to $\bS_p$ for some
$\a\in\C\setminus\R$, then it belongs to $\bS_p$ 
for any $\a\in\C\setminus\R$. Let us show that such an integral operator can
belong to $\bS_p$ for $p\le1$ only if it is zero.

Consider the operator $E:\cd(\R^2)\to\cd(\R_+\times\R_+)$ defined by
the following 
formula 
$(Ef)(s,t)\df\frac12(st)^{-1/2}f(\frac12\log s,\frac12\log t)$. 
Clearly, $E$ is an 
isomorphism. Consequently, the conjugate operator $E\,^\prime$ is an
isomorphism 
from $\cd\,^\prime(\R_+\times\R_+)$ onto
$\cd\,^\prime(\R^2)$. Clearly, $(E\,^\prime\Phi)(x,y)= 
2\Phi(e^{2x},e^{2y})e^xe^y$ if $\Phi\in L^1_{\textup{loc}}(\R_+\times\R_+)$.
Put $2\Phi(e^{2x},e^{2y})e^xe^y\df (E\,^\prime\Phi)(x,y)$ for 
$\Phi\in\cd\,^\prime(\R_+\times\R_+)$.

\begin{thm}
\label{ueq}
Let $\Phi\in\cd\,^\prime(\R_+\times\R_+)$. Then $\Phi$ determines a
bounded operator 
on $L^2(\R_+)$ if and only if the distribution $2\Phi(e^{2x},e^{2y})e^xe^y$ 
determines a bounded operator on
$L^2(\R)$. Moreover, these two operators are unitarily equivalent operators.
\end{thm}
\Pf It suffices to note that
$$
\langle2\Phi((e^{2x},e^{2y})e^xe^y,\sqrt2e^yf(e^{2y})
 \overline{\sqrt2e^xg(e^{2x})}\rangle=
\langle\Phi(s,t),f(t)\overline{g(s)}\rangle 
$$
for any $f,g\in\cd(\R_+\times\R_+)$, and the map
$h\mapsto\sqrt2e^xh(e^{2x})$ is a unitary 
operator from $L^2(\R_+)$ onto $L^2(\R)$.
$\bl$

\begin{thm}
\label{s10}
Let $h\in\cd\,^\prime(\R)$. Suppose that the distribution 
$\frac{h(x+y)}{1-2\pi{\rm i}x}$
determines an operator on $L^2(\R)$ of class $\bS_1$. Then $h=0$.
\end{thm}
\Pf By Lemma \ref{temp}, $\frac{h(x+y)}{1-2\pi{\rm i}x}\in{\cal
S}^\prime(\R^2)$. 
Consequently, $h(x+y)\in{\cal S}^\prime(\R^2)$, whence $h\in{\cal
S}^\prime(\R)$. 
Thus, there exists a distribution $\f\in\cd\,^\prime(\R_+)$ such that 
\linebreak$\f^\hs\in{\cal S}^\prime(\R)$ and $\F\f^\hs=h$
(the operation $\f\mapsto\f^\hs$ defined in \rf{ser} extends in an
obvious way to distributions $\f$).  
Lemma \ref{fue} and formula
\rf{F+-} imply that $Q_\f^-$ and thus also
$Q_\f^+$ belongs to $\bS_1$. 
By \refT{ker+}, $\f\in L^2_{\textup{loc}}(\R_+)$. 
Thus $\f=0$ by Theorem \ref{zero}.
$\bl$

We are going to prove now that for $p>1/2$ if the integral 
operator with kernel function 
$h(x+y)\bigl(\frac1{x+\a}+\frac1{x+\b}\bigr)$ belongs to $\bS_p$ for 
some $\a,\,\b\in\C\setminus\R$ with $\a+\b\not\in\R$, then it belongs to
$\bS_p$ for any $\a,\,\b\in\C\setminus\R$. We will also show that
this is not true for $p\le1/2$.

\begin{lem}
\label{anm}
The function 
$$
(x,y)\mapsto\frac1{x+y+{\rm i}}\,\chi_{[0,1]}(x),\quad x,\,y\in\R,
$$ 
is a Schur multiplier of $\bS_p$ for any $p>0$.
\end{lem}

\Pf First we prove that the function
\bay
\label{[2,2]}
(x,y)\mapsto\frac1{x+y+{\rm i}}\,\chi_{[0,1]}(x)\chi_{\R\setminus[-2,2]}(y),
\quad x,\,y\in\R,
\ey
is a Schur multiplier of $\bS_p$. We have
$$
\frac1{x+y+{\rm i}}\,\chi_{[0,1]}(x)\chi_{\R\setminus[-2,2]}(y)
=\sum\limits_{n\ge0} (-1)^n\frac{(x+{\rmi})^n}{y^{n+1}}
\chi_{[0,1]}(x)\chi_{\R\setminus[-2,2]}(y).
$$
Clearly, the $p$-multiplier norm of the $n$-th summand 
is bounded by 
$2^{-\xfrac n2}$. Consequently, the function \rf{[2,2]}
is a Schur multiplier of $\bS_p$. It remains to prove that the function
$$
(x,y)\mapsto\frac1{x+y+{\rm i}}\,\chi_{[0,1]}(x)\chi_{[-2,2]}(y),
\quad x,\,y\in\R,
$$
is a Schur multiplier of $\bS_p$.
For any $(\xi,\eta)\in[0,1]\times[-2,2]$
we can expand the function $\frac1{x+y+{\rm i}}$ in a Taylor series in a
neighborhood  
of $(\xi,\eta)$. It follows easily that for a sufficiently small
$\e>0$ the function 
$$
(x,y)\mapsto\frac1{x+y+{\rm i}}\,\chi_{[0,1]
 \cap[\xi-\e,\xi+\e]}(x)\chi_{[-2,2]\cap[\eta-\e,\eta+\e]}(y),
\quad x,\,y\in\R,
$$ 
is a Schur multiplier of $\bS_p$. It remains to choose
a finite subcover of $[0,1]\times[-2,2]$ that consists of rectangles
of the form 
$[\xi-\e,\xi+\e]\times[\eta-\e,\eta+\e]$.
$\bl$

\begin{remark}
In the same way we can prove that the function 
$$
(x,y)\mapsto\frac1{x+y+\a}\,\chi_{[\xi,\eta]}(x),\quad x,\,y\in\R,
$$ 
is a Schur multiplier of $\bS_p$ for any $p>0$ for any
$\a\in\C\setminus\R$ and for any 
$\xi,\,\eta\in\R$.
\end{remark}

\begin{cor} 
\label{anm1}
Let $\a,\,\b,\,\g\in\C$ and $\g\not\in\R$. Then the function 
$$
(x,y)\mapsto\frac{(x+\a)(y+\b)}{x+y+\g}\chi_{[0,1]}(x),\quad x,\,y\in\R,
$$
is a Schur multiplier of $\bS_p$ for any $p>0$.
\end{cor}

\Pf We have 
$$
\frac{(x+\a)(y+\b)}{x+y+\g}=(x+\a)\left(1+\frac{\b-\g-x}{x+y+\g}\right).
$$
It remains to note that 
$$
\chi_{[0,1]}(x),\quad(x+\a)\chi_{[0,1]}(x),\quad(\b-\g-x)\chi_{[0,1]}(x),
 \quad\mbox{and}\quad \frac1{x+y+\g}\chi_{[0,1]}(x)
$$ 
are Schur multipliers of $\bS_p$.
$\bl$

\begin{cor} 
\label{anm2}
Let $p>0$ and let $\a,\,\b\in\C\setminus\R$ such that
$\a+\b\not\in\R$. Suppose that the integral operator 
on $L^2(\R)$ with kernel function
$$
(x,y)\mapsto h(x+y)\left(\frac1{x+\a}+\frac1{y+\b}\right),\quad x,\,y\in\R,
$$ 
belongs to $\bS_p$. Then the integral operator with kernel function
$$
(x,y)\mapsto h(x+y)\chi_{[0,1]}(x),\quad x,\,y\in\R,
$$ 
belongs to $\bS_p$. $\bl$
\end{cor}

\begin{thm}
\label{del}
Let $p>1/2$ and let $\a_0,\,\b_0\in\C\setminus\R$ such that
$\a_0+\b_0\not\in\R$. Suppose that the 
integral operator on $L^2(\R)$ with kernel function 
$$
(x,y)\mapsto h(x+y)\left(\frac1{x+\a_0}+\frac1{y+\b_0}\right),\quad x,\,y\in\R,
$$ 
belongs to $\bS_p$. Then the integral operator with kernel function
$$
(x,y)\mapsto h(x+y)\left(\frac1{x+\a}+\frac1{y+\b}\right),\quad x,\,y\in\R,
$$ 
also belongs to $\bS_p$ for any $\a,\,\b\in\C\setminus\R$.
\end{thm}

\Pf By Corollary \ref{anm2}, the integral operator
with kernel $h(x+y)\chi_{[0,1]}(x)$ belongs to $\bS_p$. Obviously, for
any $n\in\Z$, 
$$
\|h(x+y)\chi_{[0,1]}(x)\|_{\bS_p}=\|h(x+y)\chi_{[n,n+1]}(x)\|_{\bS_p}
$$ 
(as usual we write $\|k\|_{\bS_p}$ for the $\bS_p$ norm (or quasi-norm)
of the integral operator with 
kernel $k$). Consequently,
$$
\left\|h(x+y)\left(\frac1{x+\a}-\frac1{x+\a_0}\right)
 \chi_{[n,n+1]}(x)\right\|_{\bS_p}
\le\const(1+|n|)^{-2}.
$$ 
It is now clear that the integral operator with kernel function
$$
(x,y)\mapsto h(x+y)\left(\frac1{x+\a}-\frac1{x+\a_0}\right),\quad x,\,y\in\R,
$$ 
belongs to $\bS_p$ for $p>1/2$. Similarly, we prove that the integral
operator with kernel function 
$$
(x,y)\mapsto h(x+y)\left(\frac1{y+\b}-\frac1{y+\b_0}\right),\quad x,\,y\in\R,
$$
belongs to $\bS_p$ for $p>1/2$. 
$\bl$

\begin{thm}
\label{pow}
Let $\f\in L^2_{\rm{loc}}(\R_+)$, $a\in\R\setminus\{0\}$, and $p>1/2$.
Put 
$$
\fa(t)\df\f(t^a)t^{a-1}.
$$
Then $Q_\f\in\bS_p$ if and only if 
$Q_{\fa}\in\bS_p$.
\end{thm}

\Pf Recall that $\f^\hs(t)=2\f(e^{2t})e^{2t}$ and
$\fa^\hs(t)=2\fa(e^{2t})e^{2t} =2\f(e^{2a t})e^{2a t}$. 
Consequently, $\fa^\hs(t)=\f^\hs(a t)$.
By Theorem \ref{boup}, $Q_{\fa}\in\bS_p$ if and only if the integral
operator on $L^2(\R)$ with  
kernel 
$$
(x,y)\mapsto
(\F\f^\hs)\left(\frac{x+y}a\right)\left(\dfrac1{1-2\pi{\rm i}x}
 +\dfrac1{1-2\pi{\rm i}y}\right),
\quad x,\,y\in\R,
$$
belongs to $\bS_p$, and  
thus if and only if 
the integral operator on $L^2(\R)$ with 
kernel 
$$
(x,y)\mapsto
(\F\f^\hs)(x+y)\left(\dfrac1{1-2\pi{\rm i} a x}
+\dfrac1{1-2\pi{\rm i} a y}\right),
\quad x,\,y\in\R,
$$
belongs to $\bS_p$.
By Theorem \ref{del}, this holds if and only if 
the integral operator on $L^2(\R)$ with 
kernel 
$$
(x,y)\mapsto
(\F\f^\hs)(x+y)\left(\dfrac1{1-2\pi{\rm i}x}+\dfrac1{1-2\pi{\rm i}y}\right),
\quad x,\,y\in\R,
$$
belongs to $\bS_p$.
It remains to apply Theorem \ref{boup} once more.
$\bl$

\begin{cor}
\label{C:pow}
Let $\psi\in \Lloc^2(\R)$, $a\in\R\setminus\{0\}$, and $p>1/2$.
Define 
\linebreak 
$
\psi_a(t)\df\psi(at)
$.
Then $K_\psi\in\bS_p$ if and only if 
$K_{\psi_a}\in\bS_p$. 
\qed
\end{cor}

\begin{remark}
Theorem \ref{pow} and its corollary 
do not generalize to the case $p\le1/2$.
Indeed, if $\f$ is the characteristic function of an interval, then
$Q_\f\in\bS_p$ 
for any $p>0$ but if $a\ne1$, then $Q_{\fa}\not\in\bS_p$ by Corollary
\ref{moc1}. 
It follows from the proof above that \refT{del} too does not extend to
$p\le1/2$. 
\end{remark}

\begin{remark} 
Note that if $\f\in X_p$, then $\fa(t)\in X_p$ for
any $a\in\R\setminus\{0\}$ and any $p>0$. Moreover, if $\f\in Y_p$, then
$\fa(t)\in Y_p$ for
any $a\in\R\setminus\{0\}$ and any $p>0$. Indeed, let $A>1$.
It is easy to see that
\eqref{xpi} is equivalent to the condition
$$
\sum_{n\in\Z}A^{np/2}
 \left(\int_{A^n}^{A^{n+1}}|\f(x)|^2dx\right)^{p/2}<\be,
$$
while the  
condition in Theorem \ref{T:yp1}\vvp(vii) is equivalent to
$$
\sum\limits_{n\in\Z} \norm{x\f(x)}_{BV[A^n,A^{n+1}]}^{p}<\be,
$$
which easily implies the above assertions.
\end{remark}

\begin{thm}
\label{conv}
Let $\a,\,\b\in\C\setminus\R$ and let $p>0$. Then the integral operator
on $L^2(\R)$ with kernel function
$$
(x,y)\mapsto h(x+y)\left(\frac1{x+\a}+\frac1{y+\b}\right),\quad x,\,y\in\R,
$$
belongs to $\bS_p$ if and only if convolution with the function
$h(x)(x+\a+\b)$ is an operator from $L^2(\R,(1+x^2)\,dx)$ to 
$L^2(\R,(1+x^2)^{-1}\,dx)$ of class $\bS_p$.  
\end{thm}

\Pf Clearly, the integral operator on $L^2(\R)$ with kernel function
$$
(x,y)\mapsto h(x+y)\left(\frac1{x+\a}+\frac1{y+\b}\right),\quad x,\,y\in\R,
$$
belongs to $\bS_p$ if and only if so does the integral operator with
kernel function 
$$
(x,y)\mapsto h(x-y)\frac{x-y+\a+\b}{(x+\a)(y-\b)},\quad x,\,y\in\R.
$$
To complete the proof, it remains to observe that multiplication by
$(x-\b)^{-1}$ is an 
isomorphism from $L^2(\R)$ onto $L^2(\R,(1+x^2)\,dx)$ and
multiplication by $(x+\a)^{-1}$ is an 
isomorphism from $L^2(\R,(1+x^2)^{-1}\,dx)$ onto $L^2(\R)$.
$\bl$

\begin{cor}
\label{conv1}
Let $\a,\,\b,\,\g\in\C\setminus\R$ such that $\a+\b\not\in\R$ and let
$p>1/2$. Then the integral 
operator with kernel function
$$
(x,y)\mapsto h(x+y)\left(\frac1{x+\a}+\frac1{y+\b}\right),\quad x,\,y\in\R,
$$ 
belongs to $\bS_p$ if and only if convolution with the function
$h(x)(x+\g)$ is an operator from $L^2(\R,(1+x^2)\,dx)$ to 
$L^2(\R,(1+x^2)^{-1}\,dx)$ of class $\bS_p$.  
\end{cor}

\Pf It suffices to apply Theorem \ref{del}.
$\bl$

\begin{remark}
In the same way we can prove that the following statements are equivalent
for any $\a\in\C\setminus\R$ and for any $p>0$:
\begin{enumerate}
\item
the integral operator on $L^2(\R)$ with kernel function
$$
(x,y)\mapsto k(x+y)(x+\a)^{-1},\quad x,\,y\in\R,
$$ 
belongs to $\bS_p$;\vspace*{.3cm}
\item
the integral operator on $L^2(\R)$ with kernel function
$$
(x,y)\mapsto k(x+y)(y+\a)^{-1},\quad x,\,y\in\R,
$$ 
belongs to $\bS_p$;\vspace*{.3cm}
\item
convolution with $k$ 
is an operator from $L^2(\R)$ to $L^2(\R,(1+x^2)^{-1}\,dx)$ of class
$\bS_p$;\vspace*{.3cm} 
\item
convolution with $k$
is an operator from $L^2(\R,(1+x^2)\,dx)$ to $L^2(\R)$ of class $\bS_p$.
\end{enumerate}

Let us repeat that Theorem \ref{s10} implies that the integral operator 
with kernel function $k(x+y)(x+\a)^{-1}$ can be a nonzero operator in
$\bS_p$ only if $p>1$. 
\end{remark}

\

\setcounter{equation}{0}
\section{\bf Matrix representation}
\label{S:matrix}

\

Let $\f$ be a function in $L^2_{\text{loc}}(\R)$ such that
$\f(x+1)=\f(x)$, $x\in\R$. 
Consider the operators $Q_\f^{[0,1]}$ and $Q_\f^{[0,2]}$ on $L^2[0,1]$
and $L^2[0,2]$ respectively. 
Obviously, 
\begin{equation}
\label{02a}
\left\|Q_\f^{[0,1]}\right\|_{\bS_p}\le\left\|Q_\f^{[0,2]}\right\|_{\bS_p},
\quad0<p<\be.
\end{equation}
Obviously, 
\bay
\label{02b}
\|Q_\f^{[0,2]}\|\ge\left\|\cp_{[1,2]\times[0,1]}Q_\f^{[0,2]}\right\|
=\|\f\|_{L^2[0,1]}.
\ey
It is also easy to see that
\bay
\label{02c}
\left\|Q_\f^{[0,2]}\right\|_{\bS_p}\le
C(p)\left(\left\|Q_\f^{[0,1]}\right\|_{\bS_p}+\|\f\|_{L^2[0,1]}\right),
\quad0<p<\be,
\ey
where $C(p)$ is a constant that may depend only on $p$.

\begin{thm} 
\label{QfiQpsi}
Let $0<p<\be$. Suppose that $\f$ is a function in $L^2_{\rm{loc}}(\R)$
such that $\f(x)=\f(x+1)$ and  
$\psi(x)\df\f(1-x)$. Then
$$
C_1(p)\left(\left\|Q_\f^{[0,1]}\right\|_{\bS_p}+\left\|Q_\psi^{[0,1]}
 \right\|_{\bS_p}\right)\le
\left\|Q_\f^{[0,2]}\right\|_{\bS_p}\le 
C_2(p)\left(\left\|Q_\f^{[0,1]}\right\|_{\bS_p}
 +\left\|Q_\psi^{[0,1]}\right\|_{\bS_p}\right).
$$
\end{thm}

\Pf To prove the left inequality, consider the integral operator $K$
on $L^2[0,1]$ with 
kernel function
$$
(x,y)\mapsto\f(\min\{x,y\}).
$$ 
Clearly, the operators $K$ and $Q_\psi^{[0,1]}$ are unitarily
equivalent, and so \linebreak 
$\|K\|_{\bS_p}=\bigl\| Q_\psi^{\OI}\bigr\|_{\bS_p}$ for any $p>0$. Note that
$$
\f(\min\{x,y\})+\f(\max\{x,y\})=\f(x)+\f(y).
$$
Hence, $K+Q_\f^{[0,1]}$ is the integral
operator with kernel function $(x,y)\mapsto\f(x)+\f(y)$. Thus 
$\bigl\|Q_\f^{\OI}+K\bigr\|_{\bS_p}\le C(p)\|\f\|_{L^2}$. 
Now the left inequality is
obvious. To prove the right inequality, we have to show that
$$
\|\f\|_{L^2[0,1]}\le
C(p)\left(\left\|Q_\f^{[0,1]}\right\|_{\bS_p}
 +\left\|Q_\psi^{[0,1]}\right\|_{\bS_p}\right).
$$   
Clearly, 
$$
\left(K+Q_\f^{[0,1]}\right)1=\f+\int\limits_0^1\f(t)\,dt.
$$ 
It follows that
$$
\left\|\left(K+Q_\f^{[0,1]}\right)1\right\|_{L^2[0,1]}^2=\|\f\|_{L^2[0,1]}^2+
3\left|\int\limits_0^1\f(t)\,dt\right|^2\ge
\|\f\|_{L^2[0,1]}^2.
$$
Thus, 
$$
\left\|K+Q_\f^{[0,1]}\right\|_{\bS_p}\ge
\left\|K+Q_\f^{[0,1]}\right\|\ge\|\f\|_{L^2[0,1]},
$$
and so
\bey
\|\f\|_{L^2[0,1]}
&\le&
\left\|\left(K+Q_\f^{[0,1]}\right)\right\|_{\bS_p}\\[.2cm]
&\le& 
C(p)\left(\|K\|_{\bS_p}+\left\|Q_\f^{[0,1]}\right\|_{\bS_p}\right)\\[.2cm]
&=&C(p)\left(\left\|Q_\f^{[0,1]}\right\|_{\bS_p}+\left\|Q_\psi^{[0,1]}
\right\|_{\bS_p}\right).\quad \bl
\eey

\begin{cor}
\label{Qpsi} 
Under the hypotheses of Theorem \ref{QfiQpsi}
$$
\left\|Q_\psi^{[0,2]}\right\|_{\bS_p}\le C(p)\left\|Q_\f^{[0,2]}
\right\|_{\bS_p}.\quad\bl
$$
\end{cor}

\begin{thm}
\label{Qpsia}
Let $p>0$ and $a\in\R$.
If $\f$ is a function satisfying the hypotheses of Theorem \ref{QfiQpsi} and
$\psi(x)\df\f(x-a)$, then
$$
\left\|Q_\psi^{[0,2]}\right\|_{\bS_p}
\le C(p)\left\|Q_\f^{[0,2]}\right\|_{\bS_p}.
$$
\end{thm}

\Pf Clearly, it is sufficient to consider the case $a\in(0,1)$. 
Then $Q_\psi^{[0,1]}$ is by a translation unitarily equivalent 
to $Q_\f^{[1-a,2-a]}$, and thus
$$
\left\|Q_\psi^{[0,1]}\right\|_{\bS_p}
=\left\|Q_\f^{[1-a,2-a]}\right\|_{\bS_p}
\le\left\|Q_\f^{[0,2]}\right\|_{\bS_p}.
$$
The result follows by \rf{02c} and \rf{02b}. $\bl$

Let $\bsf$ be a function on the unit circle $\T$. Put
$\T_+\df\{\zeta\in\T:\im\zeta\ge0\}$,  
$\T_-\df\{\zeta\in\T:\im\zeta<0\}$ and

\bay
\label{for}
k_{\bsf}(\z,\t)\df\left\{\begin{array}{ll}{\bsf}(\z^2),&\z\ov\t\in\T_+,\\
{\bsf}(\t^2),&\z\ov\t\in\T_-,
\end{array}\right.\quad(\z,\t)\in\T^2.
\ey
It is easy to see that the functions $\bsf$ and $k_{\bsf}$ are
equimeasurable. 
In particular, $\|{\bsf}\|_{L^2(\T)}=\|k_{\bsf}\|_{L^2(\T^2)}$.
Note also that if $\bsf$ is continuous on $\T$, then $k_{\bsf}$ is
continuous on $\T^2$. 
Let ${\bsf}\in L^2(\T)$. Denote by $K_{\bsf}$ the integral operator
on $L^2(\T)$ with kernel function 
$k_{\bsf}$.

\begin{thm} 
\label{tor}
Let $p>0$. Suppose that ${\bsf}\in L^2(\T)$ and
$\f(t)\df{\bsf}(e^{2\pi{\rm i}t})$, $t\in\R$. Then 
$$
C_1(p)\left\|Q_\f^{[0,2]}\right\|_{\bS_p}\le\|K_{\bsf}\|_{\bS_p}\le
C_2(p)\left\|Q_\f^{[0,2]}\right\|_{\bS_p},
$$
where $C_1(p)$ and $C_2(p)$ may depend only on $p$.
\end{thm}

\Pf Consider the integral operator $K$ on $L^2[0,1]$ with kernel function
\linebreak$k\in L^2([0,1]^2)$ defined by
$$
k(x,y)=\left\{\begin{array}{ll}\f(2\max\{x,y\}),&|x-y|\le1/2,\\
\f(2\min\{x,y\}),&|x-y|>1/2.
\end{array}\right.
$$
It is easy to see that $K$ is unitarily equivalent to $K_{\bsf}$.
For $\a,\b=0,1$ we consider the integral operator $K^{(\a,\b)}$ with
kernel function 
$$
(x,y)\mapsto k(x,y)\chi_{[\a/2,(1+\a)/2]}(x)\chi_{[\b/2,(1+\beta)/2]}(y).
$$
Using the substitution $(x,y)\mapsto(2x,2y)$, we find that
$2\|K^{(0,0)}\|_{\bS_p}=\left\|Q_\f^{[0,1]}\right\|_{\bS_p}$. In a
similar way we can obtain 
$2\|K^{(1,1)}\|_{\bS_p}=\left\|Q_\f^{[0,1]}\right\|_{\bS_p}$. Let
$\psi(t)\df\f(1-t)$. 
It is also easy to see that 
$2\|K^{(0,1)}\|_{\bS_p}=2\|K^{(1,0)}\|_{\bS_p}
=\left\|Q_\psi^{[0,1]}\right\|_{\bS_p}$. 
Hence,
\bey
\frac14\left(\left\|Q_\f^{[0,1]}\right\|_{\bS_p}+\left\|Q_\psi^{[0,1]}
\right\|_{\bS_p}\right)&=&
\frac14\sum_{\alpha=0}^1\sum_{\beta=0}^1\left\|K^{(\a,\b)}\right\|_{\bS_p}
\\[.2cm]
&\le&\|K_{\bsf}\|_{\bS_p}
\le C(p)\sum_{\a=0}^1\sum_{\b=0}^1\left\|K^{(\a,\b)}\right\|_{\bS_p}\\[.2cm]
&=&C(p)\left(\left\|Q_\f^{[0,1]}\right\|_{\bS_p}+\left\|Q_\psi^{[0,1]}
\right\|_{\bS_p}\right).
\eey
It remains to apply Theorem \ref{QfiQpsi}. $\bl$

We denote by $\hat f(n)$ denote the $n$th Fourier coefficient of a function
$f$ in $L^1(\T)$. For convenience we put 
$$
\hat f(n+1/2)\df0,\quad n\in\Z.
$$ 

For a function $k$ in $L^1(\T^2)$ we denote by $\{\hat
k(m,n)\}_{(m,n)\in\Z^2}$ the sequence of 
its Fourier coefficients.

Let $\bsf$ be a function on $\T$. Put 
$$
\label{for+}
k_{\bsf}^+(\z,\t)\df
\left\{\begin{array}{ll}{\bsf}(\z^2),&\z\ov\t\in\T_+,\\
0,&\z\ov\t\in\T_-,
\end{array}\right.
$$
and
$$
\label{for-}
k_{\bsf}^-(\z,\t)\df\left\{\begin{array}{ll}0,&\z\ov\t\in\T_+,\\
{\bsf}(\t^2),&\z\ov\t\in\T_-.
\end{array}\right.
$$
Clearly, $k_{\bsf}^+(\z,\t)=k_{\bsf}^-(\t,\z)$ and
$k_{\bsf}(\z,\t)=k_{\bsf}^+(\z,\t)+k_{\bsf}^-(\z,\t)$, where
$k_{\bsf}$ is defined by \rf{for}. 

\begin{lem}
\label{fc+}
Let ${\bsf}\in L^1(\T)$. Then for $(m,n)\in\Z^2$
$$
\hat k_{\bsf}^+(m,n)
=\left\{\begin{array}{ll}\frac12\hat{\bsf}(\frac{m}2),&n=0,\\[.2cm]
\frac{\rmi}{\pi n}\hat{\bsf}(\frac{m+n}2),&m,\,n\text{ are odd},\\[.2cm]
0,&\mbox{otherwise}.
\end{array}\right.
$$
\end{lem}

\Pf Let us first observe that for any $\z\in\T$ and $n\in\Z$ we have
$$
\int\limits_{\{\t\in\T:\z\ov\t\in\T_+\}}\t^{-n}\,d\m(\t)=
\left\{\begin{array}{ll}\frac12,&n=0,\\[.2cm]
\frac{\rmi}{\pi n}\z^{-n},&n\text{ is odd},\\[.2cm]
0,&\mbox{otherwise}.
\end{array}\right.
$$
It follows that
\bey
\int\limits_\T k_{\bsf}^+(\z,\t)\z^{-m}\t^{-n}\,d\m(\zeta)&=&
\int\limits_{\{\t\in\T:\z\ov\t\in\T_+\}}{\bsf}(\z^2)\z^{-m}\t^{-n}\,
d\m(\t)\\[.2cm]
&=&\left\{\begin{array}{ll}\frac12{\bsf}(\z^2)\z^{-m},&n=0,\\[.2cm]
\frac{\rmi}{\pi n}{\bsf}(\z^2)\z^{-m-n},&n\text{ is odd},\\[.2cm]
0,&\mbox{otherwise}.
\end{array}\right.
\eey
It remains to integrate the last identity in $\z$. $\bl$

\begin{cor}
\label{fc-}
Let ${\bsf}\in L^1(\T)$. Then
$$
\hat{k}_{\bsf}^-(m,n)
=\left\{\begin{array}{ll}\frac12\hat{\bsf}(\frac{n}2),&m=0,\\[.2cm]
\frac{\rmi}{\pi m}\hat{\bsf}(\frac{m+n}2),&\,m,\,n\text{ are odd},
\\[.2cm]
0,&\mbox{otherwise}.
\end{array}\right.
$$
\end{cor}

\Pf It suffices to observe that $k_{\bsf}^-(\z,\t)=k_{\bsf}^+(\t,\z)$. $\bl$

\begin{cor}
\label{fc}
Let ${\bsf}\in L^1(\T)$. Then
$$
\hat k_{\bsf}(m,n)=
\left\{\begin{array}{ll}\frac12\hat{\bsf}\left(\frac{m}2\right),& n=0,\\[.2cm]
\frac12\hat{\bsf}(\frac{n}2),&m=0,\\[.2cm]
\hat{\bsf}(0),&m=n=0,\\[.2cm]
\frac{\rmi}{\pi}\left(\frac1m+\frac1n\right)\hat{\bsf}
\left(\frac{m+n}2\right),&m,\,n\text{ are odd},\\[.2cm]
0,&\mbox{otherwise}.
\end{array}\right.
$$
\end{cor}

\Pf It suffices to observe that
$k_{\bsf}(\z,\t)=k_{\bsf}^+(\z,\t)+k_{\bsf}^+(\t,\z)$. $\bl$ 

\begin{thm}
\label{mat}
Let $p>0$. Suppose that $\f(t)=\sum\limits_{k\in\Z}a_ke^{2\pi{\rm i}kt}$ and
$\sum\limits_{k\in\Z}|a_k|^2<\be$. Then $Q_\f^{[0,2]}\in \bS_p$ if and
only if the matrix  
\bay
\label{matritsa}
\left\{a_{m+n+1}\left(\frac1{m+\frac12}+\frac1{n+\frac12}
\right)\right\}_{m,n\in\Z}
\ey
belongs to $\bS_p$.
\end{thm}

Here we identify operators on $\ell^2(\Z)$ with their matrices with
respect to the standard orthonormal basis 
of $\ell^2(\Z)$.

\Pf Consider the function $\bsf$ on $\T$ defined by
${\bsf}(z)=\sum\limits_{n\in\Z}a_nz^n$. 
By Theorem \ref{tor}, $Q_\f^{[0,2]}\in\bS_p$ if and only if
$K_{\bsf}\in\bS_p$  . 
It easy to see that the operator $K_{\bsf}$ belongs to $\bS_p$ if and
only if 
the matrix $\{\hat k_{\bsf}(m,n)\}_{m,n\in\Z}$ belongs to $\bS_p$.
Corollary 
\ref{fc} implies that $\hat k_{\bsf}(m,n)\not=0$ only if $mn$ is odd
or $mn=0$. Hence, 
it is easy to check that $\{\hat k_{\bsf}(m,n)\}_{m,n\in\Z}\in\bS_p$ if
and only if 
$\{\hat k_{\bsf}(2m+1,2n+1)\}_{m,n\in\Z}\in\bS_p$. It remains to note
that  
$$
\hat k_{\bsf}(2m+1,2n+1)
=\frac{\rmi}\pi\left(\frac1{2m+1}+\frac1{2n+1}\right)a_{m+n+1}
$$  
by Corollary \ref{fc}. $\bl$

Clearly, the same reasoning shows that $Q_\f^{[0,2]}$ is bounded on
$L^2[0,1]$ if and only if the matrix \rf{matritsa} is bounded. The
following result shows that  
the boundedness of \rf{matritsa} is equivalent to its membership of
$\bS_p$, $p>1$. 

\begin{thm}
\label{bouS2}
Let $\{a_k\}_{k\in\Z}$ be a two-sided sequence of complex numbers and let 
$$
A\df\left\{a_{m+n+1}
\left(\frac1{m+\frac12}+\frac1{n+\frac12}\right)\right\}_{m,n\in\Z}.
$$
Suppose that $p>1$.
The following are equivalent:
\begin{enumerate}
\item
$A$ is a bounded operator on $\ell^2(\Z)$;
\item
$A\in\bS_p$;
\item
$\{a_k\}_{k\in\Z}\in\ell^2(\Z)$.
\end{enumerate}
\end{thm}

\Pf Suppose that $A$ is bounded. Then the sequence 
$$
\left\{a_{n+1}\left(2+\frac{1}{n+\frac12}\right)\right\}_{n\in\Z}
$$
belongs to $\ell^2(\Z)$ which implies (iii). Clearly, (iii) is equivalent
to the fact that 
\linebreak$\f\in L^2[0,1]$. By Theorem \ref{T:p>1},
$Q_\f^{[0,2]}\in\bS_p$, and so by Theorem \ref{mat}, 
$A\in\bS_p$. The implication (ii)$\Rightarrow$(i) is trivial. $\bl$

Related results, that matrices of a roughly similar sort are bounded
if and only if they are in certain $\bS_p$ can be found in \cite{W1}. 

\begin{remark}
Note that the following identities hold:
$$
\frac{\rmi}\pi\sum\limits_{m,n\in\Z} a_{m+n+1}
\left(\frac1{m+\frac12}+\frac1{n+\frac12}\right)\z^m\t^n=
\frac{{\bsf}(\z)-{\bsf}(\t)}{\sqrt{\z\t}},
$$

$$
\frac{\rmi}\pi\sum\limits_{m,n\in\Z} \frac{a_{m+n+1}}{n+\frac12}\z^m\t^n=
\frac{{\bsf}(\z)}{\sqrt{\z\t}},
$$

$$
\frac{\rmi}\pi\sum\limits_{m,n\in\Z} \frac{a_{m+n+1}}{m+\frac12}\z^m\t^n=
-\frac{{\bsf}(\t)}{\sqrt{\z\t}},
$$
where $\sqrt{\z\t}$ is chosen so that $\overline\t\sqrt{\z\t}\in \T_+$
(the series converge in $L^2(\T^2)$).

Indeed, it suffices to note that by Corollary \ref{fc}, 
$$
\frac{\rmi}\pi\sum\limits_{m,n\in\Z} a_{m+n+1}
\left(\frac1{2m+1}+\frac1{2n+1}\right)\z^{2m+1}\t^{2n+1}=
\frac12\big(k_{\bsf}(\z,\t)-k_{\bsf}(\z,-\t)\big),
$$
by Lemma \ref{fc+},
$$
\frac{\rmi}\pi\sum\limits_{m,n\in\Z}
 \frac{a_{m+n+1}}{2n+1}\z^{2m+1}\t^{2n+1}= 
 \frac12\big(k_{\bsf}^+(\z,\t)-k_{\bsf}^+(\z,-\t)\big),
$$
and by Corollary \ref{fc-},
$$
\frac{\rmi}\pi\sum\limits_{m,n\in\Z}
 \frac{a_{m+n+1}}{2m+1}\z^{2m+1}\t^{2n+1}= 
 \frac12\big(k_{\bsf}^-(\z,\t)-k_{\bsf}^-(\z,-\t)\big).
$$
\end{remark}

\begin{remark}
Note that if $p>1/2$, then 
$$
A\df\left\{a_{m+n+1} \left(\frac1{m+\frac12}+\frac1{n+\frac12}
\right)\right\}_{m,n\in\Z}\in\bS_p
$$
if and only if
\bay
\label{matrB}
B\df\left\{a_{m+n}\left(\frac1{m+\frac12}+\frac1{n+\frac12}
\right)\right\}_{m,n\in\Z}\in\bS_p.
\ey
Indeed, put $\psi(t)\df e^{2\pi{\rm i}t}\f(t)$. By Theorem \ref{mat},
it suffices to prove that  $Q_\f^{[0,2]}\in\bS_p$ implies 
 $Q_\psi^{[0,2]}\in\bS_p$. This follows from Theorem \ref{mlp}.

Note however that for $p\le1/2$ this is not true. Indeed, if
$a_0=1$ and $a_n=0$ for $n\not=0$ (in other words, $\f(t)=1$,
$t\in\R$), then it is easy to see that 
$A$ is the zero matrix, and so it belongs to $\bS_p$ for any $p>0$.
On the other hand, the matrix $B$
has nonzero entries $-(n^2-1/4)^{-1}$ for $m=-n$, and so it belongs to
$\bS_p$ only for $p>1/2$.  
The situation is similar in the case where the restriction
$\f\big|{[0,1)}$ is the 
characteristic function of an interval in which case
$A\in\bS_p$ for any $p>0$ but $B\in\bS_p$ only for $p>1/2$, see
\refT{moc} and \refC{C:moc}. 
\end{remark}

Suppose now that $p>1/2$ and consider the following submatrices of the
matrix $B$ defined by 
\rf{matrB}:
\begin{align*}
B_1&=\left\{a_{m+n}\left(\frac1{m+\frac12}
  +\frac1{n+\frac12}\right)\right\}_{m,n\ge0},\\[.3cm]
B_2&=\left\{a_{m+n}\left(\frac1{m+\frac12}
  +\frac1{n+\frac12}\right)\right\}_{m\ge0,\,n<0},\\[.3cm]
B_3&=\left\{a_{m+n}\left(\frac1{m+\frac12}
  +\frac1{n+\frac12}\right)\right\}_{m<0,\,n\ge0},\\
\intertext{and}
B_4&=\left\{a_{m+n}\left(\frac1{m+\frac12}
  +\frac1{n+\frac12}\right)\right\}_{m,n<0}.
\end{align*}

Clearly, $B\in\bS_p$ if and only if all matrices $B_j$, $1\le j\le4$,
belong to $\bS_p$. 

It is direct that 
$B_1\in\bS_p$ if and only if the matrix
\bay
\label{B1}
\left\{a_{j+k}\frac{j+k+1}{(j+1)(k+1)}\right\}_{j,k\ge0}
\ey
belongs to $\bS_p$. Matrices of the form 
\bay
\label{wH}
\{a_{j+k}(1+j)^\a(1+k)^\b\}_{j,k\ge0}
\ey
are called \emph{weighted Hankel matrices}. It was proved in \cite{Pel2}
that if $\a>-1/2$, $\b>-1/2$, and $0<p\le1$, the matrix \rf{wH}
belongs to $\bS_p$ 
if and only if the function $\sum\limits_{n\ge0}a_nz^n$ belongs to the
Besov class 
$B_{p\,\,p}^{1/p+\a+\b}$ of functions on the unit circle
$\T$. More recent results on Schatten class properties of weighted
Hankel matrices are in \cite{RW} and \cite{W2}. 
However, in the case of interest, the 
weighted Hankel matrix \rf{B1} for $\a=\b=-1$, no characterization of
such matrices 
of class $\bS_p$ is known. In the next section we obtain some
necessary conditions for the 
matrix \rf{B1} to belong to $\bS_1$.

It is also easy to see that $B_4\in\bS_p$ if and only if the weighted
Hankel matrix 
\bay
\label{B4}
\left\{a_{-(j+k+2)}\frac{j+k+1}{(j+1)(k+1)}\right\}_{j,k\ge0}
\ey
belongs to $\bS_p$.

It can also be easily shown that $B_2\in\bS_p$ if and only if
$B_3\in\bS_p$ and this is equivalent 
to the fact that the \emph{weighted Toeplitz matrix} 
\bay
\label{B23}
\left\{a_{j-k-1}\frac{j-k}{(1+j)(1+k)}\right\}_{j,k\ge0}
\ey
belongs to $\bS_p$.

Summarizing the above, we can state the following result.

\begin{thm}
\label{wHT}
Let $1/2<p\le1$ and let $\f$ be a function in $L^2[0,1]$ of the form
$\f(t)=\sum\limits_{n\in\Z}a_ne^{2\pi{\rm i}t}$, $t\in[0,1]$.
Then $Q_\f^{[0,1]}\in\bS_p$ if and only if the matrices \eqref{B1},
\eqref{B4}, and 
\eqref{B23} belong to $\bS_p$. $\bl$
\end{thm}

In the next section we
use the results above 
to obtain 
necessary conditions for the nuclearity of operators $Q_\f$.

Let us consider now the family of functions $\{F_\l\}_{\l\in\C}$ in
$L^2([0,1)^2)$ defined by 
$$
F_\l(t,s)\df\f\left(\max\{s,t\}\right)e^{-2\pi\l{\rm i}|s-t|}+
\f\left(\min\{s,t\}\right)e^{2\pi\l{\rm i}|s-t|-2\pi\l{\rm i}}.
$$
We identify $[0,1)^2$ with $\T^2$ via the map
$(s,t)\mapsto(e^{2\pi{\rm i}s},e^{2\pi{\rm i}t})$ and we 
can consider the Fourier coefficients of functions on $[0,1)^2$.

\begin{thm} 
\label{fk}
Suppose that $\l\not\in\Z$. Then 
\bay
\label{cF}
\widehat F_\l(m,n)=\frac{1-e^{-2\pi{\rm i}\l}}{2\pi{\rm i}}\left(\frac1{\l-m}+
\frac1{\l-n}\right)\hat\f(m+n).
\ey
\end{thm}

\Pf We have 
$$
\widehat F_\l(m,n)\df
\iint\limits_{[0,1)\times[0,1)} F_\l(s,t)
e^{-2\pi{\rm i}ms-2\pi{\rm i}nt}\,dsdt=
\iint\limits_{t\ge s} +
\iint\limits_{t\le s}.
$$ 
Let us compute the first integral:
{\allowdisplaybreaks
\begin{align*}
\iint\limits_{t\ge s}&=
\int\limits_0^1\left(\int\limits_0^t\f(t)
 e^{2\pi{\rm i}\l(s-t)}e^{-2\pi{\rm i}mt-2\pi
{\rm i}ns}\,ds\right)dt\\*[.2cm] &\qquad+
\int\limits_0^1\left(\int\limits_s^1\f(s)
 e^{2\pi{\rm i}\l(t-s)-2\pi{\rm i}\l}e^{-2\pi{\rm i}mt-2\pi
{\rm i}ns}\,dt\right)ds\\[.2cm] &=
\int\limits_0^1\f(t)e^{-2\pi{\rm i}\l t-2\pi{\rm i}mt}
 \left(\int\limits_0^t e^{2\pi
{\rm i}s(\l-n)}\,ds\right)dt\\*[.2cm] &\qquad+
\int\limits_0^1\f(s)e^{-2\pi{\rm i}\l s-2\pi{\rm i}\l-2\pi{\rm i}ns}
\left(\int\limits_s^1e^{2\pi{\rm i}t(\l-m)}\,dt\right)ds\\[.2cm] &=
\int\limits_0^1\f(t)e^{-2\pi{\rm i}\l t-2\pi{\rm i}mt}
 \frac{e^{2\pi{\rm i}t(\l-n)}-1}{2\pi
{\rm i}(\l-n)}dt\\*[.2cm] &\qquad+
\int\limits_0^1\f(s)e^{-2\pi{\rm i}\l s-2\pi{\rm i}\l-2\pi{\rm i}ns}
 \frac{e^{2\pi{\rm i}\l}-
e^{2\pi{\rm i}s(\l-m)}}{2\pi{\rm i}(\l-m)}ds\\[.2cm] &=
\frac{\hat\f(m+n)}{2\pi{\rm i}(\l-n)}-\frac1{2\pi{\rm i}(\l-n)}
\int\limits_0^1\f(t)e^{-2\pi i\lambda t-2\pi imt}\,dt\\*[.2cm] &\qquad+
\frac1{2\pi{\rm i}(\l-m)}\int\limits_0^1\f(s)
 e^{-2\pi{\rm i}\l s-2\pi{\rm i}ns}\,ds-
\frac{e^{-2\pi{\rm i}\l}\hat\f(m+n)}{2\pi{\rm i}(\l-m)}.
\end{align*}
}

Similarly, 
\begin{align*}
\iint\limits_{t\le s}\,&=
\frac{\hat\f(m+n)}{2\pi{\rm i}(\l-m)}-\frac1{2\pi{\rm i}(\l-m)}
\int\limits_0^1\f(s)e^{-2\pi{\rm i}\l s-2\pi{\rm i}ns}\,ds\\
&\qquad+
\frac1{2\pi{\rm i}(\l-n)}\int\limits_0^1\f(t)e^{-2\pi{\rm i}\l
t-2\pi{\rm i}mt}\,dt- 
\frac{e^{-2\pi{\rm i}\l}\hat\f(m+n)}{2\pi{\rm i}(\l-n)}
\end{align*}
which implies \rf{cF}. $\bl$

\begin{thm}
\label{Fla} 
Let $\l\in\C$. The integral operator with kernel function $F_\l$ is
bounded on $L^2[0,1]$ if and only if 
$\f\in L^2([0,1))$.
\end{thm}

\Pf Clearly, the integral operator with kernel function $F_\l$ belongs
to $\bS_2$ if $\f\in L^2([0,1))$. 
Suppose now that the integral operator with kernel function $F_\l$ is
bounded. If $\l\in\Z$, then 
\bey
F_\l(s,t)&=&\left(\f(\max\{s,t\})e^{-4\pi{\rm i}\l\max\{s,t\}}+
\f(\min\{s,t\})e^{-4\pi{\rm i}\l\min\{s,t\}}\right)e^{2\pi{\rm i}\l(s+t)}
\\[.2cm]
&=&\left(\f(s)e^{-4\pi{\rm i}\l s}+
\f(t)e^{-4\pi{\rm i}\l t}\right)e^{2\pi{\rm i}\l(s+t)}.
\eey
Consequently, the boundedness of this operator implies that $\f\in L^2[0,1]$.
If $\l\not\in\Z$, then the boundedness of the integral operator with
kernel function $F_\l$ implies that 
$$
\sum\limits_{m\in\Z}\left|\widehat F_\l(m,0)\right|^2<+\infty
$$
and by Theorem \ref{fk} we obtain
$$
\sum\limits_{m\in\Z}|\hat\f(m)|^2\left|\frac{2\l-m}{\l-m}\right|^2<+\infty,
$$
whence $\f\in L^2[0,1]$. $\bl$

\begin{lem}
\label{multi}
Let $w,a\in\C$ and $w\not=1$. Let $p>\frac12$.
Then the function
$$
(s,t)\mapsto(w-e^{a|s-t|})^{-1}\chi_\Delta(s)\chi_\Delta(t),\quad s,\,t\in\R,
$$ 
is a Schur multiplier of $\bS_p$ if $\Delta$ is an interval of
sufficiently small length. 
\end{lem}
\ 
\Pf Clearly, it suffices to consider the case $p<1$. Note that
$w-e^{a|s-t|}$ is a Schur multiplier 
of $\bS_p(L^2(\Delta))$ by Theorem \ref{mlp}, since 
$$
w-e^{a|s-t|}=w-e^{a\max\{s,t\}}e^{-a\min\{s,t\}}.
$$ 
We have to prove that this multiplier is an isomorphism of
$\bS_p(L^2(\Delta))$ if $\Delta$ 
has sufficiently small length. 
For $\omega\in L^\infty(\Delta^2)$ we put
$$
\|\omega\|_{\Mp(\Delta)}\df\sup\|\omega k\|_{\bS_p(L^2(\Delta))},
$$
where the {\it supremum} is taken over
all integral operators with kernel $k\in L^2(\Delta^2)$ such that
$\|k\|_{\bS_p}=1$. Here by 
$\|k\|_{\bS_p}$ we mean the $\bS_p$ norm (quasi-norm if $p<1$) of the integral
operator with kernel function $k$. Obviously, it suffices to prove the
inequality 
$$
\|e^{a\max\{s,t\}}e^{-a\min\{s,t\}}-1\|_{\Mp(\Delta)}<|w-1|
$$ 
provided the length of $\Delta$ is sufficiently small. Theorem
\ref{mpe} implies that, 
for any $x_0\in\Delta$, 
$$
\lim\limits_{|\Delta|\to0}\|e^{a(\max\{s,t\}-x_0)}-1\|_{\Mp(\Delta)}=0
\quad\mbox{and}\quad
\lim\limits_{|\Delta|\to0}\|e^{-a(\min\{s,t\}-x_0)}-1\|_{\Mp(\Delta)}=0.
$$
Hence, the desired inequality is obvious. $\bl$

\begin{thm}
\label{tab}
 Suppose that $\l\not\in\Z$ and $p>1/2$. Then $Q_\f^{[0,2]}\in\bS_p$
if and only if 
the integral operator with kernel function $F_\l$ belongs to $\bS_p$.
\end{thm}

\Pf Suppose that  $Q_\f^{[0,2]}\in\bS_p$. Then the integral operators
with kernel  
functions $\f(\max\{s,t\})$ and
$\f(\min\{s,t\})$ belong to $\bS_p(L^2[0,1])$ (see Theorem
\ref{QfiQpsi}). Note that 
$e^{2\pi{\rm i}\l|s-t|}=e^{2\pi{\rm i}\l(2\max\{s,t\}-s-t)}$. It
follows now from Theorem \ref{mlp} 
that the integral operator with kernel function $F_\l$ belongs to $\bS_p$.

Suppose now that the integral operator with kernel function $F_\l$
belongs to $\bS_p$. We have to prove 
that $Q_\f^{[0,2]}\in\bS_p$. By Theorem \ref{Fla}, $\f\in L^2[0,1]$,
and so it suffices to show that 
$Q_\f^{[0,1]}\in\bS_p$. By Lemma \ref{multi}, we can
choose a positive number $\d$ such that the function
$(e^{2\pi{\rm i}\l}-e^{4\pi{\rm i}\l|s-t|})^{-1}$ belongs to $\Mp(\Delta)$
for any interval $\Delta$ of length less than $\d$. 
We can represent the interval $[0,1)$ in the form
$\bigcup\limits_{j=1}^N\Delta_j$, where the
$\D_j$ are pairwise disjoint intervals with lengths less than $\d$. Clearly,
$$
F_\l(s,t)=\f(\max\{s,t\})e^{-2\pi\l{\rm i}|s-t|}+
(\f(s)+\f(t)-\f(\max\{s,t\})e^{2\pi\l{\rm i}|s-t|-2\pi\l{\rm i}}.
$$
Let $s,\,t\in\D_j$. Then 
\begin{multline*}
\f(\max\{s,t\})
=\frac{F_\l(s,t)-(\f(s)+\f(t))e^{2\pi\l{\rm i}|s-t|-2\pi\l{\rm i}}}
{e^{-2\pi\l{\rm i}|s-t|}-e^{2\pi\l{\rm i}|s-t|-2\pi\l{\rm i}}}=\\[0.2cm]
\frac{\left(F_\l(s,t)-(\f(s)+\f(t))
e^{2\pi\l{\rm i}(\max\{s,t\}-\min\{s,t\}-1}\right)}
{e^{2\pi\l{\rm i}}-e^{4\pi\l{\rm i}|s-t|}}
e^{2\pi\l{\rm i}(1+\max\{s,t\}-\min\{s,t\})}.
\end{multline*}
Theorem \ref{mlp} and Lemma \ref{multi} imply that the integral
operator with kernel function 
$$
(s,t)\mapsto\f(\max\{t,s\})\chi_{\Delta_j}(t)\chi_{\Delta_j}(s),
\quad s,\t\in\R,
$$ 
belongs to $\bS_p$.
To complete the proof, it remains to observe that the kernel function
$$
(s,t)\mapsto\f(\max\{s,t\})-\sum\limits_{j=1}^N\f(\max\{s,t\})
\chi_{\D_j}(s)\chi_{\D_j}(t)
$$
determines a finite rank operator. $\bl$

Theorem \ref{tab} implies that if $p>1/2$, $\l_1,\,\l_2\notin\Z$, and
$F_{\l_1}\in\bS_p$, then 
$F_{\l_2}\in\bS_p$. This can also be easily deduced from the following
elementary fact: 
if $x\in\ell^2(\Z)$ and $y\in\ell^p(\Z)$ with $p\le2$, then
$\{x_{m+n}\,y_n\}_{m,n\in\Z}\in\bS_p$.

\

\setcounter{equation}{0}
\section{\bf Necessary conditions for $\bs{Q_\f\in\bS_1}$}
\label{neces}

\

In this section we obtain various necessary conditions for $Q_\f\in\bS_1$.

\begin{thm}
\label{wHsuf}
Let $\{a_n\}_{n\ge0}$ be a sequence in $\ell^2$. If the matrix
$$
\G=\left\{a_{j+k}\left(\frac1{j+\frac12}+\frac1{k+\frac12}
\right)\right\}_{j,k\ge0} 
$$
belongs to $\bS_1$, then the function
$\sum\limits_{n\ge0}\log(2+n)a_nz^n$ belongs to the Hardy class 
$H^1$.
\end{thm}

We need the following well-known lemma (see e.g., \cite{Pel1}).

\begin{lem}
\label{es1}
Suppose that the matrix $\{a_{jk}\}_{j,k\ge0}$ belongs to
$\bS_1$. Then the function  
$\sum\limits_{n\ge0}\left(\sum\limits_{j=0}^n a_{j\,\,n-j}\right)z^n$
belongs to the Hardy class 
$H^1$.
\end{lem}

\Pf It is sufficient to prove this when the matrix has rank one in which case
the result is an immediate consequence of the fact that 
$H^2\cdot H^2\subset H^1$. 
$\bl$




\begin{lem}
\label{ell}
Let $m\in\Z_+$ and let
\bay
\label{betam}
\b_n\df\sum_{j=0}^n\frac1{j+\frac12}.
\ey
Then there exists $d\in\R$ such that
$$
|\b_n-\log(2+n)-d|\le\const\frac1{1+n}.
$$
\end{lem}

\Pf We use
the following well known fact 
(see, for example, \cite[Ch.~I, (8.9)]{Z}) 
\bay
\label{eul}
\left|\sum_{j=1}^n\frac1j-\log n-\gamma\right|\le\const\cdot n^{-1},
\ey
where $\gamma$ is the Euler constant. We have
$$
\sum_{j=0}^n\frac1{j+\frac12}=2\sum_{j=0}^n\frac1{2j+1}=
2\sum_{j=1}^{2n+1}\frac1j-\sum_{j=1}^n\frac1j,
$$
and so by \rf{eul},
$$
\left|\sum_{j=0}^n\frac1{j+\frac12}-2\log(2n+1)+\log n -\gamma\right|
\le\const\cdot n^{-1}
$$
which implies the result. $\bl$

{\bf Proof of Theorem \ref{wHsuf}.} 
By Lemma \ref{es1}, we have 
\bey
\sum_{m\ge0}a_m\left(\sum_{j=0}^m
\left(\frac1{j+\frac12}+\frac1{(m-j+\frac12)}\right)\right) z^m=
2\sum_{m\ge0}a_m\left(\sum\limits_{j=0}^m\frac1{j+\frac12}\right)z^m\in H^1.
\eey
Since $\{a_n\}_{n\ge0}\in\ell^2$, it is not hard to check 
that Lemma
\ref{ell} implies that  
\linebreak$\sum\limits_{n\ge0}\log(2+n)a_nz^n\in H^1$. $\bl$

\begin{thm}
\label{L2p1}
Let $\f\in L^2[0,b]$ and $\f(x)=\sum\limits_{n\in\Z}a_ne^{2\pi{\rm i}nx/b}$.
If $Q_\f^{[0,b]}\in\bS_1$, then the functions
$\sum\limits_{n\ge0}a_n\log(2+n)z^n$ and
$\sum\limits_{n\ge0}a_{-n}\log(2+n)z^n$ in the unit disc $\dd$ 
belong to the Hardy class $H^1$.
\end{thm}

\Pf Without loss of generality we may assume that $b=1$. By Theorem
\ref{mat}, the matrices 
$$
\G=\left\{a_{j+k+1}\left(\frac1{j+\frac12}+\frac1{k+\frac12}
\right)\right\}_{j,k\ge0}
$$
and 
$$
\G=\left\{a_{-(j+k+1)}\left(\frac1{j+\frac12}+\frac1{k+\frac12}
\right)\right\}_{j,k\ge0}
$$
belong to $\bS_1$. The result follows now from Theorem \ref{wHsuf}. $\bl$

\begin{thm}
\label{nec1}
Let $I$ be a compact interval in $(0,\be)$ and let $\f$ be a function in
$\Lloc^1(\R_+)$ such that 
$Q_\f\in\bS_1$. If 
$$
a_n=\int_I\f(x) e^{-2\pi{\rm i}nx/|I|}dx,\quad n\in\Z,
$$
then the functions $\sum\limits_{n\ge0}a_n\log(2+n)z^n$ and
$\sum\limits_{n\ge0}a_{-n}\log(2+n)z^n$
belong to the Hardy class $H^1$.
\end{thm}

\Pf Since $I$ is separated away from 0, it follows that $\f\big|I\in
L^2(I)$. We can now apply a 
translation and reduce the result to Theorem \ref{L2p1}. $\bl$

\begin{cor}
\label{coef}
Under the hypotheses of either Theorem \ref{L2p1} or \ref{nec1}  the
following holds: 
\begin{enumerate}
\item $|a_n|\le\const(\log(2+|n|))^{-1}$, $n\in\Z$;\vspace*{.3cm}
\item suppose that $\{n_k\}_{k\ge0}$ is an Hadamard lacunary sequence
of positive integers, i.e., 
$$
\inf_{k\ge0}\frac{n_{k+1}}{n_k}>1,
$$
then
$$
\sum_{k\ge0}|a_{n_k}|^2(\log(1+n_k^2))^2<\be\quad\mbox{and}\quad
\sum_{k\ge0}|a_{-n_k}|^2(\log(1+n_k^2))^2<\be.
$$
\end{enumerate}
\end{cor}

\Pf (i) follows immediately from Theorem \ref{nec1} and the obvious
fact that the Fourier  
coefficients of an $H^1$ function are bounded. Finally, (ii) is an
immediate 
consequence of Theorem \ref{nec1} and Paley's inequality 
(see \cite[v.~2, Ch.~XII, (7.8)]{Z}). $\bl$ 

Note that if $I$ is a compact interval in $(0,\be)$, the restrictions
of function in $X_1$ to 
$I$ fill the space $L^2(I)$, and so the sequence of Fourier
coefficients $\{a_n\}_{n\in\Z}$ can be 
an arbitrary sequence in $\ell^2$. Thus Corollary \ref{coef} also
shows that the condition $\f\in X_1$ 
is not sufficient for $Q_\f\in\bS_1$.

Now we are going to use Theorem \ref{boup} to obtain another necessary
condition for 
$Q_\f\in\bS_1$.

We denote by ${\frak H}^1$ the Stein--Weiss space of functions $f$ in
$L^1(\R)$ such that 
$\F^{-1}(\chi_{\R_+}\F f)\in L^1(\R)$, where $\F$ is Fourier transformation.

\begin{thm}
\label{log}
Let $h\in \Lloc^2(\R)$. Suppose that the integral operator on
$L^2(\R)$ with kernel function 
$$
(x,y)\mapsto 
h^\spadesuit(x,y)\df h(x+y)\left(\frac1{x+{\rm i}}+\frac1{y-{\rm i}}\right),
\quad x,\,y\in\R,
$$
belongs to $\bS_1$. Then the Fourier transform of the function
$h(x)\log(1+x^2)$ belongs to the Stein--Weiss space ${\frak H}^1$. 
\end{thm}

\Pf Clearly, the integral operator with kernel function
$h^\spadesuit\chi_{[0,+\infty)^2}$ belongs to
$\bS_1$. Put 
$$
g(x)\df\int\limits_{\R}h^\spadesuit(t,x-t)\chi_{[0,+\infty)^2}(t,x-t) dt.
$$
We have
$$
g(x)=\left\{\begin{array}{ll}
h(x)\int\limits_0^x\left(\frac1{t+{\rm i}}+\frac1{x-t-{\rm i}}\right)
dt=h(x)\log(1+x^2),&x>0,\\[.2cm] 
0,&x<0.\end{array}\right.
$$
It follows from Theorem \ref{T:trace} that $\F g\in L^1(\R)$. In the same way 
it can be shown that the Fourier
transform of the function $h(x)\log(1+x^2)\chi_{\R_-}(x)$ belongs to $L^1(\R)$.
This implies the result. $\bl$

\begin{cor}
\label{log+} 
Let $h\in L^2_{\rm loc}(\R)$ and let $a,b\in\C\setminus\R$ such that
$a+b\notin\R$. 
Suppose that 
the integral operator on $L^2(\R)$ with kernel function
$$
h^\spadesuit_{a,b}(x,y)\df h(x+y)\left(\frac1{x+a}+\frac1{y+b}\right)
$$
belongs to $\bS_1$. Then the Fourier transform of the function
$h(x+c)\log(1+x^2)$ belongs to  ${\frak H}^1$ for any $c\in\R$.
\end{cor}

\Pf Clearly, the integral operators on $L^2(\R)$ with kernels
functions \linebreak 
$h(x+y+c)\left(\frac1{x+a+c}+\frac1{y+b}\right)$ belong to
$\bS_1$. Consequently, by 
Theorem \ref{del}, the integral operator on $L^2(\R)$ with kernel function
$h(x+y+c)\left(\frac1{x+\ii}+\frac1{y-\ii}\right)$ belongs to
$\bS_1$. It remains to apply Theorem 
\ref{log}. $\bl$ 

\begin{cor}
\label{log0}
Suppose that $h$, $a$ and $b$ satisfy the hypotheses of Corollary \ref{log+}.
Then $h(x)\log|x|\to0$ as $|x|\to\infty$. $\bl$
\end{cor}

Now Corollary \ref{log0} and Theorem \ref{boup} imply the following theorem

\begin{thm}
\label{mnc}
Let
{\em$\f\in L^2_{\text{loc}}(\R_+)$}. Suppose that
$Q_\f\in\bS_1$. Then
\bay
\label{log|x|}
\log|x|\int\limits_{\R_+}\f(t)t^{x{\rm i}} dt
\to0\quad\mbox{as}\quad|x|\to\infty . \quad\bl
\ey
\end{thm}

Note that it follows from Theorem \ref{T:x1} that $\f\in L^1(\R_+)$, and
so the integral in \rf{log|x|} is well defined. It is easy to see that
if $\f$ is an arbitrary 
$L^2$ function supported on a compact subset of $(0,\be)$, then $\f\in
X_1$. However,  
$\f$ does not have to satisfy \rf{log|x|}, and so Theorem \ref{mnc}
also implies that the 
condition $\f\in X_1$ is not sufficient for $Q_\f\in\bS_1$. 

We conclude this section with necessary conditions on the $L^1$
modulus of continuity of symbols. 
If $f$ is a function on $\T$, then its $L^1$ modulus of continuity
$\oxi f$ is defined by, in analogy with \eqref{contmod}, 
$$
\oxi f(t)\df\sup_{\z\in\T,~|1-\z|<t}
\int\limits_\T |f(\z\t)-f(\t)|\,d\m(\t),\quad t>0.
$$

The following result is possibly known to experts. We were not able to
find a reference, and 
we prove it here. 

\begin{thm}
\label{trlog}
Let $f\in L^1(\T)$ and let 
$$
g(z)=\sum_{n\in\Z}\frac{\hat f(n)}{\log(|n|+2)}z^n.
$$
Then $g\in L^1(\T)$ and 
\begin{equation}
\label{1}
\lim_{t\to0}\oxi g(t)\log \frac1t=0.
\end{equation} 
\end{thm} 

Consider the function $\bsh$ on $\T$ defined by
$$
\bsh(z)\df\sum\limits_{n\in\Z}\left(\log(|n|+2)\right)^{-1}z^n.
$$ 
It is well-known (see, for example, \cite[Ch.~V, (1.5)]{Z}) that the
series converges for 
$z\in\T\setminus\{1\}$, $\bsh\ge0$ and $\bsh\in L^1(\T)$.
We define the function $h$ on $\R$ by 
\bay
\label{h}
h(x)\df\bsh(e^{{\rm i}x})
=(\log 2)^{-1}+2\sum\limits_{n\ge1}\left(\log(n+2)\right)^{-1}\cos nx. 
\ey
Then $h$ is continuously differentiable on $\R\setminus2\pi\Z$, see 
\cite[Ch.~V, Miscellaneous theorems and examples, 7]{Z}.

\newcommand{\si}{\stackrel0{\sim}}

We use the following notation. Let $\f$ and $\psi$ be nonvanishing
functions on an interval  
$(0,\a)$. We write 
$$
\f\stackrel0{\sim}\psi,\quad\mbox{if}\quad\lim_{x\to0}\frac{\f(x)}{\psi(x)}=1.
$$

\begin{lem} 
\label{as}
Let $h$ be the function defined by \eqref{h}. Then
\bay
\label{hx}
h(x)\si\frac\pi{x(\log x)^2}
\ey
and
\bay
\label{h'}
h^\prime(x)\si-\frac\pi{x^2(\log x)^2}.
\ey
\end{lem}

\Pf \rf{hx} is proved in \cite[Ch.~V, (2.17)]{Z}. Let us prove \rf{h'}.
Using Abel's transformation, we obtain 
\begin{align*}
h(x)&=\sum_{n\ge0}
\left((\log(n+2))^{-1}-(\log(n+3))^{-1}\right)\frac{\sin(n+\frac x2)}
{\sin\frac x2}\\[.2cm]
&=\cot\frac
x2\left(\sum_{n\ge0}\left((\log(n+2))^{-1}-(\log(n+3))^{-1}
\right)\sin nx\right)\\[.2cm]
&\qquad+
\sum\limits_{n\ge0}\left((\log(n+2))^{-1}-(\log(n+3))^{-1}\right)\cos nx.
\end{align*}
Consequently, 
\begin{align*}
h^\prime(x)&=-\frac1{2\sin^2\frac x2}\left(\sum_{n\ge0}
\left((\log(n+2))^{-1}-(\log(n+3))^{-1}\right)\sin nx\right)\\[.2cm]
&\qquad+\cot\frac x2\left(\sum\limits_{n\ge0}
\left((\log(n+2))^{-1}-(\log(n+3))^{-1}\right)n\cos nx\right)\\[.2cm]
&\qquad-\sum_{n\ge0}\left((\log(n+2))^{-1}-(\log(n+3))^{-1}\right)n\sin nx
\df\Sigma_1+\Sigma_2+\Sigma_3.
\end{align*}
It remains to observe that 
$$
\Sigma_1\si-\frac\pi{x^2(\log x)^2}
$$
by \cite[Ch.~V, (2.13)]{Z}, while 
$$
\Sigma_2\si-\frac{2\pi}{x^2(\log x)^3},
$$
and 
$$
\Sigma_3\si-\frac1{x(\log x)^2}
$$
by \cite[Ch.~V, (2.18)]{Z}.
$\bl$ 

\begin{cor}
\label{Hlog}
The following inequality holds
$$
\int\limits_{\T}|\bsh(\zeta\t)-\bsh(\t)|d\m(\t)
\le C\left(\log\frac3{|\zeta-1|}\right)^{-1}
$$ 
for any $\zeta\in\T$.
\end{cor}

\Pf It suffices to prove that 
$$
\int\limits_{-\pi}^{\pi}|h(x+t)-h(x)|\,dx\le C\left(\log\frac 1t\right)^{-1}
$$
for sufficiently small positive $t$. We have
\begin{align*}
\int\limits_{-\pi}^{\pi}|h(x+t)-h(x)|\,dx
&=\int\limits_{|x|\le2t}|h(x+t)-h(x)|\,dx\\[.2cm]
&\qquad+\int\limits_{2t\le|x|\le\pi}|h(x+t)-h(x)|\,dx\\[.2cm]
&\le2\int\limits_{|x|\le3t}|h(x)|dx +
\int\limits_{2t\le|x|\le\pi}|h(x+t)-h(x)|\,dx\\[.2cm]
&\df2I_1+I_2.
\end{align*}
Using Lemma \ref{as}, we obtain
$$
I_1\le C\int\limits_0^{3t}\frac1{x(\log x)^2}dx
\le C\left(\log\frac 1t\right)^{-1}
$$
and
$$
I_2\le Ct\int\limits_{2t}^\pi\frac{dx}{x^2(\log\frac x{10})^2}
\le C\left(\log\frac 1t\right)^{-2},
$$
if $t>0$ is sufficiently small. 
$\bl$

{\bf Proof of Theorem \ref{trlog}.} Note that
$g=f*\bsh$. Consequently, $g\in L^1(\T)$.  
It follows easily from Corollary \ref{Hlog} that
$$
\oxi g(t)\le\const
\|f\|_{L^1}\left(\log\frac3{t}\right)^{-1}, \quad 0<t<2.
$$
The result follows now from the obvious fact that \rf{1} holds for
trigonometric polynomials $f$. 
$\bl$

For a function $f\in L^1(\R)$, we defined the $L^1$ modulus of
continuity $\oxi f$ in  \eqref{contmod}: 
$$
\oxi f(t)=\sup_{|s|\le t}\int\limits_\R|f(x+s)-f(x)|dx,\quad t>0.
$$
In fact this definition can be extended to functions $f$ not
necessarily in $L^1(\R)$. It is sufficient 
to assume that 
$$
\int\limits_\R|f(x+s)-f(x)|dx<\be,\quad s\in\R.
$$

In a similar way we can prove the following analog of  Theorem \ref{trlog}.

\begin{thm}
\label{analog}
Let $f\in L^1(\R)$. Then there exists a function $g\in L^1(\R)$ such that
$$
(\F f)(x)=(\F g(x))\log(|x|+2), \quad x\in\R,
$$ 
and 
$$
\lim_{t\to0}\oxi g(t)\log \frac1t=0
$$ 
\end{thm}

\Pf Indeed, let 
$$
{\frak h}(x)\df\int\limits_{\R}\left(\log(2+|t|)\right)^{-1}
e^{-2\pi{\rm i}tx}\,dt=
2\int\limits_0^\infty\left(\log(2+t)\right)^{-1}\cos(2\pi tx)\,dt.
$$
Then $\frak h$ is an even positive continuously differentiable
function on $\R\setminus\{0\}$.  
We can repeat the above reasoning to prove that
$$
{\frak h}(x)\si\frac1{2x(\log x)^2}
$$ 
and
$$
{\frak h}^\prime(x)\si-\frac1{2x^2(\log x)^2}.
$$
Moreover, $|{\frak h}(x)|\le\const\cdot\, x^{-2}$ and 
$|{\frak h}^\prime(x)|\le\const\cdot\, x^{-2}$ 
everywhere. These estimates allow us to obtain the inequality
$$
\int\limits_{\R}|{\frak h}(x+t)-{\frak h}(x)|\,dx
\le\const\left(\log\frac1t\right)^{-1}
$$
for $t\in(0,\frac12)$ and repeat the reasoning in the proof of Theorem
\ref{trlog}. $\bl$ 

Let us introduce some more notation. Set $\C_+\df\{z\in\C:\im z>0\}$
and \linebreak 
$\C_-\df\{z\in\C:\im z<0\}$. Let $f$ be a function in $\Lloc^1(\R)$ such that 
$$
\int\limits_{\R}\frac{|f(t)|}{1+|t|}\,dt<+\infty.
$$ 
Consider the Cauchy transform of $f$ defined by
$$
({\cal C} f)(\z)\df\frac1{2\pi{\rm i}}\int\limits_{\R}\frac{f(t)\,dt}{t-\z},
\quad\im\z\ne0.
$$
It is well known that $({\cal C}f)\big|\C_+$ and 
$({\cal C}f)\big|\C_-$  
are holomorphic functions of bounded characteristic in $\C_+$ and
$\C_-$ respectively, and so 
they have finite angular boundary values almost everywhere
on $\R$. Set 
$$
f_+(x)\df\lim\limits_{y\to0+}({\cal C}f)(x+{\rm i}y)\quad\mbox{and}\quad
f_-(x)\df-\lim\limits_{y\to0-}({\cal C}f)(x+{\rm i}y).
$$
By the Privalov theorem (see \cite{Pr} for the case of a Jordan curve),
$f=f_++f_-$ almost everywhere on $\R$. If 
$f\in L^1(\R)+L^2(\R)$, then 
$$
({\cal C}f)(z)=\int\limits_0^\infty(\F f)(t)e^{2\pi{\rm i}tz}\,dt,
\quad z\in\C_+,
$$ 
and
$$
({\cal C}f)(z)=-\int\limits_{-\infty}^0(\F f)(t)e^{2\pi{\rm i}tz}\,dt,
\quad z\in\C_-.
$$ 
Note that  $f_+$ does not have to be in $L^1(\R)$ for an arbitrary
function $f\in L^1(\R)$. 
In fact, if $f\in L^1(\R)$, then  $f_+\in L^1(\R)$ if and only if 
$f$ belongs to the Stein--Weiss space  ${\frak H}^1$.

\begin{thm}
\label{intlog}
Let $f\in{\frak H}^1$. Suppose that there exists a function $g\in
L^1(\R)$ such that 
$(\F f)(x)=(\F g)(x)\log(1+x^2)$ for all $x\in\R$. Then 
$$
\lim_{t\to0}\oxi {g_+}(t)\log\frac1t=0
$$ 
and
$$
\lim_{t\to0}\oxi {g_-}(t)\log\frac1t=0.
$$ 
\end{thm} 

Note, however, that the assumptions of Theorem \ref{intlog}  do not
imply that \linebreak$g_+\in L^1(\R)$ or $g_-\in L^1(\R)$.

We need some auxiliary facts. Let $\M(\R)$ be the space of finite
Borel measures on $\R$. 

\begin{lem}
\label{FL1}
Let $f\in L^1(\R)$. Suppose that $f^{\prime\prime}\in\M(\R)$ (in the
distributional sense). 
Then $\F f\in L^1(\R)$ and 
$$
\|\F f\|_{L^1(\R)}\le C\sqrt{\|f\|_{L^1(\R)}\|f^{\prime\prime}\|_{\M(\R)}}.
$$
\end{lem}
\Pf The result follows from the obvious inequality: 
$$
|(\F f)(x)|\le\min\left\{\|f\|_{L^1},
\frac{\|f^{\prime\prime}\|_{\M(\R)}}{4\pi^2x^2}\right\},
\quad x\in\R.\quad\bl
$$

\begin{cor}
\label{oh1}
Let $f\in L^1(\R)$. Suppose that $\supp\F f$ is bounded above. Then
$$
\oxi {f+}(t)\le\const\cdot\,t,\quad t>0.
$$ 
\end{cor}

\Pf It suffices to construct a function
$g_s\in L^1(\R)$ such that $\|g_s\|_{L^1(\R)}\le C|s|$ and
$f_+(x+s)-f_+(x)=(f*g_s)(x)$ for all $x\in\R$. 
Suppose that $\supp f\subset(-\infty,M]$, where $M>0$.
We may take a function $g_s$ such that 
\begin{equation*}
(\F g_s)(t)=
\begin{cases}
0, & t\le0,\\
e^{2\pi {\rm i}st}-1, & t\in[0,M], \\
e^{2\pi {\rm i}s(2M-t)}-1, & t\in[M,2M], \\
0, & t\ge 2 M.
\end{cases}
\end{equation*}
Clearly, $f_+(x+s)-f_+(x)=(f*g_s)(x)$ for all $x\in\R$. The inequality
$\|g_s\|_{L^1(\R)}\le C|s|$ 
follows from Lemma \ref{FL1} (with $C$ depending on $M$). 
$\bl$

\begin{lem}
\label{upr}
Set $\rho(t)\df2\log(2+|t|)-\log(1+t^2)$.
Then $\F\rho\in L^1$.
\end{lem}
\Pf It suffices to observe that $\rho$ is even,
$\lim\limits_{t\to\infty}\rho(t)=0$, 
$\r$ has two continuous derivatives on $(0,\be)$, and
$$
\int\limits_0^\be t|\r''(t)|dt<\be.\quad\bl
$$

\begin{lem}
\label{upr1}
Let $\f$ be an even positive function in $C^2(\R)$ 
such that $\f(x)=\log(1+x^2)$
for sufficiently 
large $|x|$. Then $\F(\f^{-1})\in L^1$.
\end{lem}
\Pf See the proof of the previous lemma.
$\bl$

{\bf Proof of Theorem \ref{intlog}.} We prove the first equality (the
proof of the second one is the same). 
Let $\psi$ be a function in $L^1(\R)$ such that $\supp\F\psi$ is
compact and 
$\supp\F\psi=1$ in a neighborhood of $0$. Then
$f=f*\psi+(f-f*\psi)$. The Fourier  
transform of the first summand has a compact support while the support
of the Fourier 
transform of the second summand does not contain $0$. Thus it is sufficient to
consider two cases.

Case 1, $\supp\F f$ is compact. The result follows from Corollary \ref{oh1}.

Case 2, $0\notin\supp\F f$.
Clearly, 
$$
(\F f_+)(x)=(\F g_+)(x)\log(1+x^2),\quad x\in\R.
$$
Lemma \ref{upr1} implies that $g_+\in L^1(\R)$. On the other hand, it
follows from 
Lemma \ref{upr} that
$(\F g_+)(x)\log(2+|x|)$ is the Fourier transform of an $L^1$-function.
It remains to apply Theorem \ref{analog}. $\bl$

\begin{thm}
\label{L2p11}
Let $\f$ be a function in $L^2[0,b]$ such that $Q_\f^{[0,b]}\in\bS_1$ and let
$\f(x)=\sum\limits_{n\in\Z}a_ne^{2\pi{\rm i}nx/b}$. 
If 
\bay
\label{fi+-}
\bsf_+(\z)\df\sum\limits_{n\ge0}a_n \z^n\quad\mbox{and}\quad 
\bsf_-(\z)\df\sum\limits_{n<0}a_n \z^n,
\ey
then
\bay
\label{of+-}
\lim_{t\to0}\oxi{\bsf_+}(t)\log\frac{1}{t}=0\quad\mbox{and}\quad
\lim_{t\to0}\oxi{\bsf_-}(t)\log\frac1t=0.
\ey
\end{thm}

\Pf The result follows immediately from Theorems \ref{L2p1} and \ref{trlog}.
$\bl$

\begin{thm}
\label{S1nec} 
Let $I$ be a compact interval in $(0,\be)$ and let $\f$ be a function in
$L^1_{\rm loc}(\R_+)$ such that
$Q_\f\in\bS_1$. If 
$$
a_n=\int_I\f(x) e^{-2\pi{\rm i}nx/|I|}dx,\quad n\in\Z,
$$
and $\bsf_+$ and $\bsf_-$ are defined by \eqref{fi+-}, then
\eqref{of+-} holds.
\end{thm}

\Pf The result is an immediate consequence of Theorem \ref{L2p11}. $\bl$

Recall that for a function $\f\in L^2_{\rm loc}(\R_+)$ the function
$\f^\hs$ is defined by \rf{ser}. 
Note that if $\qf\in\bS_1$, then by Theorem \ref{T:x1}, $\f\in L^1(\R_+)$
and thus $\f^\hs\in L^1(\R)$. 

\begin{thm}
\label{mnc0}
Let
$\f$ be a function in $L^2_{\rm{loc}}(\R_+)$ such that
$Q_\f\in\bS_1$. Then
$$ 
\lim_{t\to0}\oxi {(\f^\hs)_+}(t)\log\frac{1}{t}=0
\qquad\text{and}\qquad
\lim_{t\to0}\oxi {(\f^\hs)_-}(t)\log\frac{1}{t}=0.
$$
In particular, 
$$ 
\lim_{t\to0}\oxi {\f^\hs}(t)\log\frac{1}{t}=0.
$$
\end{thm}

\Pf The result follows from Theorem \ref{boup}, Corollary \ref{log+},
and Theorem \ref{intlog}. 
$\bl$

This result should be compared to Theorems
\refand{T:Dini2}{T:DinipR}. 
In particular, if $\f$ has compact support in $\R_+$, we see that a
Dini condition on the $L^2$ modulus of continuity is sufficient for
$\qf\in\bS_1$, while the slightly weaker condition 
$\lim_{t\to0}\oxi {\f}(t)\log\frac{1}{t}=0$ on the $L^1$ modulus
of continuity
is necessary.

\begin{thm}
\label{mnc1}
Let
$\f$ be a function in $L^2_{\rm{loc}}(\R_+)$ such that
$Q_\f\in\bS_1$. Then
$$
\lim_{a\to1}\int\limits_{\R_+}|\f(ax)-\f(x)|\,dx\cdot\log\frac{1}{|a-1|}=0.
$$
\end{thm}

\Pf By Theorem \ref{mnc0}, we have
$$
\lim_{t\to0}\int\limits_{\R}\left|\f(e^{2s+2t})e^{2s+2t}-\f(e^{2s})e^{2s}
\right|\,ds\cdot\log\frac1t=0.
$$
Substituting $e^{2s}=x$ and $e^{2t}=a$, we obtain
$$
\lim_{a\to1}\int\limits_{\R_+}|a\f(ax)-\f(x)|\,dx\cdot\log\frac{1}{|a-1|}=0.
$$
It remains to observe that by Theorem \ref{T:x1}, $\f\in L^1(\R_+)$
and obviously, 
$$
\lim_{a\to1}|a-1|\cdot\log\frac{1}{|a-1|}=0.\quad\bl
$$

\

\setcounter{equation}{0}
\section{\bf Dilation of Symbols}
\label{S:dilation}

\

Let $\f$ be a function in $L^2_{\text{loc}}(\R)$ such that
$\f(x+1)=\f(x)$, $x\in\R$.   
For $a>0$ we define the function
$\f_a$ on $[0,1]$ by $\f_a(x)\df\f(ax)$ for $x\in[0,1)$. We are going
to obtain in this section 
upper and lower estimates for $\left\|Q_{\f_a}^{[0,1]}\right\|_{\bS_p}$.

Note that we can extend $\f_a$ to $\R$ as a $1$-periodic function on
$\R$. Using an obvious estimate,
see \eqref{02a}--\eqref{02c}, 
\bay
\label{012}
C_1\left(\left\|Q_{\f_a}^{[0,1]}\right\|_{\bS_p}+\|\f_a\|_{L^2[0,1]}\right)
\le\left\|Q_{\f_a}^{[0,2]}\right\|_{\bS_p}
\le C_2\left(\left\|Q_{\f_a}^{[0,1]}\right\|_{\bS_p}
 +\|\f_a\|_{L^2[0,1]}\right),
\ey
we can reduce the estimation of $\|Q_{\f_a}^{[0,1]}\|_{\bS_p}$ to that
of $\|Q_{\f_a}^{[0,2]}\|_{\bS_p}$. 
We can consider the Fourier coefficients of $\f_a$ defined by
$$
\hat\f_a(n)\df\int_0^1\f_a(t)e^{-2\pi{\rm i}nt}dt,\quad n\in\Z.
$$

\begin{thm}
\label{spn}
Let $\f$ be a $1$-periodic function in $L^2_{\rm{loc}}(\R)$ and let
$a>0$. Suppose that 
$\f$ has bounded variation on $[0,1]$. Then 
$$
\|Q_{\f_a}^{[0,2]}\|_{\bS_1}\le C(\f)\log(2+a)
$$ 
and 
$$
\|Q_{\f_a}^{[0,2]}\|_{\bS_p}\le C(\f)(1+a)^{1/p-1},\quad 1/2<p<1.
$$
\end{thm}

\Pf The result follows from Theorem \ref{mpe}.
$\bl$

\begin{thm}
\label{s1-}
Let $\f$ be a nonconstant $1$-periodic function in
$L^2_{\rm{loc}}(\R)$. Then for $a\ge1$ 
$$
\|Q_{\f_a}^{[0,1]}\|_{\bS_1}\ge C(\f)\log(2+a).
$$
\end{thm}

\Pf It follows from \rf{012} that it is sufficient to prove that
\bay
\label{fia}
\|Q_{\f_a}^{[0,2]}\|_{\bS_1}\ge C(\f)\log(2+a).
\ey
First we consider the case where $a$ is an integer. There exists an
integer $k\in\Z\setminus\{0\}$ 
such that $\hat\f(k)\not=0$. By Corollary \ref{coef},  
$$
|\hat\f_a(l)|\le C(\log (2+|l|))^{-1}\|Q_{\f_a}\|_{\bS_1}.
$$ 
Substituting $l=ak$, we obtain \rf{fia},
since $\hat\f_a(ak)=\hat\f(k)$. 

Let now $a$ be an arbitrary
number in $(1,\infty)$.  For any $\sigma\in[1,2]$ there exists
$k_\sigma\in\Z\setminus\{0\}$ such that 
$\hat\f_\sigma(k_\sigma)\not=0$. Consequently, there exists a
neighborhood $U_\sigma$ of $\sigma$ 
such that $\hat\f_\tau(k_\sigma)\not=0$ for any $\tau\in
U_\sigma$. The first part of the proof 
allows us to obtain the required estimate for any $a>1$ such that
$a/N\in U_\sigma$ for some positive 
integer $N$. To complete the proof, we can choose a finite subcover
$U_{\sigma_j}$ of $[1,2]$. $\bl$  

\begin{thm}
\label{sp-}
Let $\f$ be a nonconstant $1$-periodic function in $L^2_{\rm{loc}}(\R)$
and let $0<p<1$. Then 
for $a\ge1$
$$
\left\|Q_{\f_a}^{[0,1]}\right\|_{\bS_p}\ge C(\f)a^{1/p-1}.
$$
\end{thm}

\Pf It suffices to consider the case when $a$ is an even integer.  The
general case may be reduced to 
this special case in the same way as in the proof of Theorem
\ref{s1-}. With any kernel $k$ on the 
square $[0,1)^2$ and any integer $n\ge1$ we associate the kernel
$k^{[n]}$ defined by 
$$
k^{[n]}(x,y)=n^2\int\limits_{\frac jn}^{\frac{j+1}n}
\int\limits_{\frac ln}^{\frac{l+1}n}
k(t,s)\,dtds,\quad \mbox{if}\quad
x\in\left[\frac jn,\frac{j+1}n\right)
\quad\mbox{and}\quad y\in\left[\frac ln,\frac{l+1}n\right).
$$
Clearly, $\|k^{[n]}\|_{\bS_p}\le\|k\|_{\bS_p}$ for any positive $p$
(recall that $\|k\|_{\bS_p}$ means the $S_p$-norm (or quasinorm) of the
integral operator with kernel $k$).

Suppose that 
$$
\int_0^1\f(t)\,dt\not=2\int_0^1t\f(t)\,dt.
$$ 
Put $k_n(x,y)\df\f\left(n\max\{x,y\}\right)$ for $x,y\in[0,1)$. Clearly,
$$
\left\|k_{2n}^{[n]}-k_{2n}^{[2n]}\right\|_{\bS_p}\le2^{1/p}\|k_{2n}\|_{\bS_p}.
$$
It is not hard 
to see that on $\left[\frac
jn,\frac{j+1}n\right)\times\left[\frac ln,\frac{l+1}n\right)$ 
$$
k_{2n}^{[n]}=\left\{\begin{array}{ll}\int\limits_0^1\f(t)\,dt,&j\ne l,\\[.2cm]
\frac12\int\limits_0^1\f(t)\,dt+\int\limits_0^1t\f(t)\,dt,&j=l.
\end{array}\right.
$$
Next, on $\left[\frac j{2n},\frac{j+1}{2n}\right)\times
\left[\frac l{2n},\frac{l+1}{2n}\right)$
$$
k_{2n}^{[2n]}=\left\{\begin{array}{ll}\int\limits_0^1\f(t)\,dt,&j\ne l,\\[.2cm]
2\int\limits_0^1t\f(t)\,dt,&j=l.
\end{array}\right.
$$
Thus, the kernel $k_{2n}^{[n]}-k_{2n}^{[2n]}$ vanishes outside the
``diagonal'' 
$\bigcup\limits_{j=0}^{n-1}\left[\frac jn,\frac{j+1}n\right)
\times\left[\frac jn,\frac{j+1}n\right)$. 
Clearly, for $x,\,y\in\left[0,\frac{n-1}n)\times[0,\frac{n-1}n\right)$
we have 
$$k_{2n}^{[n]}(x,y)-k_{2n}^{[2n]}(x,y)=k_{2n}^{[n]}
\left(x+\frac1n,y+\frac1n\right)-k_{2n}^{[2n]}
\left(x+\frac1n,y+\frac1n\right).
$$
Consequently,
$$
\left\|k_{2n}^{[n]}-k_{2n}^{[2n]}\right\|_{\bS_p}=
n^{1/p}\left\|\left(k_{2n}^{[n]}-k_{2n}^{[2n]}\right)\chi_{[0,\frac1n)
\times[0,\frac1n)}\right\|_{\bS_p}.
$$
We have 
$$
k_{2n}^{[n]}-k_{2n}^{[2n]}=\frac12\int_0^1\f(t)\,dt-\int_0^1t\f(t)\,dt
$$ on
$\left(\left[0,\frac1{2n}\right)\times\left[0,\frac1{2n}\right)\right)
\cup\left(\left[\frac1{2n},\frac1n\right)
\times\left[\frac1{2n},\frac1n\right)\right)$
and
$$
k_{2n}^{[n]}-k_{2n}^{[2n]}=\int_0^1t\f(t)\,dt-\frac12\int_0^1\f(t)\,dt
$$
on
$\left(\left[0,\frac1{2n}\right)\times\left[\frac1{2n},\frac1n\right)\right)
\cup\left(\left[\frac1{2n},\frac1n\right)
\times\left[0,\frac1{2n}\right)\right)$.
Now it is easy to see that 
$$
\left\|k_{2n}^{[n]}-k_{2n}^{[2n]}\right\|_{\bS_p}=cn^{1/p-1}
$$ 
for some nonzero $c$, since
$\int\limits_0^1\f(t)\,dt\not=2\int\limits_0^1t\f(t)\,dt$.  

Suppose now that $\f$ is an arbitrary nonconstant
$1$-periodic function. It suffices to prove that there exists $b\in\R$
such that $\int\limits_0^1 
\f(t-b)\,dt\not=2\int\limits_0^1t\f(t-b)\,dt$.  Suppose that 
\bay
\label{ft-b}
\int\limits_0^1\f(t-b)\,dt=2\int\limits_0^1t\f(t-b)\,dt,\quad b\in\R.
\ey
Let $h$ be the $1$-periodic function such that $h(t)=2t-1$ for $t\in[0,1)$.
Clearly, $\hat h(n)\ne0$ for $n\ne0$. Thus it follows from \rf{ft-b}
that $\f$ is constant. 
$\bl$

\

\section*{Acknowledgement}
\label{ack}
This research was begun by the first and third authors.
The two other authors joined them during the
Conference on Function Spaces, Interpolation Theory, and related
topics in honour of Jaak Peetre on his 65th
birthday 
at
Lund University, Sweden,
August 17--22, 2000,
where three of the authors were present and several of the results
were proved.
We thank Jonathan Arazy and Vladimir Maz'ya for helpful comments.

\

\end{document}